\documentclass[letterpaper,11pt]{amsart}%
\usepackage{amsmath,amsfonts,amssymb,amsthm, comment}
\usepackage{mathrsfs}
\usepackage{mathscinet}
\usepackage{graphicx}
\usepackage{enumerate}
\usepackage{algorithm}
\usepackage[noend]{algpseudocode}
\usepackage{color}
\usepackage[numbers]{natbib}
\usepackage{xspace}
\usepackage{amsmath}
\usepackage{amsfonts}
\usepackage{amssymb}
\usepackage{multirow}
\usepackage{pdflscape}
\usepackage{setspace}
\usepackage[margin=1in]{geometry}%
\providecommand{\U}[1]{\protect\rule{.1in}{.1in}}

\providecommand{\U}[1]{\protect\rule{.1in}{.1in}}

\providecommand{\U}[1]{\protect\rule{.1in}{.1in}}

\newtheorem{theorem}{Theorem}

\newtheorem{corollary}{Corollary}

\newtheorem{definition}{Definition}
\newtheorem{example}{Example}
\newtheorem{lemma}{Lemma}

\theoremstyle{remark}

\begin{document}
\title[Sample Out-of-Sample Inference]{Sample Out-of-Sample Inference Based on Wasserstein Distance}
\author[Blanchet and Kang]{{\ {\large J\MakeLowercase{ose} B\MakeLowercase{lanchet}}} \hspace{20pt}
{\large Y\MakeLowercase{ang} K\MakeLowercase{ang}}}
\address{ \\
{\large \MakeLowercase{{jose.blanchet@stanford.edu  \hspace{1cm}
yangkang@stat.columbia.edu} } }}

\begin{abstract}
We present a novel inference approach that we call Sample Out-of-Sample (or
SOS) inference. The approach can be used widely, ranging from semi-supervised
learning to stress testing, and it is fundamental in the application of
data-driven Distributionally Robust Optimization (DRO). Our method enables
measuring the impact of plausible out-of-sample scenarios in a given
performance measure of interest, such as a financial loss. The methodology is
inspired by Empirical Likelihood (EL), but we optimize the empirical
Wasserstein distance (instead of the empirical likelihood) induced by
observations. From a methodological standpoint, our analysis of the asymptotic
behavior of the induced Wasserstein-distance profile function shows dramatic
qualitative differences relative to EL. For instance, in contrast to EL, which
typically yields chi-squared weak convergence limits, our asymptotic
distributions are often not chi-squared. Also, the rates of convergence that
we obtain have some dependence on the dimension in a non-trivial way but remain controlled as the dimension increases.

\end{abstract}
\maketitle



\doublespacing

\textbf{Subject classifications}: Non-parametric statistics; Probability;
Distributionally Robust Optimization\newline\textbf{Area of review}:
Stochastic Models.\newline\bigskip\noindent\rule{16cm}{0.4pt} \clearpage

\section{Introduction}

The goal of this paper is to introduce a novel methodology for non-parametric
inference allows incorporating the adverse impact of out-of-sample
scenarios. We call the procedure Sample Out-of-Sample (SOS)) inference. Our method is general, and we discuss several applications,
including Distributionally Robust Optimization (DRO), semi-supervised learning,
and a novel stress-testing framework. We use the DRO framework in the
introduction to put our contributions in perspective. We elaborate on
semi-supervised learning and stress-testing applications in Section
\ref{Sect_Motivations}.

A data-driven DRO problem takes the form
\begin{equation}
\min_{\theta\in R^{d}}\max_{P\in\mathcal{U}_{\delta}\left(  P_{n}\right)
}E_{P}[\mathcal{L}\left(  \theta,X\right)  ], \label{General_DRO_F}%
\end{equation}
where $\mathcal{L}:\mathbb{R}^{d\times l}\rightarrow\lbrack0,\infty)$ is a cost (or
loss)\ function, $X\in\mathbb{R}^{l}$ is a random element, and $\theta
\in\mathbb{R}^{d}$ is a decision. Often, $\mathcal{L}\left(  \cdot,x\right)  $
is assumed to be strictly convex and smooth (e.g. twice differentiable) and we
will assume this throughout our motivating discussion. The notation
$E_{P}\left(  \cdot\right)  $ denotes the expectation operator associated to
the probability measure $P$. We use $P_{n}$ to denote the empirical measure
corresponding to $\left\{  X_{i}\right\}  _{i=1}^{n}$ independent identical distributed (i.i.d.) observations
that follow the distribution $P_{\ast}$. The set $\mathcal{U}_{\delta}\left(
P_{n}\right)  $ is the distributional uncertainty set. The parameter
$\delta>0$ is called the \textquotedblleft size of the distributional
uncertainty\textquotedblright\ so that the family of sets $(\mathcal{U}%
_{\delta}\left(  P_{n}\right)  :\delta\geq0)$ is increasing (in the sense of
inclusion) as a $\delta>0$ increases and so that for $\delta=0$,
$\mathcal{U}_{0}\left(  P_{n}\right)  =\{P_{n}\}$. Therefore, intuitively,
$P_{n}$ is the \textquotedblleft center\textquotedblright\ of the
distributional uncertainty region and $\delta>0$ can be thought of as its
\textquotedblleft radius.\textquotedblright

Ideally, one would like to compute $\theta_{\ast}=\arg\min E_{P_{\ast}%
}[\mathcal{L}\left(  \theta,X\right)  ]$, but $P_{\ast}$ is unknown.
Therefore, the intuition behind formulation (\ref{General_DRO_F}) is that one
is interested in choosing a decision $\theta$, which performs well uniformly
over a range of models that constitute reasonable (or plausible) variations
of the data (encoded by $P_{n}$).

We are interested in variations of the empirical distribution $P_{n}$ (the
elements in $\mathcal{U}_{\delta}\left(  P_{n}\right)  $) that systematically
explore the impact of \textit{out-of-sample scenarios} in the loss function
$\mathcal{L}\left(  \cdot\right)  $. Therefore, $P\in\mathcal{U}_{\delta
}\left(  P_{n}\right)  $ should not be supported only on the underlying data
set. Instead, we are interested in a framework that admits models in
$P\in\mathcal{U}_{\delta}\left(  P_{n}\right)$ that may be supported
outside the sample $\{X_{i}\}_{i=1}^{n}$. Because of this out-of-sample
exploration feature, we choose $\mathcal{U}_{\delta}\left(  P_{n}\right)  $
based on the Wasserstein distance of order 2, which is explained in Section
\ref{Sect_Basic_Defs}. We shall also discuss different alternative norms that
are supported by our analysis and discuss how these can be calibrated in a
data-driven way. 

Distributionally robust optimization formulations such as (\ref{General_DRO_F}%
)\ based on the Wasserstein distances have been studied recently in a wide
range of settings, especially in applications to machine learning and
artificial intelligence, see for example,
\cite{shafieezadeh-abadeh_distributionally_2015,MohajerinEsfahani2017,ZHAO2018262,blanchet_quantifying_2016,gao2016distributionally,blanchet2016robust,yang2017convex,sinha2018certifiable,gao2018robust,volpi2018generalizing,wolfram2018,blanchet2017data,blanchet2017doubly}%
.

All of these studies focus on the setting in which the support of the
distributions inside $\mathcal{U}_{\delta}\left(  P_{n}\right)  $ is
$\mathbb{R}^{d}$. Moreover, within the current literature, only
\cite{blanchet2016robust} studies the optimal selection of the parameter
$\delta$ by defining a natural optimization criterion. The work of
\cite{blanchet2016robust} also shows that such criterion recovers choices
that have been argued to be effective for recovery in machine learning
settings for which a DRO\ representation can be posed.

In contrast, compared to \cite{blanchet2016robust}, \textit{ our work is the
first one that studies the statistical implications of choosing the support of
the members of the distributional uncertainty }$P\in U_{\delta}\left(
P_{n}\right)  $\textit{ in a data-driven way.} One of our main contributions of this paper consists in providing a comprehensive study of an optimal
data-driven choice of uncertainty size, $\delta$, when the support of the
members in $\mathcal{U}_{\delta}\left(  P_{n}\right)  $ is obtained from an
arbitrary random sample whose size is increasing with $n$.

More generally, our contributions can be viewed in the lens of a novel
inference framework that we call SOS inference, based on the
analysis of the so-called SOS profile function for estimating equations.

In the DRO framework, we consider enriching the empirical data set
$\mathcal{X}_{n}=\{X_{i}\}_{i=1}^{n}$ (which is assumed to be i.i.d.)\ by
including a set of scenarios $\{Y_{i}\}_{i=1}^{m}$ (which is also assumed to
be i.i.d.), with $m=[\kappa n]$ for some $\kappa\in\lbrack0,\infty)$. The
$Y_{i}$'s and the $X_{i}$'s are not assumed share the same distribution. In
order to unify the notation we write $Z_{i}=X_{i}$ for $i=1,...,n$,
$Z_{n+k}=Y_{k}$ for $k=1,...,m$ and set $\mathcal{Z}_{n+m}=\{Z_{j}%
\}_{j=1}^{n+m}$. (We use $\mathbf{P}$ to denote the probability measure
supporting the infinite sequences $\{X_{i}\}_{i\geq1}$ and $\{Y_{i}\}_{i\geq
1}$, where the support of $\mathbf{P}$ is dense in the support of the
underlying sampling distribution.)

In order to emphasize the difference between the analysis in
\cite{blanchet2016robust} and our analysis here, we write $\mathcal{U}%
_{\delta}\left(  P_{n};\mathbb{R}^{l+1}\right)  $ to denote the full support
case (studied in \cite{blanchet2016robust}) and $\mathcal{U}_{\delta}\left(
P_{n};\mathcal{Z}_{n+m}\right)  $ for the uncertainty set considered in our
current setting.

Let us describe the optimality criterion introduced in
\cite{blanchet2016robust} for choosing $\delta$. Here we restrict the support
on the observed sequence and we would expect larger $\delta$ due to the extra
constraint. Since the set $\mathcal{U}_{\delta}\left(  P_{n}\right)  $ is
interpreted as the set of plausible variations of the data, then the set
\begin{equation}
\Lambda_{\delta}\left(  P_{n}\right)  =\{\theta:\theta=\arg\min E_{P}%
[\mathcal{L}\left(  \theta,X\right)  ]\text{ for }P\in\mathcal{U}_{\delta
}\left(  P_{n};\mathcal{Z}_{n+m}\right)  \} \label{CR_a}%
\end{equation}
corresponds to the set of plausible decisions, those that  are compatible with
the distributional uncertainty region. Note that $\Lambda_{\delta}\left(
P_{n}\right)  $ is a random set that can be interpreted as a confidence
region. The criterion that we utilize is the following
\begin{equation}
\min\{\delta:\mathbf{P}\left(  \theta_{\ast}\in\Lambda_{\delta}\left(
P_{n}\right)  \right)  \geq\alpha\}, \label{Crit_delta}%
\end{equation}
where $\alpha$ is a desired confidence level.

To analyze (\ref{Crit_delta}), we first argue that
\begin{equation}
\left\{  \theta_{\ast}\in\Lambda_{\delta}\left(  P_{n}\right)  \right\}
=\{R_{n}^{W}\left(  \theta_{\ast}\right)  \leq\delta\}, \label{Equiv_R}%
\end{equation}
for a suitable function, $R_{n}^{W}\left(  \cdot\right)  $, which we call the Sample-out-of-Sample (SOS) profile function. In simple words,
$R_{n}^{W}\left(  \theta_{\ast}\right)  $ can be computed directly in terms of
the shortest Wasserstein distance between $P_{n}$ and the set of probability
models $P\in\mathcal{U}_{\delta}\left(  P_{n};\mathcal{Z}_{n+m}\right)  $ for
which $E_{P}\left[  \nabla_{\theta}\mathcal{L}\left(  \theta_{\ast},x\right)
\right]  =0$.

As a consequence of (\ref{Equiv_R}), the optimal $\delta$ solving
(\ref{Crit_delta}) is simply the $\alpha$-quantile of $R_{n}^{W}\left(
\theta_{\ast}\right)  $.

In general, we can use our methodology to test the hypothesis that
$\theta_{\ast}$ satisfies $E_{P_{\ast}}\left(  h\left(  \theta_{\ast
},x\right)  \right)  =0$, simply replacing $\nabla_{\theta}\mathcal{L}\left(
\theta,x\right)  $ by $h\left(  \theta,x\right)  $ in the definition of the
SOS\ profile function. The hypothesis is rejected for high values of the
statistics $R_{n}^{W}\left(  \theta_{\ast}\right)  $. Thus, it is important to
compute the asymptotic distribution $R_{n}^{W}\left(  \theta_{\ast}\right)  $.

Our contributions are then stated at this level of generality (i.e., asymptotic
analysis of $R_{n}^{W}\left(  \theta_{\ast}\right)  $ for the purpose of
hypothesis testing). In the end, this paper involves two main methodological contributions:

\bigskip

\textbf{A)} First, we characterize the asymptotic distribution of $R_{n}%
^{W}\left(  \theta_{\ast}\right)  $ as $n\rightarrow\infty$; see Theorem
\ref{SOSTheoremMean}, Theorem \ref{SOSGeneralImplicit}, and Theorem \ref{SOSGeneralExplicit}. We explain how to compute the
asymptotic limiting distributions in Section
\ref{Section_Evaluating_Distribution}.

\bigskip

\textbf{B)} Second, we discuss various extensions that we believe are
natural to study in order to define DRO optimal transport cost functions.
These include implicit DRO\ formulations and plug-in estimators. We illustrate
the extensions in the empirical result section (Section \ref{Section_Applications}). For example, writing
$\theta_{\ast}=\left(  \gamma_{\ast},v_{\ast}\right)  $ we develop the
asymptotic distribution of $R_{n}^{W}\left(  \gamma_{\ast},\bar{v}_{n}\right)
$, where $\bar{v}_{n}$ is a suitable consistent plug-in estimator for
$v_{\ast}$ as $n\rightarrow\infty$ ; see Corollary
\ref{SOSGeneralExplicitPlugIn}. The construction of $\bar{v}_{n}$ may be based
on standard empirical estimators. This extension may be used in the context of
stochastic optimization with constraints, as illustrated in Section
\ref{Section_Applications}.

The theory that we develop in this paper parallels the main fundamental
results obtained in the context of Empirical Likelihood (EL), introduced by
Art Owen in \cite{owen_empirical_1988,owen_empirical_1990,owen_empirical_2001}%
. In fact, the construction of the function $R_{n}^{W}\left(  \cdot\right)$
borrows a great deal of inspiration from the empirical likelihood profile
function and its extensions based on divergence criteria, rather than the
likelihood function (see \cite{owen_empirical_2001}), and also see \cite{bayraksan2015data} for a comprehensive review of divergence-based distributional uncertainty sets in optimization, many of which are amenable to EL-based analysis. There are, however,
several important features of our framework that, we believe, add
significant value to the non-parametric inference literature. 

Before we discuss these features, we want to emphasize that our motivation is not to disprove the appropriateness of divergence approaches. The DRO community is actively investigating the advantages of various choices of uncertainty sets. Our discussion should be seen as a step in this direction. The most likely picture to eventually emerge is that divergence and Wasserstein approaches complement each other depending on issues such as convenience and tractability. For the purpose of using out-of-sample scenarios to inform the uncertainty set, we believe the Wasserstein distance is a natural choice, as we shall explain.

First, using divergence-based criteria (as it is typically done in standard EL settings) carries implicit support assumptions that seem unnatural in our setting as the sample size increases. For example, it is not difficult to see that a divergence-based distance between the empirical measure based on $n$ i.i.d. samples and that of $m=[\kappa n]$ i.i.d. samples (both from the same distribution) may not converge to zero. In our setting, this suggests that under divergence-type constructions, it requires a large uncertainty set to include distributions that one may reasonably and intuitively see as relatively small perturbation of the data.  So, choosing a large-sized uncertainty to accommodate these small perturbations may inflate the estimates artificially, just because the populations are large but unbalanced. Alternatively, if the size of uncertainty is small (which is expected under the null hypothesis as the sample size increases), the proportion of mass allocated outside the support of the empirical measure decreases to zero, so the overwhelming proportion of the mass in the models contained in the uncertainty set is concentrated in the support of the baseline model. Hence, we believe that the direct use of the EL framework may not be suitable in our setting. Additional out-of-sample issues that arise from using divergence criteria for data-driven distributional robust optimization
(closely related to EL) are noted in the stochastic optimization literature (see  \cite{esfahani_data-driven_2015}), and see also  \cite{wang_likelihood_2009,ben-tal_robust_2013} for related work.

Second, from a methodological standpoint, the mathematical techniques needed
to understand the asymptotic behavior of $R_{n}^{W}\left(  \theta_{\ast
}\right)  $ are qualitatively different from those arising typically in the
context of EL. We will show that if $l\geq3$, then the following weak
convergence limit holds (under suitable assumptions on $\mathcal{L}\left(
\cdot\right)  $),%
\[
n^{1/2+3/(2l+2)}R_{n}^{W}\left(  \theta_{\ast}\right)  \Rightarrow R\left(
\theta_{\ast}\right)  ,
\]
as $n\rightarrow\infty$. Note that the scaling depends on the dimension of the
random vector $X$ in a very particular way. In contrast, the Empirical
Likelihood Profile function is always scaled linearly in $n$ and the
asymptotic limiting distribution is generally a chi-squared distribution with
appropriate degrees of freedom and a constant scaling factor.

In our case, $R\left(  \theta_{\ast}\right)  $ can be explicitly
characterized, depending on the dimension in a non-trivial way, but it is no
longer a suitably scaled chi-squared distribution. When $l=1$, we obtain a
similar limiting distribution as in the EL case.
The intuition here is that a sample of order $O(n)$ provides enough coverage of the space since the optimal transport plan will displace points at distance $O(1/n^{1/2})$. The case $l=2$,
interestingly, requires a special analysis. In this case, the scaling remains
linear in $n$ (as in the case $l=1$), although the limiting distribution is
not exactly chi-squared, but a suitable quadratic form of a multivariate
Gaussian random vector. For the case $l\geq3$ the limiting distribution is not
a quadratic transformation of a multivariate Gaussian, but a more complex (yet
still explicit) polynomial function depending on the dimension.

At a high level, some of the qualitative distinctions in the methodology
arise because of the linear programming formulation underlying the SOS
function, which will typically lead to corner solutions (i.e., basic feasible
solutions in the language of linear programming). The high level intuition of
the scaling is associated with the interplay between the linear programming
formulation and the coverage of a sample of size $n$ in a space in $l$
dimension. A high-level intuition is given in more detail in Section
\ref{Sect_Intuitive}. In contrast to the analysis of the SOS function, in the
EL analysis of the profile function, the optimal solutions are amenable to a
smooth perturbation analysis as $n\rightarrow\infty$ using a Taylor expansion
of second (and higher) order terms. The lack of a continuously differentiable
derivative (of the optimal solution as a function of $\theta$) requires a
different type of analysis relative to the approach (traced back to the
classical Wilks' theorem as in \cite{wilks_large-sample_1938}), which lies at the
core of EL analysis.

The high-level intuition developed in Section \ref{Sect_Intuitive} also
underscores the distinction between our development here and the analysis in
\cite{blanchet2016robust}. In contrast to our development here, the scaling in
\cite{blanchet2016robust} is always dimension independent. This is because the
issue involving the coverage of the random scenarios in the support of the
alternative distributions is not a feature that needs to be considered.
Moreover, the current setting introduces a correlation structure in the
optimal transportation map, which is not present in the analysis of
\cite{blanchet2016robust}. This is because the feasible transportation
locations are now given by a random sample. To this end, we take advantage of
recent sample-path martingale inequalities. The use of these inequalities is
showcased in the technical Section \ref{Subsection_Technical_Lemmas} and we
believe that these techniques may be applicable more broadly in non-parametric
statistical analysis.

The rest of the paper is organized as follows. In Section
\ref{Sect_Motivations} we discuss semi-supervised learning and stress-testing
applications that motivate the formulation in which the support of
$P\in\mathcal{U}_{\delta}\left(  P_{n};\mathcal{Z}_{n+m}\right)  $ is
data-driven. Basic definitions, including a review of the Wasserstein
distance, are given in Section \ref{Sect_Basic_Defs}. Our main technical
results are described in Section \ref{Section_Main_Results}. We include
applications of our results to settings such as stochastic optimization, risk
analysis, and semi-supervised learning in Section \ref{Section_Applications}.
A short section including conclusions and additional discussions is given in
Section \ref{Section_Conclusions}. Finally, our technical development is given
in Section \ref{Section_Methodology}, starting with a high-level intuition of
the nature of our results and scaling in Section \ref{Sect_Intuitive}.

\section{Motivating Settings\label{Sect_Motivations}}

\subsection{Semi-supervised Learning Applications}

The setting of semi-supervised learning can be used to illustrate our
framework. Consider a classification problem that takes the form
$\mathcal{D}_{n}=\left\{  \left(  X_{i},Y_{i}\right)  \right\}  _{i=1}^{n}$
and $Y_{i}\in\{-1,1\}$ is the $i$-th response variable and $X_{i}\in R^{l}$ is
the $i$-th predictor. For concreteness, let us consider the logistic
regression setting in which
\[
P(Y_{i}=1|X_{i})=\frac{\exp\left(  Y_{i}\beta_{\ast}^{T}X_{i}\right)  }%
{1+\exp\left(  Y_{i}\beta_{\ast}^{T}X_{i}\right)  }=1-P(Y_{i}=-1|X_{i}).
\]
Suppose that we have access to an unlabeled data set $\left\{  X_{i}^{\prime
}\right\}  _{i=1}^{m}$ and we are interested in using this data in a
meaningful way for estimating $\beta_{\ast}$. This is the semi-supervised
learning setting arising in cases in which obtaining responses or labels
for every individual may be costly.
\begin{figure}[pth]
	\centering
	\includegraphics[width=15cm]{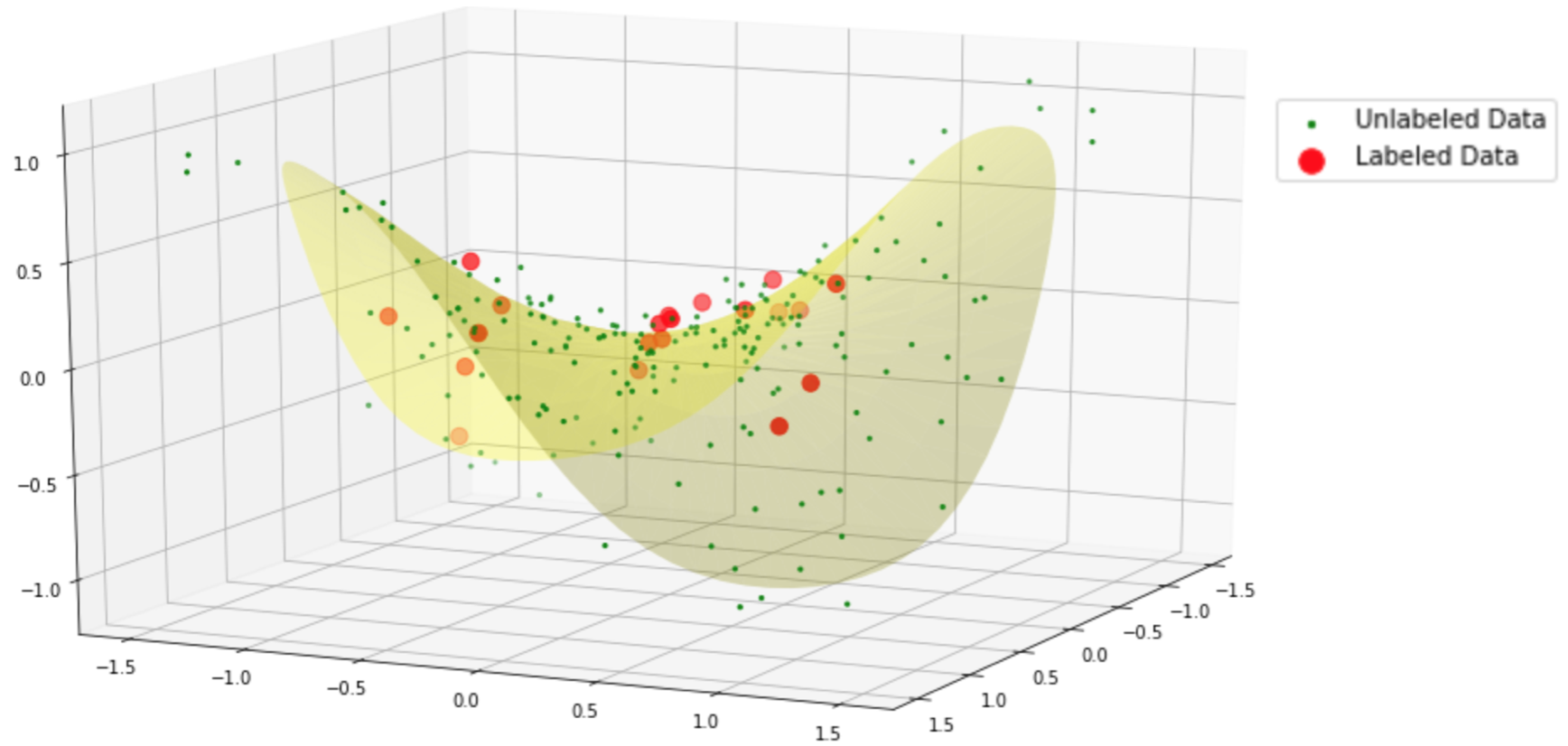}
	\caption{
	An illustrative example showing that the unlabeled observations (with green dots) can be used to provide a 
	proxy for the underlying manifold (the yellow surface) in which the
	the predictive variables lie; whereas the labeled data points (red dots) are not sufficient to provide such information. 
}%
	\label{Fig_DRO_intui1}%
\end{figure}

If the predictive variables are contained inside a lower-dimensional manifold embedded in the underlying ambient space, our intuition is that unlabeled data can be used as a proxy to profile precisely such a lower dimensional manifold. Thus it is natural to impose a DRO formulation that enhances statistical performance by quantifying the impact of out-of-sample scenarios that lay in the relevant lower-dimensional manifold.  This intuition is illustrated in Figure \ref{Fig_DRO_intui1}. 
The work of \cite{blanchet2017semi} proposes combining both the
labeled and unlabeled data by forming the set $\mathcal{X}_{n,m}=\mathcal{D}%
_{n}\cup\left(  \left\{  X_{i}^{\prime}\right\}  _{i=1}^{m}\times\left\{
-1, 1\right\}  _{i=1}^{n}\right)  $ (i.e., the original data set is enriched by
considering the unlabeled data with all the possible responses recorded by the
labeled data).

Then, \cite{blanchet2017semi} considers a DRO formulation for estimating $\beta_{\ast}$ in which the
distributional uncertainty region is defined in terms of the Wasserstein
distance. The DRO formulation proposed in \cite{blanchet2017semi}
is equivalent to the problem%
\begin{equation}
\min_{\beta}\max_{\mathcal{U}_{\delta}\left(  P_{n};\mathcal{X}_{n,m}\right)
}E_{P}\left[  \mathcal{L}\left(  X,Y,\beta\right)  \right]  ,
\label{SSL_1}%
\end{equation}
which corresponds to (\ref{General_DRO_F}).

The formulation of \cite{blanchet2017semi} (i.e. (\ref{SSL_1})) is
of significant interest because it is a natural semi-supervised learning
extension version of regularized linear regression, which is a highly popular
supervised machine learning estimator (see \cite{hastie_elements_2005}). In
particular, it is shown in \cite{blanchet2016robust}, see also
\cite{shafieezadeh-abadeh_distributionally_2015}, that replacing
$\mathcal{U}_{\delta}\left(  P_{n};\mathcal{X}_{n,m}\right)  $ by
$\ \mathcal{U}_{\delta}\left(  P_{n};\mathbb{R}^{l}\right)  $ in (\ref{SSL_1})
one recovers exactly regularized logistic regression and $\delta$ corresponds
exactly to the regularization parameter. This connection between Wasserstein
DRO and mainstream supervised learning estimators has been established for a
large class of methods, including square-root Lasso
(\cite{blanchet2016robust}), support vector machines
(\cite{blanchet2016robust}), group Lasso (\cite{blanchet2017distributionally2}%
), adaptive Lasso (\cite{blanchet2017data}), etc.

The methods developed in this paper provide the theoretical underpinning for
the choice of the uncertainty size $\delta$ in the context of (\ref{SSL_1}),
which yields regularized estimators that are informed by the unlabeled data in
a meaningful way.

\subsection{Novel Stress-testing Framework\label{Sect_SS}}

Consider the following stress-testing exercise. An insurance company wishes to
estimate a certain expectation of interest, say $\mathbb{E}_{\mathbb{P}^{\ast
}}(L(X))$, where $X$ might represent one or several risk factors, $L\left(
X\right)  $ is the corresponding financial loss and $\mathbb{P}^{\ast}\left(
\cdot\right)  $ is the underlying probability measure which may be unknown.

The insurance company may estimate $\mathbb{E}^{\ast}\left(  L(X)\right)  $
based on $n$ i.i.d. empirical
samples $X_{1},...,X_{n}\in\mathbb{R}^{l}$. However, the regulator (or
auditor) is also interested in quantifying the potential financial loss based
on stress scenarios, say an i.i.d. sample $Y_{1},...,Y_{m}\in\mathbb{R}^{l}$,
where $m=\left[  \kappa n\right]  $ with $\kappa\in\lbrack0,\infty)$. It may
be natural to choose $\kappa=1$ so that the amount of information provided by
the regulator and the company is balanced, but this is not necessary.

The scenarios provided by the regulator may or may not come from the same
distribution as the $X_{i}$'s. In fact, typically they will come from a
different distribution. The regulator's beliefs are captured by the
distribution of the $Y_{i}$'s. These beliefs may, in turn, be informed by the
knowledge that is accessible only by the regulator and not by the insurance
company. The regulator may not necessarily question the fact that the
historical data from the $X_{i}$'s follows distribution $\mathbb{P}_{\ast
}\left(  \cdot\right)  $, but the regulator might be concerned that the insurance
company lacks additional information to assess the overall risk exposure better.

On the one hand, the insurance company clearly knows well its idiosyncratic
risk exposures, so the data represented by the $X_{n}$'s, arising from a model
with such idiosyncratic information is meaningful and should be considered
carefully. On the other hand, it is also correct that the regulator possesses additional information that should be considered in evaluating the potential impact of scenarios that may not be appropriately captured by the data of the insurance company.

How does one incorporate both the $X_{i}$'s and the $Y_{i}$'s in a meaningful
way for the purposes of evaluating the risk of the company?

The methodology developed in this paper allows incorporating both the empirical data of the insurance company and the stress scenarios provided by the regulator into a Distributionally Robust Performance Analysis (DRPA) formulation (closely related to Distributionally Robust Optimization -- DRO) as we describe next.

Define $Z_{k}=X_{k}$ for $k=1,\ldots,n$ and $Z_{n+k}=Y_{k}$ for $k=1,\ldots,m$
(i.e., merge both the empirical samples and the stress scenarios into a set
$\mathcal{Z}_{n+m}=\{Z_{1},...,Z_{n+m}\}$). We let
\[
P_{n}\left(  dx\right)  =n^{-1}\sum_{k=1}^{n}\delta_{\{X_{k}\}}\left(
dx\right)
\]
be the empirical distribution of the data generated by the insurance company.
A natural estimate for $\mathbb{E}^{\ast}\left(  L(X)\right)  $ based on the
insurance company's data is given by
\[
\mathbb{E}_{P_{n}}\left(  L(X)\right)  =n^{-1}\sum_{k=1}^{n}L\left(
X_{k}\right)  .
\]
Now, let $\mathcal{P}\left(  \mathcal{Z}_{n+m}\right)  $ be the set of all
probability distributions with support on $\mathcal{Z}_{n+m}$. Our DRPA
approach consists in providing estimates for $\mathbb{E}_{P}\left(
L(X)\right)  $ via%
\begin{equation}
\theta_{-}\left(  \delta\right)  ,\theta_{+}\left(  \delta\right)
=\underset{P\in\mathcal{U}_{\delta}\left(  P_{n};\mathcal{Z}_{n,m}\right)
}{\min\; ,\;\max}\mathbb{E}_{P}\left(  L(X)\right)  . \label{DRO_00}%
\end{equation}

We believe that the DRPA formulation (\ref{DRO_00}) provides a reasonable
approach for combining both the insurance company's information and the
regulator's beliefs. We do not disregard the data coming from the insurance
company (in fact, the empirical distribution $P_{n}$ is placed at the center
of the uncertainty set), but we also capture the potential impact of
out-of-sample scenarios based on the regulator's beliefs.

Formulation (\ref{DRO_00}) is closely related to (\ref{General_DRO_F}) and the
methodology that we present in this paper can be used to find an optimal
choice for $\delta$. In particular, an equivalent way of representing the
range $[\theta_{-}\left(  \delta\right)  ,\theta_{+}\left(  \delta\right)  ]$
is in terms of a suitably defined SOS\ profile function\ (or \textquotedblleft
SOS function\textquotedblright), $R_{n}^{W}\left(  \cdot\right)  $, as we
shall see, so that%
\begin{equation}
\lbrack\min\{\theta:R_{n}^{W}\left(  \theta\right)  \leq\delta\},\max
\{\theta:R_{n}^{W}\left(  \theta\right)  \leq\delta\}]=[\theta_{-}\left(
\delta\right)  ,\theta_{+}\left(  \delta\right)  ]. \label{CI_01}%
\end{equation}
Therefore, the study of the function $R_{n}^{W}\left(  \cdot\right)  $ is a
key in the analysis of (\ref{DRO_00}) and the selection of $\delta$ based on
statistical principles, and this leads us to our contributions A)-B) described
in the Introduction.

We emphasize, however, that our choice of $\delta$ is purely statistical. That is, we operate under the blanket assumption that the risk is correctly computed solely with the bank's internal data as the sample size grows to infinity. Under this assumption there is less and less need for scenarios as the sample size of the internal data increases. In practice, the sample size is always finite and, in the end, the choice of regulatory capital is the result of an informed negotiation between the regulator and the bank. We provide a tool that helps to inform this discussion because it statistically combines both elements (internal data and external scenarios) in a way that is consistent with the guidelines described in \cite{fed_stress_2019} for generating stress scenarios. However, non-statistical criteria (e.g., social cost based) may also be used to choose $\delta$, leading to, for instance, hybrid methods that would build on our current development. However, these types of hybrid choices would require additional modeling elements that are beyond the scope of our statistical treatment.

\section{Basic Definitions\label{Sect_Basic_Defs}}

Throughout our development we adopt the convention that all vectors we
consider are expressed as columns, so, for example, $x^{T}=\left(
x_{1},...,x_{l}\right)  $ is a row vector in $\mathbb{R}^{l}$ (here we use
$x^{T}$ to denote the transpose of $x$). Also, given a random variable
$W\in\mathbb{R}^{d}$ so that $\mathbb{E}\left(  W\right)  =0$ and
$\mathbb{E}\left(  \left\Vert W\right\Vert _{2}^{2}\right)  <\infty$, we use
$Var\left(  W\right)  =\mathbb{E}\left(  WW^{T}\right)  $ to denote the
covariance matrix of $W$.

\subsection{On Wasserstein Distance and Distributional Uncertainty}

As we mentioned in the introduction, we utilize the Wasserstein distance of
order 2 to describe the distributional uncertainty region. We consider two
closed subsets of $\mathbb{R}^{l}$, namely $\mathcal{X}$ and $\mathcal{Z}$. We
use the notation $\mathcal{P}\left(  \mathcal{X}\times\mathcal{Z}\right)  $ to
denote all the Borel probability measures $\pi$ with support on $\mathcal{X}%
\times\mathcal{Z}$. Any $\pi\in\mathcal{P}\left(  \mathcal{X}\times
\mathcal{Z}\right)  $ can be thought of as the joint distribution of a pair of
random vectors $\left(  X,Z\right)  $. We use the notation $\pi_{X}$ to denote
the marginal distribution of $X$ under $\pi$; similarly, $\pi_{Z}$ is the
marginal distribution of $Z$ under $\pi$.

The Wasserstein distance (of order $2$) between the Borel probability measures
$\mu$ and $\upsilon$, supported on $\mathcal{X}$ and $\mathcal{Z}$,
respectively, is defined as $\sqrt{D\left(  \mu,\upsilon\right)  }$, where%
\[
D\left(  \mu,\upsilon\right)  =\inf\{\int\int\left\Vert x-z\right\Vert
_{2}^{2}\pi\left(  dx,dz\right)  :\pi\in\mathcal{P}\left(  \mathcal{X}%
\times\mathcal{Z}\right)  ,\pi_{X}=\mu,\pi_{Z}=v\}.
\]
In simple words, the square of the Wasserstein distance of order 2 (under the
Euclidean metric) is defined as the minimum cost of transporting the mass
encoded by $\mu$ into the mass encoded by $\upsilon$; computing the
unitary-cost-per-transportation of a unit of mass from $x$ to $y$ as the
square of the Euclidean distance between the source ($x$) and destination ($y$).

Our results can be directly adapted to the situation in which the Euclidean
metric is replaced by the so-called Mahalanobis distance, namely, $\left\Vert
x-y\right\Vert _{A}^{2}=\left(  x-y\right)  ^{T}A\left(  x-y\right)  $ for any
positive definite matrix $A$. The use of this distance and procedures to fit
$A$ for classification tasks based on manifold learning tools are studied in
\cite{blanchet2017data}. In order to simplify the notation and the exposition we continue
with the standard Euclidean metric throughout our development, corresponding
to $A=I$.

In the sequel, $\mathcal{X}$ and $\mathcal{Z}$ are finite cardinality sets.
Therefore, in this case, the evaluation of $D\left(  \mu,\upsilon\right)  $ is
a finite dimensional linear programming problem and so, conceptually,
computing $D\left(  \mu,\upsilon\right)  $ is straightforward. The Wasserstein
distance is defined in great generality (for arbitrary metric spaces) as the
solution of the Monge-Kantorovich problem with the cost-per-transportation
defined in terms of the underlying metric. We refer the reader to
\cite{villani_optimal_2008} for more information on Wasserstein distances.
Because we focus on the finite-cardinality case,  it is enough with elementary
notions of finite dimensional linear programming to understand the definition
we use in this paper.

The distributional uncertainty set, $\mathcal{U}_{\delta}\left(  P_{n}\right)
$, mentioned in the Introduction to motivate our contributions can then be
defined by choosing $\mathcal{X=X}_{n}$ and $\mathcal{Z}=\mathcal{Z}_{n+m}$
and letting
\[
\mathcal{U}_{\delta}\left(  P_{n}\right)  =\mathcal{U}_{\delta}\left(
P_{n};\mathcal{Z}_{n+m}\right)  =\{P:D\left(  P_{n},P\right)  \leq\delta\}.
\]

\subsection{The SOS Profile Function}

To motivate the definition of the SOS\ Profile function, once again, we return
to the DRO\ framework defined in the Introduction. We note from (\ref{CR_a})
that%
\[
\Lambda_{\delta}\left(  P_{n}\right)  =\{\theta:E_{P}\left[  h\left(
X,\theta\right)  \right]  =0\text{ for }P\in\mathcal{U}_{\delta}\left(
P_{n}\right)  \},
\]
where $h\left(  X,\theta\right)  =\nabla_{\theta}\mathcal{L}\left(
X,\theta\right)  $. So (by convexity) we have that $\theta_{\ast}\in
\Lambda_{\delta}\left(  P_{n}\right)  $ if and only if there exists
$P\in\mathcal{U}_{\delta}\left(  P_{n}\right)  $ such that
\begin{equation}
E_{P}\left[  h\left(  X,\theta_{\ast}\right)  \right]  =0. \label{Eq_Opt}%
\end{equation}

Let $R_{n}^{W}\left(  \theta_{\ast}\right)  $ be the smallest transportation
cost (measured by $D\left(  P_{n},P\right)  $) between $P_{n}$ and any member
$P\in\mathcal{P}\left(  \mathcal{Z}_{n+m}\right)  $ for which (\ref{Eq_Opt}) is true.
It is easy to reason that $R_{n}^{W}\left(  \theta_{\ast}\right)  \leq\delta$
if and only if $\theta_{\ast}\in\Lambda_{\delta}\left(  P_{n}\right)  $.
Formally, we have the following definition for the SOS profile function
$R_{n}^{W}\left(  \theta\right)  $, namely%
\begin{equation}
R_{n}^{W}\left(  \theta\right)  =\min\{D\left(  P_{n},P\right)  :E_{P}\left[
h\left(  X,\theta\right)  \right]  =0\}. \label{Def_SOS_Profile}%
\end{equation}

The goal of this paper is to study the behavior of $R_{n}^{W}\left(
\theta_{\ast}\right)  $ under the estimating equation assumption
\begin{equation}
E_{P_{\ast}}\left[  h\left(  X,\theta_{\ast}\right)  \right]  =0,
\label{Est_Eq_A}%
\end{equation}
and the $\left\{  X_{i}\right\}  _{i=1}^{n}$ being an i.i.d. sample from
$P_{\ast}$. We will formulate our results in terms of the estimating equation
(\ref{Est_Eq_A}) for general $h\left(  \cdot\right)  $ (not necessarily
arising from an optimization problem).

We consider this more general framework because we believe that our results
may be applicable to inference settings other than DRO, for instance, the
stress-testing framework described earlier. In fact, we now return to such
setting to explain how to use the SOS\ profile function in this case.

\subsubsection{The SOS Profile function for stress-testing setting}

In the stress-testing setting described earlier, we wish to select $\delta$
just as large to guarantee that $\theta_{\ast}:=\mathbb{E}_{\mathbb{P}^{\ast}%
}(L(X))\in\lbrack\theta_{-}\left(  \delta\right)  ,\theta_{+}\left(
\delta\right)  ]$ with a certain degree of confidence, which we shall denote
by $\alpha$.

Therefore, because of equation (\ref{CI_01}), we are interested in choosing
the smallest $\delta$ so that
\begin{equation}
P\{\theta_{\ast}\in\lbrack\theta_{-}\left(  \delta\right)  ,\theta_{+}\left(
\delta\right)  ]\}=\mathbf{P}\{R_{n}^{W}\left(  \theta_{\ast}\right)
\leq\delta\}=\alpha. \label{EQ_CONF_INT}%
\end{equation}
In other words, $\delta$ is chosen to be the $\alpha$-quantile of the random
variable
\[
R_{n}^{W}\left(  \theta_{\ast}\right)  =\min\{D\left(  P_{n},P\right)
:E_{P}\left[  L\left(  X\right)  -\theta_{\ast}\right]  =0\}.
\]
Note that this formulation is a particular case of the one introduced in
(\ref{Est_Eq_A}) by letting $h\left(  \theta,x\right)  =L\left(  x\right)
-\theta$. For pedagogical reasons, we will present our results first for the
SOS profile function for means  (i.e., assuming that $L\left(  x\right)
=x$) and later we move to more general estimating equations.

\section{Main Results\label{Section_Main_Results}}

\subsection{SOS Function for Means}

We state the following underlying assumptions throughout this subsection.

\textbf{A1):} Let us write $\mathcal{X}_{n}=\{X_{1},...,X_{n}\}\subset$
$\mathbb{R}^{l}$ to denote an i.i.d. sample from a continuous distribution.
Therefore, the cardinality of the set $\mathcal{X}_{n}$ is $n$.

\textbf{A2): }We also consider an independent i.i.d. sample $\mathcal{Y}%
_{m}=\{Y_{1},...,Y_{m}\}\subset$ $\mathbb{R}^{l}$ from a continuous
distribution. Throughout our discussion we shall assume that $m=\left[  \kappa
n\right]  $ with $\kappa\in\lbrack0,\infty)$.

\textbf{A3): }Assume that $\mathbb{E}\left\Vert X_{1}\right\Vert _{2}%
^{2}+\mathbb{E}\left\Vert Y_{1}\right\Vert _{2}^{2}<\infty$.

\textbf{A4): }If $l=1$ we assume that $X_{i}$ and $Y_{i}$ have positive
densities $f_{X}\left(  \cdot\right)  $ and $f_{Y}\left(  \cdot\right)  $. If
$l\geq2$ we assume that $X_{i}$ and $Y_{i}$ have differentiable positive
densities $f_{X}\left(  \cdot\right)  $ and $f_{Y}\left(  \cdot\right)  $,
with bounded gradients.

Define $\mathcal{Z}_{n+m}=\{Z_{1},...,Z_{n+m}\}=\mathcal{X}_{n}\cup
\mathcal{Y}_{m}$, with $Z_{k}=X_{k}$ for $k=1,...,n$, and $Z_{n+j}=Y_{j}$ for
$j=1,...,m$. For any closed set $\mathcal{C}$ let us write $\mathcal{P}\left(
\mathcal{C}\right)  $ to denote the set of probability measures supported on
$\mathcal{C}$. Therefore, in particular, a typical element $\upsilon_{n}%
\in\mathcal{P}\left(  \mathcal{Z}_{n+m}\right)  $ takes the form
\[
\upsilon_{n}\left(  dz\right)  =\sum_{k=1}^{n+m}v\left(  k\right)
\delta_{Z_{k}}\left(  dz\right)  ,
\]
where $\delta_{Z_{k}}\left(  dz\right)  $ is a Dirac measure centered at
$Z_{k}$. Now, we shall use $\mu_{n}\in\mathcal{P}\left(  \mathcal{X}%
_{n}\right)  $ to denote the empirical measure associated to $\mathcal{X}_{n}%
$, that is,
\[
\mu_{n}\left(  dx\right)  =\frac{1}{n}\sum_{i=1}^{n}\delta_{X_{i}}\left(
dx\right)  .
\]

Given any $\pi\in\mathcal{P}\left(  \mathcal{X}_{n}\times\mathcal{Z}%
_{n+m}\right)  $ we write $\pi_{X}\in\mathcal{P}\left(  \mathcal{X}%
_{n}\right)  $ to denote the marginal distribution with respect to the first
coordinate, namely $\pi_{X}\left(  dx\right)  =\int_{z\in\mathcal{Z}_{n+m}}%
\pi\left(  dx,dz\right)  $ and, likewise, we define $\pi_{Z}\in\mathcal{P}%
\left(  \mathcal{Z}_{n}\right)  $ as $\pi_{Z}\left(  dz\right)  =\int%
_{x\in\mathcal{X}_{n}}\pi\left(  dx,dz\right)  $.

We have the following formal definition of the SOS function for estimating means.

\begin{definition}
The SOS function, $R_{n}^{W}\left(  \cdot\right)  $, to estimate $\theta
_{\ast}=E\left(  X\right)  $ is defined as%
\begin{align}
R_{n}^{W}\left(  \theta_{\ast}\right)   &  =\inf\{\int\int\left\Vert
x-z\right\Vert _{2}^{2}\pi\left(  dx,dz\right)  :\label{EWP_1}\\
&  \text{s.t. }\left.  \pi\in\mathcal{P}\left(  \mathcal{X}_{n}\times
\mathcal{Z}_{n+m}\right)  ,\pi_{X}=\mu_{n},\pi_{Z}=v_{n},\int zv_{n}\left(
dz\right)  =\theta_{\ast}\right.  \},\nonumber\\
&  =\inf\{\int\int\left\Vert x-z\right\Vert _{2}^{2}\pi\left(  dx,dz\right)
:\nonumber\\
&  \text{s.t. }\left.  \pi\in\mathcal{P}\left(  \mathcal{X}_{n}\times
\mathcal{Z}_{n+m}\right)  ,\pi_{X}=\mu_{n},\int z\pi_{Z}\left(  dz\right)
=\theta_{\ast}\}\right.  .\nonumber
\end{align}
(Here and throughout the paper, s.t. is an abbreviation for \textquotedblleft subject
to.\textquotedblright)
\end{definition}

We now state the following asymptotic distributional result for the SOS function.

\begin{theorem}
[SOS Profile Function Analysis for Means]\label{SOSTheoremMean} In addition to
Assumptions A1)-A4), suppose that the covariance matrix of $X$, $Var\left(
X\right)  $, exists. The following asymptotic result follows

\begin{itemize}
\item When $l=1$,
\[
nR_{n}^{W}(\theta_{\ast})\Rightarrow\sigma^{2}\chi_{1}^{2}%
\]
where $\sigma^{2}=Var\left(  X\right)  $.

\item When $l=2$, define $\tilde{Z}\sim N\left(  0,Var\left(  X\right)
\right)  \in\mathbb{R}^{l}$, then
\[
nR_{n}^{W}\left(  \theta_{\ast}\right)  \Rightarrow\rho\left(  \tilde
{Z}\right)  \left(  2-\tilde{\eta}\left(  \tilde{Z}\right)  \rho\left(
\tilde{Z}\right)  \right)  \left\Vert \tilde{Z}\right\Vert _{2}^{2}%
\]
where $\rho:=\rho\left(  \tilde{Z}\right)  $ is the unique solution to the
equation
\[
\frac{1}{\rho}=\tilde{g}\left(  \rho\tilde{Z}\right)  ,
\]
and $\tilde{g}:\mathbb{R}^{l}\rightarrow\mathbb{R}$ is a deterministic
function defined as
\[
\tilde{g}\left(  x\right)  =\mathbb{P}\left( \tau(0)\leq\left\Vert x\right\Vert
_{2}^{2}\right)  ,
\]
where $\tau$ is a random variable satisfying
\[
\mathbb{P}\left(  \tau>t\right)  =\mathbb{E}\left[  \exp\left(  -\left(
f_{X}\left(  X_{1}\right)  +\kappa f_{Y}\left(  X_{1}\right)  \right)  {\pi
}t\right)  \right]  .
\]
And the function $\tilde{\eta}:\mathbb{R}^{l}\rightarrow\mathbb{R}$ is a
deterministic function given as
\[
\tilde{\eta}\left(  x\right)  =\mathbb{E}\left[  \max\left(  1-{\tau(0)
}/{\left\Vert x\right\Vert _{2}^{2}},0\right)  \right]  .
\]

\item When $l\geq3$,
\[
n^{1/2+\frac{3}{2l+2}}R_{n}^{W}\left(  \theta_{\ast}\right)  \Rightarrow
\frac{2l+2}{l+2}\frac{\left\Vert \tilde{Z}\right\Vert _{2}^{1+\frac{1}{l+1}}%
}{\left(  \mathbb{E}\left[  \frac{\pi^{l/2}}{\Gamma(l/2+1)}\left(
f_{X}\left(  X_{1}\right)  +\kappa f_{Y}\left(  X_{1}\right)  \right)
\right]  \right)  ^{\frac{1}{l+1}}}%
\]
where $\tilde{Z}\sim N\left(  0,Var\left(  X\right)  \right)  \in
\mathbb{R}^{l}$.
\end{itemize}
\end{theorem}

\subsubsection{More on the limiting distribution}
The limiting distributions that we obtain are explicitly
characterized. They depend on parameters that are meaningful in the application settings that we shall discuss. For example, the distribution from which stress scenarios are generated or the distribution of the predictors of the unlabeled data clearly play a key role in the limiting distribution. These parameters dictate the \textquotedblleft spread\textquotedblright of the distribution and, consequently, the size of
quantiles. So, the parameters that appear in our limit theorems readily
affect the uncertainty size in a quantifiable way.

In order to make this point relatively more tangible, consider the following example based on simulated data.

\begin{figure}[pth]
	\centering
	\includegraphics[width=10cm]{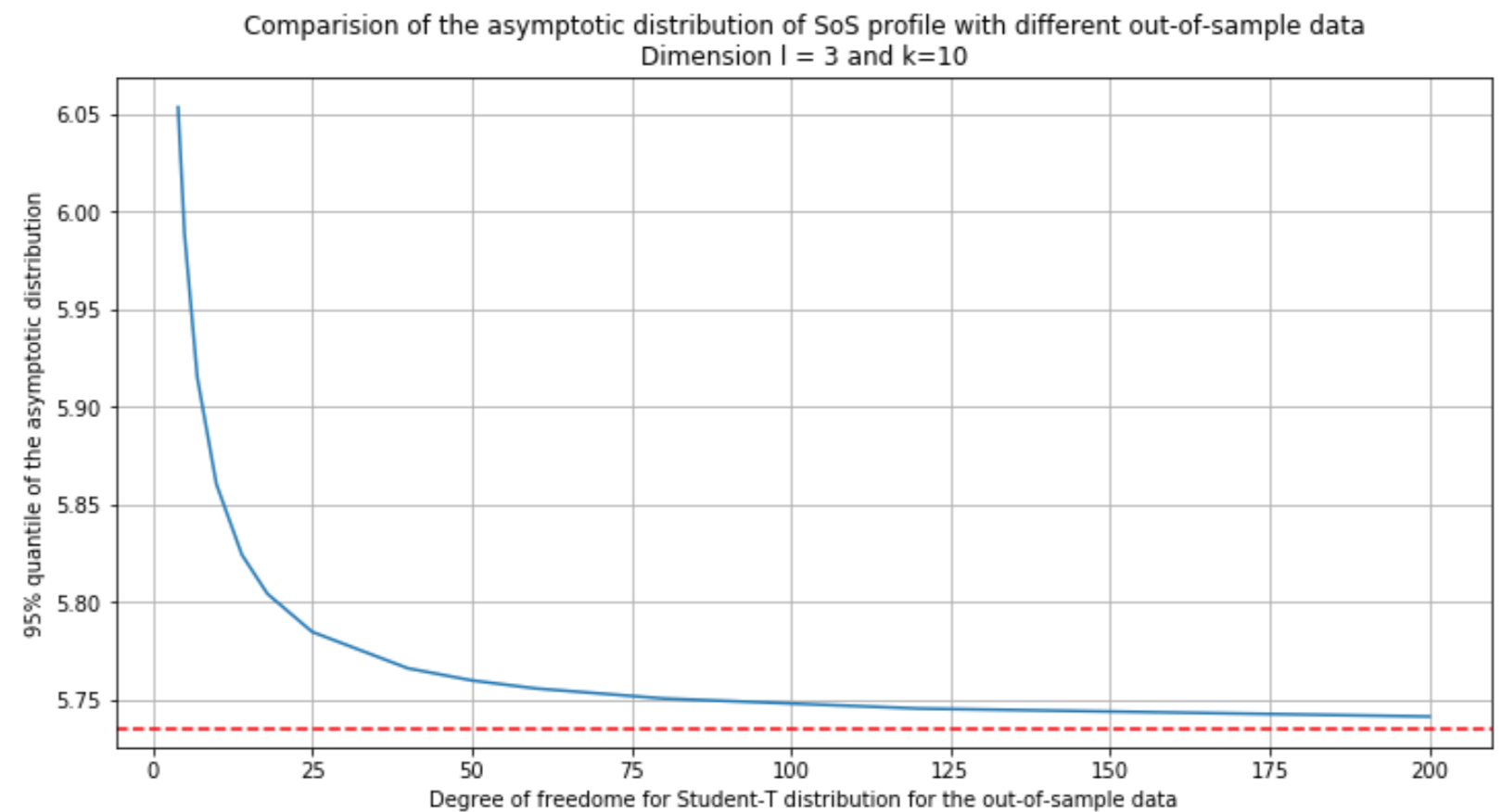}
	\caption{
	$95\%$ quantile for the SOS profile function asymptotic distribution (dimension being 3 and the $\kappa=10$) with different degree of freedom for the Student-t distribution in stress scenarios. The in-sample data is standard Gaussian. The red dashed line illustrates the situation in which stress scenarios are also chosen to be standard Gaussian. }
	\label{Fig_quantile_df}%
\end{figure}

Our asymptotic theorem gives different asymptotic distributions for different degrees of freedom (d.f.) in the Student-t distribution. If we select the $95\%$ quantile for the construction of our robust risk valuation interval, we can see that the higher the d.f., the smaller the quantile, as we show in Figure \ref{Fig_quantile_df}. So, an increase in the likelihood of more extreme scenarios provided by the regulator translates directly into a larger confidence region for the risk or a larger size in the uncertainty region, in a precisely quantifiable way thanks to our results.
The SOS profile function is the distance between the empirical distribution and the manifold determined by the estimating equation(s). If the in-sample data and the stress-scenario data are more similar, we would expect smaller quantiles (this corresponds to the setting in which the d.f. is large for the Student-t distribution), and we will observe larger quantiles when the two distributions are different from each other (this is the setting in which the d.f. is small for Student-t).

\subsubsection{Evaluating the Limiting Distribution}

\label{Section_Evaluating_Distribution}

In Theorem \ref{SOSTheoremMean} and in the rest of our results, the limiting
distribution depends on parameters that might be unknown. For example, take
the case $l\geq3$ in Theorem \ref{SOSTheoremMean}. We obtain that
\begin{equation}
n^{1/2+\frac{3}{2l+2}}R_{n}^{W}\left(  \theta_{\ast}\right)  \Rightarrow
\frac{2l+2}{l+2}\frac{\left\Vert \tilde{Z}\right\Vert _{2}^{1+\frac{1}{l+1}}%
}{\left(  c_{0}\right)  ^{1/\left(  l+1\right)  }}, \label{RHS_A}%
\end{equation}
where
\[
c_{0}:=\mathbb{E}\left[  \frac{\pi^{l/2}}{\Gamma(l/2+1)}\left(  f_{X}\left(
X_{1}\right)  +\kappa f_{Y}\left(  X_{1}\right)  \right)  \right]
\]
and $\tilde{Z}\sim N\left(  0,Var\left(  X\right)  \right)  $. This situation
is quite standard when developing asymptotic distributions for hypothesis
testing and the remedy is to simply use any consistent plug-in estimator to
estimate the unknown quantities. For instance, we can use
\[
\Sigma_{n}=\frac{1}{n}\sum_{j=1}^{n}\left(  X_{j}-E_{P_{n}}\left(
X_{j}\right)  \right)  \left(  X_{j}-E_{P_{n}}\left(  X_{j}\right)  \right)
^{T}%
\]
instead of $\Sigma=Var\left(  X\right)  $. We can also use any consistent
estimator (converging on compact sets and with rapid decay at infinity) for
the densities of $f_{X}\left(  \cdot\right)  $ and $f_{Y}\left(  \cdot\right)
$, say $f_{X}^{\left(  n\right)  }\left(  \cdot\right)  $ and $f_{Y}^{\left(
n\right)  }\left(  \cdot\right)  $, respectively, and estimate $c_{0}$ via%
\[
c_{0}\left(  n\right)  =E_{P_{n}}\left[  \frac{\pi^{l/2}}{\Gamma
(l/2+1)}\left(  f_{X}^{\left(  n\right)  }\left(  X\right)  +\kappa
f_{Y}^{\left(  n\right)  }\left(  X_{1}\right)  \right)  \right]  ,
\]
which is consistent as $n\rightarrow\infty$. Because the asymptotic
distribution in (\ref{RHS_A}) is continuous in $c_{0}$ and $\Sigma$, it
follows that estimating quantiles based on the plug-in estimators
$c_{0}\left(  n\right)  $ and $\Sigma_{n}$ in place of $c_{0}$ and $\Sigma$
leads to asymptotically equivalent specifications for the asymptotic quantiles
of $R_{n}^{W}\left(  \theta_{\ast}\right) $. These quantiles, in turn, can be
estimated by Monte Carlo using the asymptotic limits, with the plug-in
estimators in place. A completely analogous approach can be followed for the
asymptotic distributions obtained in the developments that we discuss next.

\subsection{SOS Function for Estimating Equations}

Throughout this subsection we assume that \textbf{A1)} and \textbf{A2)} are in
force. Let us assume that $h:R^{d}\times R^{l}\rightarrow R^{q}$ and $q\leq d$. We also impose the following assumptions.

\textbf{B1)} Assume $\theta_{\ast}\in R^{d}$ satisfies
\[
\mathbb{E}\left(  h\left(  \theta_{\ast},X\right)  \right)  =0.
\]

\textbf{B2)} Furthermore, suppose that
\[
\mathbb{E}\left\Vert h\left(  \theta_{\ast},X\right)  \right\Vert _{2}%
^{2}<\infty, \text{ and }\mathbb{E}\left\Vert h\left(  \theta_{\ast},Y\right)  \right\Vert _{2}%
^{2}<\infty.
\]

Our goal is to estimate $\theta_{\ast}$ under two reasonable SOS function
formulations, which we shall discuss. These are \textquotedblleft
implicit\textquotedblright\ or \textquotedblleft indirect\textquotedblright%
\ and \textquotedblleft explicit\textquotedblright\ or \textquotedblleft
direct\textquotedblright\ formulations, we will explain their nature next.

\subsubsection{Implicit SOS Formulation for Estimating Equations}

The first SOS function form for estimating equations is the following; we call
it Implicit SOS or Indirect SOS function because the Wasserstein distance is
applied to $h\left(  \theta,X_{i}\right)  $ and $h\left(  \theta,Z_{k}\right)
$ and thus it implicitly or indirectly induces a notion of proximity among the samples.

\begin{definition}
[Implicit SOS Profile Function for Estimating Equations]%
\label{ImplicitSOSFormulation} Let us write $\mathcal{X}_{n}^{h}\left(
\theta_{\ast}\right)  =\{h\left(  \theta_{\ast},X_{i}\right)  :X_{i}%
\in\mathcal{X}_{n}\}$ and $\mathcal{Z}_{n+m}^{h}\left(  \theta_{\ast}\right)
=\{h\left(  \theta_{\ast},Z_{k}\right)  :Z_{k}\in\mathcal{Z}_{n+m}\}$ then
\begin{align}
R_{n}^{W}(\theta_{\ast})  &  =\inf\{\int\int\left\Vert h\left(  \theta_{\ast
},x\right)  -h\left(  \theta_{\ast},z\right)  \right\Vert _{2}^{2}\pi\left(
dx,dz\right)  :\label{WassersteinProfileFunctionGeneralized}\\
&  \text{s.t. }\left.  \pi\in\mathcal{P}\left(  \mathcal{X}_{n}^{h}\left(
\theta_{\ast}\right)  \times\mathcal{Z}_{n+m}^{h}\left(  \theta_{\ast}\right)
\right)  ,\pi_{X}=\mu_{n},\int h\left(  \theta_{\ast},z\right)  \pi_{Z}\left(
dz\right)  =0\}\right.  .\nonumber
\end{align}

\end{definition}

The Implicit SOS formulation might lead to dimension reductions if $l$ ( the
dimension of the ambient space of $X$) is large. In addition, the presence of
$h\left(  \cdot\right)  $ in the distance evaluation allows the procedure to
use the available information in a more efficient way. For instance, if
$h\left(  \theta,x\right)  =\left\vert x\right\vert -\theta$, then the sign of
$x$ is irrelevant for the estimation problem and this will have the effect of
increasing the power of the Implicit SOS function relative to the explicit counterpart.

The analysis of the Implicit SOS function follows as a direct consequence of
Theorem \ref{SOSTheoremMean}; just redefine $X_{i}\leftarrow h\left(
\theta_{\ast},X_{i}\right)  $, $Z_{k}\leftarrow h\left(  \theta_{\ast}%
,Z_{k}\right)  $, and apply Theorem \ref{SOSTheoremMean} directly. Thus the
proof of the next result is omitted.

\begin{theorem}
[Implicit SOS Profile Function Analysis]\label{SOSGeneralImplicit}Let us
denote $g_{X}(\cdot)$ as the density for $h\left(  \theta_{\ast},X_{i}\right)
\in R^{q}$ and $g_{Y}(\cdot)$ for the density of $h\left(  \theta_{\ast}%
,Y_{i}\right)  \in R^{q}$. Then, the Wasserstein profile function defined in
Equation (\ref{WassersteinProfileFunctionGeneralized}) has the following
asymptotic results:

\begin{itemize}
\item When $q=1$,
\[
nR_{n}^{W}(\theta_{\ast})\Rightarrow Var\left(  h\left(  \theta_{\ast}%
,X_{1}\right)  \right)  \chi_{1}^{2}%
\]

\item When $q=2$, if $\tilde{Z}\sim N\left(  0,Var\left(  h\left(
\theta_{\ast},X\right)  \right)  \right)  \in R^{q}$ then%
\[
nR_{n}^{W}(\theta_{\ast})\Rightarrow\rho\left(  \tilde{Z}\right)  \left[
2-\eta\left(  \tilde{Z}\right)  \rho\left(  \tilde{Z}\right)  \right]
\left\Vert \tilde{Z}\right\Vert _{2}^{2},
\]
where $\rho\left(  \tilde{Z}\right)  $ is the unique solution to the equation
\[
\frac{1}{\rho}=\tilde{g}\left(  \rho\tilde{Z}\right)  ,
\]
and $\tilde{g}:\mathbb{R}^{q}\rightarrow\mathbb{R}$ is a deterministic
function defined as
\[
\tilde{g}\left(  x\right)  =\mathbb{P}\left(  \left\Vert x\right\Vert _{2}%
^{2}\geq\tau(0)\right)  ,
\]
where $\tau$ is a random variable satisfying
\[
\mathbb{P}\left[  \tau>t\right]  =\mathbb{E}\left[  \exp\left(  -\left[
g_{X}\left(  h\left(  \theta_{\ast},X_{1}\right)  \right)  +\kappa
g_{Y}\left(  h\left(  \theta_{\ast},X_{1}\right)  \right)  \right]  {\pi
}t\right)  \right]  .
\]
And the function $\tilde{\eta}:\mathbb{R}^{q}\rightarrow\mathbb{R}$ is a
deterministic continuous function given as
\[
\tilde{\eta}\left(  x\right)  =\mathbb{E}\left[  \max\left(  1-{\tau(0)
}/{\left\Vert x\right\Vert _{2}^{2}},0\right)  \right]  .
\]

\item When $q\geq3$,
\[
n^{1/2+\frac{3}{2q+2}}R_{n}^{W}(\theta_{\ast})\Rightarrow\frac{2q+2}{q+2}%
\frac{\left\Vert \tilde{Z}\right\Vert _{2}^{1+\frac{1}{q+1}}}{\left(
\mathbb{E}\left[  \frac{\pi^{q/2}}{\Gamma(q/2+1)}\left(  g_{X}\left(  h\left(
\theta_{\ast},X_{1}\right)  \right)  +\kappa g_{Y}\left(  h\left(
\theta_{\ast},X_{1}\right)  \right)  \right)  \right]  \right)  ^{\frac
{1}{q+1}}}%
\]
where $\tilde{Z}\sim N\left(  0,Var\left(  h\left(  \theta_{\ast},X\right)
\right)  \right)  \in\mathbb{R}^{q}$.
\end{itemize}
\end{theorem}

\subsubsection{Explicit SOS Formulation for Estimating Equations}

The second SOS function form we call Explicit SOS function because the
Wasserstein distance is explicitly or directly applied to the samples and the scenarios.

\begin{definition}
[Explicit SOS Profile Function for Estimating Equations]%
\label{ExplicitSOSFormulation}%
\begin{align}
\label{SOSFunctionGeneralExplicit}R_{n}^{W}(\theta_{\ast}) =  &  \inf
\{\int\int\left\Vert x-z\right\Vert _{2}^{2}\pi\left(  dx,dz\right)  :\\
&  \text{s.t. }\left.  \pi\in\mathcal{P}\left(  \mathcal{X}_{n}\times
\mathcal{Z}_{(n+m)}\right)  ,\pi_{X}=\mu_{n},\int h\left(  \theta_{\ast
},z\right)  \pi_{Z}\left(  dz\right)  =0\}\right.  .\nonumber
\end{align}

\end{definition}
Both the implicit and explicit SOS formulations have their merits. We have discussed the
merit of the implicit SOS formulation. For the Explicit SOS formulation,
consider the stress-testing application discussed in Section \ref{Sect_SS}.
The interest of an auditor or a regulator might be on the impact of scenarios
on a specific performance measure of interest. One might think that the
regulator applies the same stress scenarios to different insurance companies
or banks, and therefore the function $h\left(  \cdot\right)  $ is unique to
each insurance company. The regulator is interested in the impact of stress-testing scenarios on the structure of the bank (modeled by $h\left(
\cdot\right)  $). In this setting, the Explicit SOS formulation appears more appropriate.

While the analysis of the Explicit SOS formulation is also largely based on
the techniques developed for Theorem \ref{SOSTheoremMean}, it does require
some additional assumptions that are not immediately clear without examining
the proof of Theorem \ref{SOSTheoremMean}. In particular, in addition to
\textbf{A1)}, \textbf{A2)}, \textbf{B1)} and \textbf{B2)}, here we impose the
following assumptions.

\textbf{BE1)} Assume that the derivative of $h\left(  \theta_{\ast},x\right)
$ with respect to (w.r.t.) $x$, $D_{x}h\left(  \theta_{\ast},\cdot\right)
:R^{l}\rightarrow R^{q\times l}$, is continuous function of $x$ and the second
derivative w.r.t. $x$ is bounded, i.e., $\left\Vert D_{x}^{2}h\left(
\theta_{\ast},\cdot\right)  \right\Vert <\tilde{K}$ for all $x$.

\textbf{BE2)} Define $V_{i}=D_{x}h\left(  \theta_{\ast},X_{i}\right)  \cdot
D_{x}h\left(  \theta_{\ast},X_{i}\right)  ^{T}\in R^{q\times q}$ and assume
that $\Upsilon=\mathbb{E}\left(  V_{i}\right)  $ is strictly positive definite.

We provide the proof of the next result in our technical Section \ref{Sec_additional_proof}.

\begin{theorem}
[Explicit SOS Profile Function Analysis]\label{SOSGeneralExplicit} Under
assumptions A1)-A2), B1)-B2) and BE1)-BE2), we have that
(\ref{SOSFunctionGeneralExplicit}) satisfies

\begin{itemize}
\item When $l=1$,
\[
nR_{n}^{W}(\theta_{\ast})\Rightarrow\tilde{Z}^{T}\Upsilon^{-1}\tilde{Z}%
\]
where $\tilde{Z}\sim N\left(  0,Var\left(  h\left(  \theta_{\ast},X\right)
\right)  \right)  \in\mathbb{R}^{q}$.

\item Assume that $l=2$. Let $\tilde{Z}\sim N\left(  0,Var\left(  h\left(
\theta_{\ast},X\right)  \right)  \right)  \in\mathbb{R}^{q}$. It is possible
to uniquely define deterministic continuous mapping, $\tilde{\zeta}%
:\mathbb{R}^{q}\rightarrow\mathbb{R}^{q}$, such that $\tilde{\zeta}\left(
z\right)  $ is defined via%
\[
z=-\mathbb{E}\left[  V_{1}I\left(  \tau\leq\tilde{\zeta}^{T}\left(  z\right)
V_{1}\tilde{\zeta}\left(  z\right)  \right)  \right]  \tilde{\zeta}\left(
z\right)  ,
\]
where $\tau$ is independent of $V_{1}$ satisfying
\[
\mathbb{P}\left(  \tau>t\right)  =\mathbb{E}\left(  \exp\left(  -\left[
f_{X}\left(  X_{1}\right)  +\kappa f_{Y}\left(  X_{1}\right)  \right]  {\pi
}t\right)  \right)  .
\]
Then, we have that,
\[
nR_{n}^{W}(\theta_{\ast})\Rightarrow-2\tilde{Z}^{T}\tilde{\zeta}\left(
\tilde{Z}\right)  -\tilde{\zeta}^{T}\left(  \tilde{Z}\right)  \tilde{G}\left(
\tilde{\zeta}\left(  \tilde{Z}\right)  \right)  \tilde{\zeta}\left(  \tilde
{Z}\right)  ,
\]
where $\tilde{G}:\mathbb{R}^{q}\rightarrow\mathbb{R}^{q\times q}$ is a
deterministic continuous mapping defined as
\[
\tilde{G}\left(  \zeta\right)  =\mathbb{E}\left[  V_{1}\max\left(  1-{\tau
}/{\left(  \zeta^{T}V_{1}\zeta\right)  },0\right)  \right]  .
\]

\item Suppose that $l\geq3$. It is possible to uniquely define deterministic
continuous mapping $\tilde{\zeta}:\mathbb{R}^{q}\rightarrow\mathbb{R}^{q}$,
such that
\[
z=-\mathbb{E}\left[  \frac{\pi^{l/2}\left(  f_{X}\left(  X_{1}\right)  +\kappa
f_{Y}\left(  X_{1}\right)  \right)  }{\Gamma(l/2+1)}V_{1}\cdot\left(
\tilde{\zeta}^{T}\left(  z\right)  V_{1}\tilde{\zeta}\left(  z\right)
\right)  ^{l}\right]  \tilde{\zeta}\left(  z\right)  ,
\]
(note that $V_{1}$ is a function of $X_{1}$, so these are correlated).
Moreover,
\[
n^{1/2+\frac{3}{2l+2}}R_{n}^{W}(\theta_{\ast})\text{ }\Rightarrow-2\tilde
{Z}^{T}\tilde{\zeta}\left(  \tilde{Z}\right)  -\frac{2}{l+2}\tilde{G}\left(
\tilde{Z}\right)  ,
\]
where $\tilde{Z}\sim N\left(  0,Var\left(  h\left(  \theta_{\ast},X\right)
\right)  \right)  \in\mathbb{R}^{q}$ and $\tilde{G}:\mathbb{R}^{q}%
\rightarrow\mathbb{R}$ is a deterministic continuous function defined as
\[
\tilde{G}\left(  \zeta\right)  =\mathbb{E}\left[  \frac{\pi^{l/2}}%
{\Gamma(l/2+1)}\left(  f_{X}\left(  X_{1}\right)  +\kappa f_{Y}\left(
X_{1}\right)  \right)  \left(  \zeta^{T}V_{1}\zeta\right)  ^{l/2+1}\right]  .
\]

\end{itemize}
\end{theorem}

We should observe that unlike the implicit formulation, the rate of
convergence will only depend on the dimension of data $X_{i}\in\mathbb{R}^{l}%
$, but the shape of asymptotic distribution is determined by the estimating
functions $h\left(  \theta_{\ast},X_{i}\right)  \in\mathbb{R}^{q}$.

\subsection{Plug-in Estimators for SOS Functions}

In many situations, for example in the context of stochastic optimization, we
are interested in a specific parameter $\theta_{\ast}=\left(  \gamma_{\ast
},\nu_{\ast}\right)  \in R^{d+p}$ such that $\mathbb{E}\left[  h\left(
\gamma_{\ast},\nu_{\ast},X\right)  \right]  =0$, where $\nu_{\ast}\in R^{p}$
is the nuisance parameter (for example Lagrange multipliers in the setting of
constrained optimization).

We shall discuss a method that allows us to deal with the nuisance parameter
using a plug-in estimator, while taking advantage of the SOS framework for the
estimation of $\gamma_{\ast}$. After we state our assumptions we will provide
the results in this section, and the proofs, which follow closely those of
Theorem \ref{SOSGeneralImplicit} and Theorem \ref{SOSGeneralExplicit}, will be given
in Section \ref{Section_Methodology}.

Throughout this subsection, let us suppose that $h\left(  \gamma,\nu,x\right)
\in\mathbb{R}^{q}$. In addition, we impose the following assumptions.

\textbf{C1)} Given $\gamma_{\ast}$ there is a unique $\nu_{\ast}\in R^{p}$
such that
\begin{equation}
\mathbb{E}\left[  h\left(  \gamma_{\ast},\nu,X\right)  \right]  =0
\label{Eq_h2}%
\end{equation}
and, given $\nu_{\ast}$, we also assume that $\gamma_{\ast}$ satisfies%
\begin{equation}
\mathbb{E}\left[  h\left(  \gamma,\nu_{\ast},X\right)  \right]  =0.
\label{Eq_h1}%
\end{equation}

\textbf{C2)} We have access to a suitable estimator $v_{n}$ such that the
sequence
\[
\left\{  n^{1/2}\left(  v_{n}-\nu_{\ast}\right)  \right\}  _{n=1}^{\infty
}\text{ is tight,}%
\]
and
\[
\frac{1}{\sqrt{n}}\sum_{i=1}^{n}h\left(  \gamma_{\ast},v_{n},X_{i}\right)
\Rightarrow\tilde{Z}^{\prime},
\]
for some random variable $\tilde{Z}^{\prime}$, as $n\rightarrow\infty$.

\textbf{C3)} Assume that $h\left(  \gamma,{\cdot},x\right)  $ is continuously
differentiable a.e. (almost everywhere with respect to the Lebesgue measure)
in some neighborhood $\mathcal{V}$ around $v_{\ast}$.

\textbf{C4)} Suppose that there is a function $M\left(  \cdot\right)
:R^{l}\rightarrow(0,\infty)$ satisfying that%
\begin{align*}
\left\Vert h\left(  \gamma_{\ast},{\nu},x\right)  \right\Vert _{2}^{2}  &
\leq M(x)\text{ for a.e. }\nu\in\mathcal{V},\\
\left\Vert D_{\nu}h\left(  \gamma_{\ast},{\nu},x\right)  \right\Vert _{2}^{2}
&  \leq M(x)\text{ for a.e. }\nu\in\mathcal{V},
\end{align*}
and $E\left(  M\left(  X_{1}\right)  \right)  <\infty$ and $E\left(  M\left(  Y_{1}\right)  \right)  <\infty$.

\subsubsection{\textbf{Plug-in Estimators for Implicit SOS Functions}}

We are interested in studying the plug-in implicit SOS function (or implicit
pseudo-SOS profile function) given by
\begin{align}
R_{n}^{W}(\gamma_{\ast})  &  =\inf\{\int\int\left\Vert h\left(  \gamma_{\ast
},v_{n},x\right)  -h\left(  \gamma_{\ast},v_{n},z\right)  \right\Vert _{2}%
^{2}\pi\left(  dx,dz\right)  :\label{PlugInWPEE}\\
&  \text{s.t. }\left.  \pi\in\mathcal{P}\left(  \mathcal{X}_{n}^{h}\left(
\gamma_{\ast},v_{n}\right)  \times\mathcal{Z}_{n+m}^h\left(  \gamma_{\ast
},v_{n}\right)  \right)  ,\pi_{X}=\mu_{n},\int h\left(  \gamma_{\ast}%
,v_{n},z\right)  \pi_{Z}\left(  dz\right)  =0\},\right. \nonumber
\end{align}
where
\[
\mathcal{X}_{n}^{h}\left(  \gamma_{\ast},v_{n}\right)  =\{h\left(
\gamma_{\ast},v_{n},x\right)  :x\in\mathcal{X}_{n}\},\text{ }\mathcal{Z}%
_{n+m}^{h}\left(  \gamma_{\ast},v_{n}\right)  =\{h\left(  \gamma_{\ast}%
,v_{n},z\right)  :z\in\mathcal{Z}_{(n+m)}\}.
\]
We typically will use (\ref{Eq_h2}) to find a plug-in estimator $v_{n}$. Under
suitable assumptions on the consistency and convergence rate of the plug-in
estimator, we have an asymptotic result for (\ref{PlugInWPEE}), as we indicate next.

\begin{corollary}
[Plug-in for Implicit SOS Formulation]\label{SOSGeneralImplicitPlugIn}

Assume \textbf{A1)}-\textbf{A2)} and \textbf{C1)}-\textbf{C4)} hold.
Moreover, suppose we denote $g_{X}(\cdot)$ as the density for $h\left(
\gamma_{\ast},v_{\ast},X_{i}\right)  \in R^{q}$ and $g_{Y}(\cdot)$ for the
density of $h\left(  \gamma_{\ast},v_{\ast},Y_{i}\right)  \in R^{q}$. We
notice $\tilde{Z}^{\prime}\in\mathbb{R}^{q}$ is defined in C2). We obtain that
(\ref{PlugInWPEE}) has following asymptotic behavior:

\begin{itemize}
\item When $q=1$,
\[
nR_{n}^{W}(\gamma_{\ast})\Rightarrow\left(  \tilde{Z}^{\prime}\right)  ^{2}.
\]

\item When $q=2$,
\[
nR_{n}^{W}(\gamma_{\ast})\Rightarrow\rho\left(  \tilde{Z}^{\prime}\right)
\left[  2-\tilde{\eta}\left(  \tilde{Z}^{\prime}\right)  \rho\left(  \tilde
{Z}^{\prime}\right)  \right]  \left\Vert \tilde{Z}^{\prime}\right\Vert
_{2}^{2}%
\]
where $\rho\left(  \tilde{Z}^{\prime}\right)  $ is the unique solution to the equation
\[
\frac{1}{\rho}=\tilde{g}\left(  \rho\tilde{Z}^{\prime}\right)  ,
\]
and $\tilde{g}:\mathbb{R}^{q}\rightarrow\mathbb{R}$ is a deterministic
continuous function defined as
\[
\tilde{g}\left(  x\right)  =\mathbb{P}\left(  \left\Vert x\right\Vert _{2}%
^{2}\geq\tau\right)  .
\]
The function $\tilde{\eta}:\mathbb{R}^{q}\rightarrow\mathbb{R}$ is a
deterministic continuous function defined as
\[
\tilde{\eta}\left(  x\right)  =\mathbb{E}\left[  \max\left(  1-{\tau
}/{\left\Vert x\right\Vert _{2}^{2}},0\right)  \right]  .
\]
Moreover, $\tau$
satisfies
\[
\mathbb{P}\left[  \tau>t\right]  =\mathbb{E}\left[  \exp\left(  -\left[
g_{X}\left(  h\left(  \gamma_{\ast},\nu_{\ast},X_{1}\right)  \right)  +\kappa
g_{Y}\left(  h\left(  \gamma_{\ast},\nu_{\ast},X_{1}\right)  \right)  \right]
{\pi}t\right)  \right]  .
\]

\item When $q\geq3$,
\[
n^{1/2+\frac{3}{2q+2}}R_{n}^{W}(\gamma_{\ast})\Rightarrow\frac{2q+2}{q+2}%
\frac{\left\Vert \tilde{Z}^{\prime}\right\Vert _{2}^{1+\frac{1}{q+1}}}{\left(
\mathbb{E}\left[  \frac{\pi^{q/2}}{\Gamma(q/2+1)}\left(  g_{X}\left(  h\left(
\gamma_{\ast},\nu_{\ast},X_{1}\right)  \right)  +\kappa g_{Y}\left(  h\left(
\gamma_{\ast},\nu_{\ast},X_{1}\right)  \right)  \right)  \right]  \right)
^{\frac{1}{q+1}}}.
\]

\end{itemize}
\end{corollary}

\subsubsection{\textbf{Plug-in Estimators for Explicit SOS Functions}}

We can also analyze plug-in estimators for Explicit SOS profile functions. We
now define the explicit plug-in (or pseudo) SOS function based on
(\ref{SOSFunctionGeneralExplicit})\ as simply plugging in the nuisance
parameter:
\begin{align}
R_{n}^{W}(\gamma_{\ast})  &  =\inf\big\{\int\int\left\Vert x-z\right\Vert
_{2}^{2}\pi\left(  dx,dz\right)  :\label{PlugESOS}\\
&  \text{s.t. }\left.  \pi\in\mathcal{P}\left(  \mathcal{X}_{n}\times
\mathcal{Z}_{(n+m)}\right)  ,\pi_{X}=\mu_{n},\int h\left(  \gamma_{\ast}%
,v_{n},z\right)  \pi_{Z}\left(  dz\right)  =0\big\}.\right. \nonumber
\end{align}

In addition to \textbf{C1)} to \textbf{C4)}\ introduced at the beginning of
this subsection, we shall impose the following additional assumptions:

\textbf{C5)} Define $\bar{V}_{i}\left(  v_{\ast}\right)  =D_{x}h\left(
\gamma_{\ast},{\nu}_{\ast},X_{i}\right)  \cdot D_{x}h\left(  \gamma_{\ast
},{\nu}_{\ast},X_{i}\right)  ^{T}$ and assume that $\bar{\Upsilon}%
=\mathbb{E}\left(  \bar{V}_{i}\right)  $ is strictly positive definite.

\textbf{C6) }The function $M\left(  \cdot\right)  $ from condition
\textbf{C4)} also satisfies%
\begin{align*}
\left\Vert D_{x}h\left(  \gamma_{\ast},{\nu},x\right)  \right\Vert _{2} ^{2}
&  \leq M(x)\text{ for a.e. }\nu\in\mathcal{V}.\\
\left\Vert D_{\nu}D_{x}h\left(  \gamma_{\ast},{\nu},x\right)  \right\Vert _{2}
^{2}  &  \leq M(x)\text{ for a.e. }\nu\in\mathcal{V}.
\end{align*}

\textbf{C7)} The second derivative w.r.t. $x$ exist and bounded, i.e.,
$\left\Vert D^{2}_{x}h\left(  \gamma_{\ast},{\nu},x\right)  \right\Vert <
\tilde{K} \text{ for a.e. }\nu\in\mathcal{V}$ and all $x$.

\begin{corollary}
[Plug-in for Explicit SOS Formulation]\label{SOSGeneralExplicitPlugIn}
Let $X_{i}\in\mathbb{R}^{l}$, $h\left(  \gamma,\nu,x\right)  \in\mathbb{R}^{q}$, and 
assume that A1)-A2) and C1)-C7) hold. We notice $\tilde{Z}^{\prime}$ is
defined in C2). Then, the SOS profile function defined in Equation
(\ref{PlugESOS}) has the following asymptotic properties:

\begin{itemize}
\item When $l=1$,
\[
nR_{n}^{W}(\gamma_{\ast})\Rightarrow\tilde{Z}^{\prime T}\bar{\Upsilon}%
^{-1}\tilde{Z}^{\prime}.
\]

\item Suppose that $l=2$. It is possible to uniquely define deterministic
continuous mapping $\tilde{\zeta}:\mathbb{R}^{q}\rightarrow\mathbb{R}^{q}$,
such that
\[
z=-\mathbb{E}\left[  \bar{V}_{1}I\left(  \tau\leq\tilde{\zeta}^{T}\left(
z\right)  \bar{V}_{1}\tilde{\zeta}\left(  z\right)  \right)  \right]
\tilde{\zeta}\left(  z\right)  ,
\]
where $\tau$ is independent of $\bar{V}_{1}$ and it satisfies%
\[
\mathbb{P}\left(  \tau>t\right)  =\mathbb{E}\left(  \exp\left(  -\left[
f_{X}\left(  X_{1}\right)  +\kappa f_{Y}\left(  X_{1}\right)  \right]  {\pi
}t\right)  \right)  .
\]
Furthermore,
\[
nR_{n}^{W}(\gamma_{\ast})\Rightarrow-2\tilde{\zeta}^{T}\left(  \tilde
{Z}^{\prime}\right)  \tilde{Z}^{\prime}-\tilde{\zeta}^{T}\left(  \tilde
{Z}^{\prime}\right)  \tilde{G}\left(  \tilde{\zeta}\left(  \tilde{Z}^{\prime
}\right)  \right)  \tilde{\zeta}\left(  \tilde{Z}^{\prime}\right)  ,
\]
where $\tilde{G}:\mathbb{R}^{q}\rightarrow\mathbb{R}^{q\times q}$ is a
deterministic continuous mapping defined as
\[
\tilde{G}\left(  \zeta\right)  =\mathbb{E}\left[  \bar{V}_{1}\max\left(
1-{\tau}/{\left(  \zeta^{T}\bar{V}_{1}\zeta\right)  },0\right)  \right]  .
\]

\item Assume that $l\geq3$. A deterministic and continuous mapping
$\tilde{\zeta}:\mathbb{R}^{q}\rightarrow\mathbb{R}^{q}$ can be defined
uniquely so that
\[
z=-\mathbb{E}\left[  \frac{\pi^{l/2}\left(  f_{X}\left(  X_{1}\right)  +\kappa
f_{Y}\left(  X_{1}\right)  \right)  }{\Gamma(l/2+1)}\bar{V}_{1}\left(
\tilde{\zeta}^{T}\left(  z\right)  \bar{V}_{1}\tilde{\zeta}\left(  z\right)
\right)  ^{l}\right]  \tilde{\zeta}\left(  z\right)
\]
(note that $\bar{V}_{1}$ is a function of $X_{1}$). Moreover,
\[
n^{1/2+\frac{3}{2l+2}}R_{n}^{W}(\gamma_{\ast})\Rightarrow-2\tilde{\zeta}%
^{T}\left(  \tilde{Z}^{\prime}\right)  \tilde{Z}^{\prime}-\frac{2}{l+2}%
\tilde{G}\left(  \tilde{\zeta}\left(  \tilde{Z}^{\prime}\right)  \right)  ,
\]
where $\tilde{G}:\mathbb{R}^{q}\rightarrow\mathbb{R}$ is a deterministic
continuous function defined as
\[
\tilde{G}\left(  \zeta\right)  =\mathbb{E}\left[  \frac{\pi^{l/2}}%
{\Gamma(l/2+1)}\left(  f_{X}\left(  X_{1}\right)  +\kappa f_{Y}\left(
X_{1}\right)  \right)  \left(  \zeta^{T}\bar{V}_{1}\zeta\right)
^{l/2+1}\right]  .
\]

\end{itemize}
\end{corollary}

\section{Application to Stochastic Optimization and Stress
Testing\label{Section_Applications}}

We will provide an application of the SOS inference framework to quantify
model uncertainty in the context of stochastic programming. Motivating
applications include the evaluation of Conditional Value at Risk (C-VaR) and
semi-supervised learning settings, as we shall discuss in the examples below.

We are interested in the value function of a stochastic programming problem
formulation via
\begin{align}
\label{OPT_APP}C_{\ast}=  &  \min_{\theta}\text{ }\text{ }\mathbb{E}\left[
m(\theta,X)\right] \\
&  s.t.\text{ }\mathbb{E}[\phi(\theta,X)]\leq0.\nonumber
\end{align}
We assume that the objective function $\psi(\theta)=\mathbb{E}\left[
m(\theta,X)\right]  $ is a convex function in $\theta$; while the constraints
$\mathbb{E}[\phi(\theta,X)]\leq0$ specify a convex region in $\theta$; for
example we shall assume that $\phi(\cdot,x)$ is a convex function for any $x$.

Following \cite{blanchet_quantifying_2016}, the goal is to estimate the optimal
value function using the SOS\ formulation and we will apply a plug-in
estimator for $\theta_{\ast}$ (which is treated as a nuisance parameter).
Subsequently, when introducing the Lagrangian relaxation of \eqref{OPT_APP} we
will be able to also introduce a plug-in estimator for the associated Lagrange
multiplier. Therefore, for simplicity, we shall focus on the unconstrained
minimization problem $C_{\ast}=\min_{\theta}\left\{  \mathbb{E}\left[
m(\theta,X)\right]  \right\}  $.

The authors in \cite{lam_quantifying_2015,lam_empirical_2016} provide a
discussion for some potential approaches to derive nonparametric confidence
interval (including Empirical Likelihood, a Bayesian approach, bootstrap and
the delta method). In \cite{lam_quantifying_2015,lam_empirical_2016} it is
argued that the Empirical Likelihood method tends to have superior finite
sample performance, and \cite{blanchet_quantify_2016} provides an optimal (in
certain sense) specification for the Empirical Likelihood approach. More
importantly, in \cite{blanchet_quantify_2016} an approach combining Empirical
Likelihood and a plug-in estimator for the optimizer is introduced, which
avoids solving a non-convex optimization problem introduced in the discussion
of \cite{lam_quantifying_2015}.

Our goal in this section is to derive a plug-in estimator based on the SOS
inference approach introduced in Section \ref{Section_Main_Results}. The
approach that we introduce next is the analog of the plug-in strategy
discussed in \cite{blanchet_quantify_2016} in order to find a robustified
confidence interval for $C_{\ast}$.

The following corollary plays the key role in specifying confidence interval
for $C_{\ast}$. The result is a direct extension of Corollary
\ref{SOSGeneralImplicitPlugIn} and Corollary \ref{SOSGeneralExplicitPlugIn},
provided the following assumptions are in place.

We define $M\left(  \theta\right) = \mathbb{E}\left[m(\theta,X)\right]$, and the assumptions are

\textbf{D1): }Assume $m\left(  \cdot\right)  $ is convex differentiable in
$\theta$, then $M\left(  \theta\right)$ is also convex differentiable. We assume there is a unique optimizer $\theta_{\ast}$n for  $M\left(  \theta\right)$.

\textbf{D2): }Assume that $m\left(  \cdot\right)  $ is strongly convex at
$\theta_{\ast}$, that is, there exist $\delta>0$, such
that for every $\theta$ 
\[
M\left(  \theta\right)  \geq M\left(  \theta_{\ast}\right)  +\delta\left\Vert
\theta-\theta_{\ast}\right\Vert _{2}^{2}.
\]

\begin{corollary}
\label{SOSGeneralPlugIn} Let us consider stochastic programming problem
$C_{\ast}=\min_{\theta}M\left(  \theta\right)  =\min_{\theta}\mathbb{E}\left[
m(\theta,X)\right]  $. Assume that D1)-D2) hold. We consider the estimating
equations to be the derivative condition and value function condition
\[
\mathbb{E}\left[  m(\theta_{\ast},X)-C_{\ast}\right]  =0,\text{ and
}\mathbb{E}\left[  D_{\theta}m\left(  \theta_{\ast},X\right)  \right]  =0.
\]
For simplicity, let us denote $h\left(  \theta_{\ast},C_{\ast},x\right)
=\left(  m(\theta_{\ast},x)-C_{\ast},D_{\theta}m\left(  \theta_{\ast
},x\right)  ^{T}\right)  ^{T}$. We are interested in $C_{\ast}$ only and
consider a sample average approximation (SAA) estimator for $\theta_{\ast}$ to
be $\hat{\theta}_{SAA}$. For $h\left(  \cdot,C_{\ast},x\right)  $ we assume
C1)-C7) hold. Let us denote $U\sim{N}\left(  0,\text{Var}\left(  m\left(
\theta_{\ast},X\right)  \right)  \right)  \in\mathbb{R}$ and $U(0)=\left(
U,\vec{0}\right)  ^{T}\in\mathbb{R}^{d+1}$. Recalling the implicit and explicit
formulations for general estimating equation SOS function defined in
Definition \ref{ImplicitSOSFormulation} and Definition
\ref{ExplicitSOSFormulation}, we have the following asymptotic
results.\newline For the implicit SOS formulation, we have

\begin{itemize}
\item When $d=1$ (estimating equation dimension is $d+1=2$)
\[
n{R}_{n}^{W}(C_{\ast})\Rightarrow\rho\left(  U\right)  \left[  2-\tilde{\eta
}\left(  U\right)  \rho\left(  U\right)  \right]  U^{2,}%
\]
where $\rho\left(  U\right)  $ is the unique solution to
\[
\frac{1}{\rho}=\tilde{g}\left(  \rho U\right)  ,
\]
and $\tilde{g}:\mathbb{R}\rightarrow\mathbb{R}$ is a deterministic continuous
function defined as
\[
\tilde{g}\left(  x\right)  =\mathbb{P}\left[  x^{2}\geq\tau\right]  .
\]
$\tilde{\eta}\left(  x\right)  $ is also a deterministic function, defined as
\[
\tilde{\eta}\left(  x\right)  =\mathbb{E}\left[  \max\left(  1-{\tau}/{x^{2}%
},0\right)  \right]  ,
\]
and $\tau$ satisfies
\[
\mathbb{P}\left[  \tau>t\right]  =E\left(  \exp\left(  -\left(g_X\left(  h\left(
\theta_{\ast},C_{\ast},X_{1}\right)  + \kappa g_Y\left(  h\left(
\theta_{\ast},C_{\ast},X_{1}\right) \right) \right)  {\pi}t\right)  \right)  \right) .
\]

\item When $d\ge2$,
\begin{align*}
n^{1/2+\frac{3}{2d+4}}{R}_{n}^{W}( C_{\ast}) \Rightarrow\frac{2d+4}{d+3}
\frac{\left|  \left|  U\right|  \right|_2  ^{1+\frac{1}{d+2}}}{\mathbb{E}\left[
\frac{\pi^{(d+1)/2}} {\Gamma((d+3)/2)}\left(g_{X}\left(  h\left(  \theta_{\ast
},C_{\ast},X_{1}\right)  \right)  + g_{Y}\left(  h\left(  \theta_{\ast
},C_{\ast},X_{1}\right)  \right) \right) \right]  ^{\frac{1}{d+2}}}.
\end{align*}

\end{itemize}

For the explicit formulation, we have the following asymptotic results (we use
${\zeta}_{[1]}$ to denote the first element of vector $\zeta$)

\begin{itemize}
\item When $l=1$,
\begin{align*}
n{R}_{n}^{W}( C_{\ast})\Rightarrow v_{1,1}U^{2},
\end{align*}
where $v_{1,1}$ is the $\left(  1,1\right)  $ element of matrix $\Upsilon^{-1}
$.

\item Suppose that $l=2$. It is possible to uniquely define deterministic
continuous mapping $\tilde{\zeta}:\mathbb{R}^{q}\rightarrow\mathbb{R}^{q}$,
such that
\[
z=-\mathbb{E}\left[  \bar{V}_{1}I\left(  \tau\leq\tilde{\zeta}^{T}\left(
z\right)  \bar{V}_{1}\tilde{\zeta}\left(  z\right)  \right)  \right]
\tilde{\zeta}\left(  z\right)  ,
\]
where $\tau$ is independent of $U$ satisfying
\[
\mathbb{P}\left(  \tau>t\right)  =\mathbb{E}\left(  \exp\left(  -\left[
f_{X}\left(  X_{1}\right)  +\kappa f_{Y}\left(  X_{1}\right)  \right]  {\pi
}t\right)  \right)  .
\]
Furthermore,
\[
nR_{n}^{W}(C_{\ast}) \Rightarrow-2U\tilde{\zeta}_{[1]} -\tilde{\zeta}%
^{T}\left(  U(0)\right)  \tilde{G}\left(  \tilde{\zeta}\left(  U(0)\right)\right)  \tilde{\zeta
}\left(  U(0)\right)  ,
\]
where $\tilde{G}:\mathbb{R}^{q}\to\mathbb{R}^{q\times q}$ is a deterministic
continuous mapping defined as
\[
\tilde{G}\left(  \zeta\right)  = \mathbb{E}\left[  \bar{V}_{1} \max\left(
1-\frac{\tau}{\zeta^{T}\bar{V}_{1}\zeta},0\right)  \right]  ,
\]

and $U$ is independent with $\bar{V}_{1}$ and $\tau$.

\item Assume that $l\geq3$. A continuous function $\tilde{\zeta}
:\mathbb{R}^{q}\rightarrow\mathbb{R}^{q}$ can be defined uniquely so that
\[
z=-\mathbb{E}\left[  \frac{\pi^{l/2}\left(  f_{X}\left(  X_{1}\right)  +\kappa
f_{Y}\left(  X_{1}\right)  \right)  }{\Gamma(l/2+1)}\bar{V}_{1}\left(
\tilde{\zeta}^{T}\left(  z\right)  \bar{V}_{1}\tilde{\zeta}\left(  z\right)
\right)  ^{l}\right]  \tilde{\zeta}\left(  z\right)
\]
(note that $\bar{V}_{1}$ is a function of $X_{1}$). Moreover,
\[
n^{1/2+\frac{3}{2l+2}}R_{n}^{W}(C_{\ast})\Rightarrow-2 U \tilde{\zeta
}_{[1]} -\frac{2}{l+2}\tilde{G}\left(  \tilde{\zeta}\left(  U(0)\right)
\right)  ,
\]
where $\tilde{G}:\mathbb{R}^{q}\to\mathbb{R}$ is a deterministic function
given as
\[
\tilde{G}\left(  \zeta\right)  = \mathbb{E}\left[  \frac{\pi^{l/2}}%
{\Gamma(l/2+1)}\left( f_{X}\left(  X_{1}\right) + \kappa f_{Y}\left(  X_{1}\right)\right)  \left(  {\zeta}^{T}\bar{V}%
_{1}{\zeta}\right)  ^{l/2+1}\right]  ,
\]
and $U$ and $X_{1}$ are independent.
\end{itemize}
\end{corollary}

As indicated earlier, the corollary is a special case of Corollary
\ref{SOSGeneralImplicitPlugIn} and Corollary \ref{SOSGeneralExplicitPlugIn},
so the proof is omitted. The estimating equations correspond to the first
order optimality condition (i.e., the first derivative equal to zero) and the corresponding optimal value equation. We use sample average
approximation estimator as the underlying plug-in estimator.

We notice that for sample average approximation, under assumptions D1)-D2), it
has been shown in \cite{ruszczynski_stochastic_2003,shapiro_lectures_2014}
that the optimizer $\hat{\theta}_{SAA}$ and the optimal value function
$\frac{1}{n}\sum_{i=1}^{n}m\left(  \hat{\theta}_{SAA},X_{i}\right)  $ satisfy
\begin{align*}
&  \hat{\theta}_{SAA}-\theta_{\ast}=O\left(  1/n^{1/2}\right) \\
&  \frac{1}{n}\sum_{i=1}^{n}\nabla_{\theta}m\left(  \hat{\theta}_{SAA}%
,X_{i}\right)  =0,\\
&  \frac{1}{\sqrt{n}}\sum_{i=1}^{n}\left(  m\left(  \hat{\theta}_{SAA}%
,X_{i}\right)  -C_{\ast}\right)  \Rightarrow N\left(  0,\text{Var}\left(
m\left(  \theta_{\ast},X\right)  \right)  \right)  .
\end{align*}
Therefore, Corollary \ref{SOSGeneralExplicitPlugIn} and Corollary
\ref{SOSGeneralImplicitPlugIn} apply.

\bigskip

Similar to the derivation in \cite{blanchet_quantify_2016} in the setting of
Empirical Likelihood, for the plug-in estimator derived from sample average
approximation, if we denote $n^{1/2+3/(2d+4)}R_{n}^{W(implicit)}(C_{\ast
})\Rightarrow R_{0}^{(implicit)}$ and $n^{1/2+3/(2l+2)}R_{n}^{W(explicit)}%
(C_{\ast})\Rightarrow R_{0}^{(explicit)}$, we can specify a robust $95\%$
confidence interval for $C_{\ast}$ under both explicit and implicit
formulation by:
\[
\text{CI}^{(\cdot)}\left(  C_{\ast}\right)  =\left\{  C\in\mathbb{R}%
\middle|n^{\alpha}R_{n}^{W(\cdot)}\left(  C\right)  \leq R_{0}^{(\cdot
)}\left(  95\%\right)  \right\}
\]
where $\alpha$ depends on the formulation and dimension as in Corollary
\ref{SOSGeneralPlugIn} and $R_{0}^{(\cdot)}\left(  95\%\right)  $ is the upper
$95\%$ quantile for $R_{0}^{(explicit)}$ (or $R_{0}^{(implicit)}$). The
upper/lower bound of confidence interval ($C_{up}^{\left(  \cdot\right)  }%
$/$C_{l0}^{\left(  \cdot\right)  }$) can be found by solving the linear
programming problem
\begin{align*}
C_{up}^{\left(  \cdot\right)  }/C_{lo}^{\left(  \cdot\right)  }=  &  \max
_{\pi(i,j)}\text{ }/\text{ }\min_{\pi(i,j)}\{\sum_{i,j=1}^{n}\pi
(i,j)m(\hat{\theta}_{SAA},X_{i})\\
&  s.t.\text{ }\pi(i,j)\geq0\;\sum_{j=1}^{n}\pi(i,j)=1/n;\sum_{i,j=1}^{n}%
\pi(i,j)\left\Vert X_{i}-X_{j}\right\Vert _{2}^{2}\leq \frac{R_{0}^{(\cdot)}\left(
95\%\right)}{n^\alpha}  \}.
\end{align*}

\bigskip

Next, we are going to provide a numerical example in quantifying C-VaR using
the methodology we developed above.

\begin{example}
[Quantify the uncertainty of Conditional Value at Risk (C-VaR)]In this example
we would like to find an SOS-based $95\%$ confidence interval for
conditional value at risk with $90\%$ level. The conditional value at risk
with $\alpha-$level is given as solving the stochastic programming problem:
\[
\text{C-VaR}(\alpha)=\inf_{\theta}\left\{  \theta+\frac{1}{1-\alpha}
\mathbb{E}\left[  \left(  \sum_{k=1}^{l}X^{(k)}-\theta\right)  ^{+}\right]
\right\}  .
\]
We shall test our method using simulated data under different distributional
assumptions. We sample i.i.d. observations $\left\{  X_{i}\right\}
_{i=1}^{n}\subset\mathbb{R}^{l}$. We will apply the SOS inference procedure to
provide a non-parametric confidence interval for $\text{C-VaR}(90\%)$. In
order to verify the coverage probability we use data simulated from normal
distribution and Laplace (double exponential) distributions. We consider the
case $l=4$. For the normal distribution setting we assume $X_{i}\sim{N}\left(
0,I_{4\times4}\right)  $, while for Laplace distribution we consider for each
$k=1,...,4$, $X_{i}^{k}\sim Laplace(0,1)$ and all of these random variables
are independent. For these two cases, we can calculate the solution in closed
form; for the normal setting the optimizer is $\theta^{\star}=2.5632$ and
optimal value function is $\text{C-VaR}(0.9)=3.510 $; for Laplace setting the
optimizer is $\theta^{\star}=3.497$ with optimal value function equal to
$\text{C-VaR}(0.9)=5.066$.

In this example, we have three approaches in which our SOS procedure can
be applied: 1) implicit SOS formulation (ISOS); 2) explicit SOS formulation
assuming data being of dimension $l$ (ESOS-O), i.e. $X_{i}=\left(  X_{i}^{(1)}
,\dots,X_{i}^{(l)}\right)  ^{T}\in\mathbb{R}^{l}$; 3) explicit formulation
assuming data being of dimension $1$ (ESOS-C), i.e. $X_{i}=X_{i}^{(1)}
+\ldots+X_{i}^{(l)}\in\mathbb{R}$. We compare our methods with empirical
likelihood method (EL) in \cite{blanchet_quantify_2016}, nonparametric
bootstrap method (BT), and central limit theorem-based Delta method (CLT) discussed in Theorem
5.7 \cite{shapiro_lectures_2014}. We consider four settings $n = 20, 50, 100
\text{ and } 500$. For each setting, we repeat the experiment $N=1000$ times,
and note down the empirical coverage probability, mean of upper and lower
bounds, and the mean and standard deviation of the interval width for each
method. The results are summarized in Table 1 for Normal distribution and
Table 2 for Laplace distribution below.

We can observe that the three SOS-based approaches seem to have comparable
coverage probabilities in most cases, for both generating distributions, in comparison to the EL,
bootstrap, and delta method. However, for small sample situations
($n=20$) EL and all of the SOS-based approaches appear to perform better than the rest. It is discussed in \cite{lam_quantifying_2015} that EL has
better finite sample performance compared to delta method and bootstrap. 
We can also notice that all empirical SOS methods tend to have smaller
variance compared to others, especially for relatively large sample sizes
($n=100,500$). Between the three SOS methods, we can see that explicit
formulations work better compared to implicit, which follows our discussion
after Definition \ref{ExplicitSOSFormulation}. For the two
explicit-formulation methods, since we know the data affects the objective
function in the form $X_{i} ^{(1)}+\ldots+X_{i}^{(l)}$, we would expect better
performance if we combined the data into a single dimension. The numerical
results validate our intuition. \renewcommand{\arraystretch}{0.25}
\begin{table}[th]
\label{TableCVaRGaussian} {\small \ \centering
\begin{tabular}
[c]{l|l|lllll}%
$n$ & Method &
\begin{tabular}
[c]{@{}l}%
Coverage\\
Probability
\end{tabular}
&
\begin{tabular}
[c]{@{}l}%
Mean Lower\\
Bound
\end{tabular}
&
\begin{tabular}
[c]{@{}l}%
Mean Upper\\
Bound
\end{tabular}
&
\begin{tabular}
[c]{@{}l}%
Mean Interval\\
Length
\end{tabular}
&
\begin{tabular}
[c]{@{}l}%
S.D. of\\
Length
\end{tabular}
\\\hline
\multirow{6}{*}{$20$} & ESOS-C & $79.8\%$ & $2.59$ & $4.68$ & $2.09$ &
$0.79$\\
& ESOS-O & $73.4\%$ & $2.55$ & $4.65$ & $2.10$ & $1.21$\\
& ISOS & $70.8\%$ & $2.34$ & $4.87$ & $2.53$ & $0.82$\\
& EL & $71.7\%$ & $2.61$ & $5.18$ & $2.57$ & $1.92$\\
& BT & $55.6\%$ & $1.76$ & $3.88$ & $2.12$ & $1.23$\\
& CLT & $71.8\%$ & $2.01$ & $4.52$ & $2.51$ & $1.87$\\\hline
\multirow{6}{*}{$50$} & ESOS-C & $93.3\%$ & $2.67$ & $4.57$ & $1.90$ &
$0.30$\\
& ESOS-O & $91.0\%$ & $2.63$ & $4.54$ & $1.91$ & $0.57$\\
& ISOS & $87.3\%$ & $2.70$ & $4.75$ & $2.05$ & $0.56$\\
& EL & $89.2\%$ & $2.81$ & $4.78$ & $1.96$ & $0.83$\\
& BT & $82.7\%$ & $2.30$ & $4.25$ & $1.95$ & $0.77$\\
& CLT & $86.6\%$ & $2.47$ & $4.44$ & $1.97$ & $0.78$\\\hline
\multirow{6}{*}{$100$} & ESOS-C & $92.8\%$ & $2.84$ & $4.20$ & $1.36$ & $0.08
$\\
& ESOS-O & $92.4\%$ & $2.80$ & $4.22$ & $1.42$ & $0.23$\\
& ISOS & $91.3\%$ & $2.89$ & $4.32$ & $1.53$ & $0.25$\\
& EL & $91.4\%$ & $2.94$ & $4.46$ & $1.52$ & $0.43$\\
& BT & $90.1\%$ & $2.67$ & $4.16$ & $1.49$ & $0.41$\\
& CLT & $90.4\%$ & $2.75$ & $4.17$ & $1.42$ & $0.39$\\\hline
\multirow{6}{*}{$500$} & ESOS-C & $95.3\%$ & $3.16$ & $3.85$ & $0.69$ & $0.01
$\\
& ESOS-O & $94.9\%$ & $3.14$ & $3.77$ & $0.63$ & $0.05$\\
& ISOS & $91.2\%$ & $3.19$ & $3.88$ & $0.79$ & $0.03$\\
& EL & $93.9\%$ & $3.20$ & $3.93$ & $0.73$ & $0.08$\\
& BT & $94.2\%$ & $3.16$ & $3.84$ & $0.68$ & $0.07$\\
& CLT & $94.7\%$ & $3.17$ & $3.84$ & $0.67$ & $0.08$\\\cline{1-7}%
\end{tabular}
}\caption{\textbf{$\alpha=0.9-$Conditional Value at Risk with Gaussian Data.}
The data $X$ is simulated from 4-dim standard Gaussian distribution, while
each dimension is independent. We consider sample size $n=20,50,100, \text{
and }500$. We repeat the experiments $N=1000$ times and record the coverage
probability for the confidence interval (CI), the average upper and lower
bound for CI, also the average length and standard deviation for CI. ESOS-C is
the explicit formulation of SOS with combined data, ESOS-O stands for
explicit-SOS with original data, ISOS is the implicit SOS, EL stands for
empirical likelihood, BT is short for nonparametric bootstrap, and CLT is the
asymptotic CI method.}%
\end{table}

\renewcommand{\arraystretch}{0.25} \begin{table}[th]
{\small \ \label{TableCVaRLaplace}\centering
\begin{tabular}
[c]{l|l|lllll}%
$n$ & Method &
\begin{tabular}
[c]{@{}l}%
Coverage\\
Probability
\end{tabular}
&
\begin{tabular}
[c]{@{}l}%
Mean Lower\\
Bound
\end{tabular}
&
\begin{tabular}
[c]{@{}l}%
Mean Upper\\
Bound
\end{tabular}
&
\begin{tabular}
[c]{@{}l}%
Mean Interval\\
Length
\end{tabular}
&
\begin{tabular}
[c]{@{}l}%
S.D. of\\
Length
\end{tabular}
\\\hline
\multirow{6}{*}{$20$} & ESOS-C & $78.2\%$ & $3.57$ & $6.89$ & $3.32$ &
$1.10$\\
& ESOS-O & $73.8\%$ & $3.48$ & $7.10$ & $3.62$ & $1.91$\\
& ISOS & $73.1\%$ & $3.87$ & $7.55$ & $3.68$ & $1.16$\\
& EL & $72.3\%$ & $3.56$ & $8.00$ & $4.44$ & $3.30$\\
& BT & $58.1\%$ & $2.40$ & $6.01$ & $3.61$ & $2.40$\\
& CLT & $70.5\%$ & $2.53$ & $6.90$ & $4.37$ & $3.24$\\\hline
\multirow{6}{*}{$50$} & ESOS-C & $89.4\%$ & $3.78$ & $6.64$ & $2.86$ &
$0.42$\\
& ESOS-O & $89.3\%$ & $3.69$ & $6.78$ & $3.09$ & $0.89$\\
& ISOS & $80.1\%$ & $4.21$ & $7.17$ & $2.96$ & $0.63$\\
& EL & $86.2\%$ & $3.89$ & $7.43$ & $3.53$ & $1.66$\\
& BT & $80.5\%$ & $3.15$ & $6.58$ & $3.43$ & $1.54$\\
& CLT & $83.6\%$ & $3.29$ & $6.64$ & $3.35$ & $1.54$\\\hline
\multirow{6}{*}{$100$} & ESOS-C & $91.9\%$ & $3.93$ & $6.22$ & $2.29$ & $0.14
$\\
& ESOS-O & $90.8\%$ & $3.88$ & $6.30$ & $2.42$ & $0.43$\\
& IISOS & $86.6\%$ & $4.30$ & $6.78$ & $2.44$ & $0.36$\\
& EL & $89.9\%$ & $4.10$ & $6.66$ & $2.56$ & $0.86$\\
& BT & $86.2\%$ & $3.71$ & $6.16$ & $2.45$ & $0.81$\\
& CLT & $87.6\%$ & $3.76$ & $6.17$ & $2.41$ & $0.79$\\\hline
\multirow{6}{*}{$500$} & ESOS-C & $94.7\%$ & $4.53$ & $5.62$ & $1.09$ & $0.06
$\\
& ESOS-O & $94.3\%$ & $4.46$ & $5.59$ & $1.13$ & $0.08$\\
& ISOS & $92.1\%$ & $4.43$ & $5.61$ & $1.17$ & $0.13$\\
& EL & $94.0\%$ & $4.53$ & $5.78$ & $1.25$ & $0.18$\\
& BT & $92.2\%$ & $4.46$ & $5.58$ & $1.12$ & $0.16$\\
& CLT & $93.1\%$ & $4.45$ & $5.48$ & $1.13$ & $0.15$\\\cline{1-7}%
\end{tabular}
}
\par
{\small \ }\caption{ \textbf{$\alpha=0.9-$Conditional Value at Risk with
Laplace Data.} The data $X$ is simulated from 4-dim standard Laplace
distribution, while each dimension is independent. We consider sample size
$n=20,50,100, \text{ and }500$. We repeat the experiments $N=1000$ times and
record the coverage probability for the confidence interval (CI), the average
upper and lower bound for CI, also the average length and standard deviation
for CI. ESOS-C is the explicit formulation of SOS with combined data, ESOS-O
stands for explicit-SOS with original data, ISOS is the implicit SOS, EL
stands for empirical likelihood, BT is short for nonparametric bootstrap, and
CLT is the asymptotic CI method.}%
\end{table}

\end{example}

In addition, we report the computational time for our calculation in Table
\ref{TableComputation}. The different formulations of SOS-based methods
share the same computation cost, thus we only report the case for implicit
SOS. We report the average calculating time in seconds with thousands of
experiments, where the experiments are implemented in Python with Scipy
optimizers and our machine is equipped with an Intel i7 3.5Ghz processor and
16GB memory. Our SOS based method requires solving the C-VaR optimization
problem once, then solve the linear programming. The EL based method is
similar, with solving the C-VaR optimization problem once, it then solves a
convex optimization problem. Finally, the bootstrap based method requires
solving the C-VaR optimization repetitively. We can observe that for the
example we consider, our SOS-based method does not face computational challenges
compared with other methods. \begin{table}[th]
\label{TableComputation}
\begin{tabular}
[c]{l|llll}
& 20 & 50 & 100 & 500\\\hline
ISOS & 0.042 & 0.108 & 0.613 & 14.069\\
EL & 0.018 & 0.069 & 0.401 & 7.272\\
BT & 0.099 & 1.038 & 2.085 & 18.023
\end{tabular}
{\small \ }\caption{ \textbf{Computational Cost for Our C-VaR examples.} The
average computational time in seconds for different algorithms with different
sample sizes. }%
\end{table}

\begin{example}
[Semi-supervised Learning]We consider the DRO formulation for Semi-supervised
Learning (SSL) as suggested in \cite{blanchet2017semi}. We
formulate the data-driven DRO problem and compare the results for choosing the
distributional uncertainty size with the above asymptotic results of SOS
function as suggested in Corollary \ref{SOSGeneralExplicit}. We consider the
MiniBooNE data set from UCI machine learning data base \cite{blake1998uci}. We
consider logistic regression as our baseline model and form SSL-DRO
formulation. For each iteration, we randomly split the data into labeled
training set with size $n=30$, unlabeled training set with size $N-n=5000$,
and testing set with size $n=125034$. We compare the choice of the uncertainty
size using 5-fold cross-validation and SoS asymptotic results. We also include
the results for logistic regression and regularized logistic regression as
reference. We report the average training error and testing error as
log-exponential loss and testing accuracy as accurate classification rate. The
mean and standard deviation of the training error, testing error, and testing
accuracy are evaluated via 500 independent experiments. The details are
included in Table \ref{Table-SSL}.

\begin{table}[th]
\centering
\par%
\begin{tabular}
[c]{c|ccc}
& Training Error & Testing Error & Testing Accuracy\\\hline
Logistic Regression & $0 \pm0$ & $18.2\pm10.0$ & $.678 \pm.059$\\
LRL1 with CV & $.401 \pm.167$ & $.910 \pm.131$ & $.717 \pm.041$\\
DRO-SSL with CV & $.287 \pm.047$ & $.609 \pm.054$ & $.710 \pm.032$\\
DRO-SSL with SoS & $.304 \pm. 045$ & $.682 \pm.048$ & $.709 \pm.028$\\\hline
\end{tabular}
\caption{ \textbf{Numerical Results for Semi-supervised Learning.}}

\label{Table-SSL}%
\end{table}
\end{example}

\section{Conclusions and Discussion\label{Section_Conclusions}}

This paper introduces a methodology inspired by Empirical Likelihood, but in
which the likelihood ratio function is replaced by a Wasserstein distance. The
method that we propose is motivated by the problem of systematically finding
estimators that incorporate out-of-sample performance in their design.

In turn, as a motivation for the need to find these types of estimators we
discussed applications to stress testing and semi-supervised learning, which
have been discussed in the body of this paper.  Another way in which we can
justify our framework is as an approximation approach to solving the problem
\[
\min_{\theta\in R^{l}}\max_{P\in\mathcal{U}_{\delta}\left(  P_{n}%
;\mathbb{R}^{d}\right)  }E_{P}[\mathcal{L}\left(  X,\theta\right)  ].
\]
It turns out that in great generality (see \cite{esfahani_data-driven_2015})
\[
\max_{P\in\mathcal{U}_{\delta}\left(  P_{n};\mathbb{R}^{d}\right)  }%
E_{P}[\mathcal{L}\left(  X,\theta\right)  ]=\min_{\lambda\geq0}\{\lambda
\delta+E_{P_{n}}[f\left(  X,\theta;\lambda\right)  ]\},
\]
where $f\left(  x,\theta;\lambda\right)  $ is defined as the solution of an
optimization problem involving a parameter $y\in\mathbb{R}^{d}$ which we refer
to as the \textquotedblleft inner optimization problem.\textquotedblright The
inner optimization problem is typically not convex and therefore it is
challenging to solve. There are cases in which the inner optimization problem
can be solved in closed form, however, and many of those cases have been
documented in the literature in \cite{esfahani_data-driven_2015}. Our results
can be used to suitably calibrate an alternative formulation that may be more
tractable given that $y\in\mathbb{R}^{d}$ is replaced by $y\in\mathcal{Z}%
_{n,m}$.

There are a number of structural properties in our procedure that are worth
investigating and that we plan to explore in future work. For instance, we believe the
choice of a particular cost in optimal transport distance deserves
substantial analysis. In this paper we have chosen the L$_{2}$ Wasserstein
metric to illustrate our results. The methodology that we propose can be
extended to cover other Wasserstein metrics, so on the technical side our work
provides the foundations for such extensions. However, it is the impact of
such selection that appears to also bring about interesting connections. This
already is made evident from our work \cite{blanchet2016robust} in which we
see that the connections that we mentioned earlier in this discussion (to
LASSO and SVM) are made after carefully choosing a natural Wasserstein metric.

In addition, given the parallel philosophy underpinning the method that we
proposed (based on Empirical Likelihood), the results described in this paper open up a
significant amount of research opportunities that are parallel to the
substantial literature produced in the area of Empirical Likelihood during the
last three decades. We mention, in particular, applications to regression
problems (see
\cite{owen_empirical_1991,chen_accuracy_1993,wang_empirical_2001,zhao_empirical_2008,chen_review_2009,murphy_likelihood_1995,li_nonparametric_1996,
hollander_likelihood_1997,li_semiparametric_1997,einmahl_confidence_1999,
wang_likelihood_2009,zhou_empirical_2015}), machine learning (see \cite{duchi2016statistics,hu2018does, duchi2018learning,blanchet2019distributionally}), econometrics (see
\cite{newey_higher_2004,bravo_empirical_2004,kitamura_empirical_2006,
antoine_efficient_2007,guggenberger_finite_2008,imbens_generalized_2012}),
and additional recent work on stochastic optimization (see
\cite{lam_quantifying_2015,lam_empirical_2016,blanchet_quantify_2016}). The
methodology we propose could be extended to the above applications by simply
replacing the Empirical Likelihood function by the SOS function and by
applying asymptotic theorems developed in this paper (or natural extensions).

\section{Methodological Development\label{Section_Methodology}}

We shall analyze the limiting distribution of the SOS profile function for
means first. In order to gain some intuition let us perform some basic
manipulations. First, without loss of generality we assume $\theta_{\ast}=0$,
otherwise, we can let $\tilde{X_{i}}=X_{i}-\theta_{\ast}$ and apply the
analysis to the $\tilde{X_{i}}$'s.

\subsection{The Dual Problem and High-Level Understanding of
Results\label{Sect_Intuitive}}

\paragraph{The Dual Problem}

Let us revisit the definition of (\ref{EWP_1}) and write it as a linear
programming problem,
\begin{align}
R_{n}^{W}(\theta_{\ast})  &  =\min_{ \pi(i,j) \geq0 } \quad\sum_{i=1}^{n}%
\sum_{j=1}^{m+n}\pi(i,j)\left\Vert X_{i}-Z_{j}\right\Vert _{2}^{2}%
\label{Primal_1}\\
&  \text{s.t. }\left\{
\begin{array}
[c]{rcl}%
\sum_{j=1}^{(m+n)}\pi(i,j)= 1/n,\text{ for all }i &  & \\
\sum_{j=1}^{(m+n)} \left(  \sum_{i=1}^{n}\pi(i,j)\right)  Z_{j}=0 &  &
\end{array}
\right.  .\nonumber
\end{align}

We know with probability one when $n\rightarrow\infty$, $\vec{0}$ is in the
convex hull of $Z_{j}$, thus the original linear programming problem is
feasible for all $n$ large enough with probability one. Applying the strong
duality theorem for linear programming problem, see for example,
\cite{luenberger_introduction_1973}, we can write (\ref{Primal_1}) in the dual
formulation as
\begin{align*}
R_{n}^{W}(\theta_{\ast})  &  =\max_{\lambda,\tilde{\gamma}_{i}}\left\{
-\frac{1}{n}\sum_{i=1}^{n}\tilde{\gamma}_{i}\right\} \\
&  s.t.\quad\tilde{\gamma}_{i}+\left\Vert X_{i}-Z_{j}\right\Vert _{2}%
^{2}-\lambda^{T}Z_{j}\geq0\text{ for all }i,j.
\end{align*}
Let us define ${\gamma}_{i}=\tilde{\gamma}_{i}-\lambda^{T}Z_{i}$. By the
constraint in the above optimization problem, if we take $i=j$, we have
$\tilde{\gamma}_{i}\geq\lambda^{T}Z_{i}$, which is equivalent to $\gamma
_{i}\geq0$. Then, we can write the optimization problem in $\gamma_{i}$'s as
\begin{align*}
R_{n}^{W}(\theta_{\ast})  &  =\max_{\lambda,\gamma_{i}\geq0}\left\{
-\lambda^{T}\bar{X}_{n}-\frac{1}{n}\sum_{i=1}^{n}\gamma_{i}\right\} \\
\text{s.t.}\;-\lambda^{T}X_{i}-\gamma_{i}  &  \leq-\lambda^{T}Z_{j}+\left\Vert
X_{i}-Z_{j}\right\Vert _{2}^{2},\text{ for all }i,j.
\end{align*}

We can further simplify the constraints by minimizing over $j$, while keeping
$i$ fixed, therefore arriving to the simplified dual formulation
\begin{align}
R_{n}^{W}(\theta_{\ast})  &  =\max_{\lambda,\gamma_{i}\geq0}\left\{
-\lambda^{T}\bar{X}_{n}-\frac{1}{n}\sum_{i=1}^{n}\gamma_{i}\right\}
\label{WPDualFormulation}\\
\text{s.t.}\;-\lambda^{T}X_{i}-\gamma_{i}  &  \leq\inf_{j}\left\{  -
\lambda^{T}Z_{j}+\left\Vert X_{i}-Z_{j}\right\Vert _{2}^{2}\right\}  ,\text{
for all }i\text{.}\nonumber
\end{align}

\paragraph{\textbf{High-Level Intuitive Analysis}}

At this point we can perform a high-level analysis which can help us guide our
intuition about our result. First, consider an approximation performed by
freeing the $Z_{j}$ in the constraints of (\ref{WPDualFormulation}), in this
portion the reader can appreciate that the assumption that $X_{j}$ has a
density yields%
\begin{equation}
\inf_{j}\left\{  \left\Vert Z_{j}-(X_{i}+\lambda/2)\right\Vert _{2}%
^{2}\right\}  =\epsilon_{n}\left(  i\right)  , \label{Eq_DET}%
\end{equation}
where error $\epsilon_{n}\left(  i\right)  $ is small as $n\rightarrow\infty$
and it will be discussed momentarily. Equation (\ref{Eq_DET}) is equivalent
to
\[
\inf_{j}\left\{  -\lambda^{T}Z_{j}+\left\Vert X_{i}-Z_{j}\right\Vert _{2}%
^{2}\right\}  =-\lambda^{T}X_{i}-\left\Vert \lambda\right\Vert _{2}%
^{2}/4+\epsilon_{n}\left(  i\right)  .
\]
Hence, the $i$-th constraint in (\ref{WPDualFormulation}) takes the form
\[
-\lambda^{T}X_{i}-\gamma_{i}\leq-\lambda^{T}X_{i}-\left\Vert \lambda
\right\Vert _{2}^{2}/4+\epsilon_{n}\left(  i\right)  ,
\]
and thus (\ref{WPDualFormulation}) can ultimately be written as%
\begin{align}
R_{n}^{W}(\theta_{\ast})  &  =-\min_{\lambda,\gamma_{i}\geq0}\left\{
\lambda^{T}\bar{X}_{n}+\frac{1}{n}\sum_{i=1}^{n}\gamma_{i}\right\}
\label{Simp_Dual_P1}\\
\text{s.t.}\;\gamma_{i}  &  \geq\left(  1-\epsilon_{n}\left(  i\right)
\right)  \left\Vert \lambda\right\Vert _{2}^{2}/4\text{ for all }%
i\text{.}\nonumber
\end{align}

Now, observe that if $Z_{j}$ was free, then the optimal choice in
(\ref{Eq_DET}) would be $a_{\ast}\left(  i\right)  =X_{i}+\lambda/2$. \ 

Consider the case $l=1$, in this case it is not difficult to convince
ourselves (because of the existence of a density) that $\epsilon_{n}\left(
i\right)  =O_{p}\left(  1/n\right)  $ as $n\rightarrow\infty$ (basically with
a probability which is bounded away from zero there will be a point in the
sample $\{Z_{1},...,Z_{m+n}\}\backslash X_{i}$ which is within $O_{p}\left(
1/n\right)  $ distance of $a_{\ast}\left(  i\right)  $). Then it is intuitive
to expect the approximation
\[
R_{n}^{W}(\theta_{\ast})=-\min_{\lambda}\left\{  \lambda\bar{X}_{n}+\left(
1+O_{p}\left(  1/n\right)  \right)  \lambda^{2}/4\right\}  ,
\]
which formally yields an optimal selection
\[
\lambda_{\ast}=-\frac{\bar{X}_{n}}{\left(  1/2+O_{p}\left(  1/n\right)
\right)  }=-2\bar{X}_{n}+O_{p}\left(  1/n^{3/2}\right)  ,
\]
and therefore we expect, due to the Central Limit Theorem (CLT), that%
\begin{equation}
nR_{n}^{W}(\theta_{\ast})=n\bar{X}_{n}^{2}+nO_{p}\left(  1/n^{3/2}\right)
\Rightarrow Var\left(  X\right)  \chi_{1}^{2}, \label{Chi_SQ_INT_L}%
\end{equation}
as $n\rightarrow\infty$. This analysis will be made rigorous in the next subsection.

Let us continue our discussion in order to elucidate why the rate of
convergence in the asymptotic distribution of $R_{n}^{W}(\theta_{\ast})$
depends on the dimension. Such dependence arises due to the presence of the
error term $\epsilon_{n}\left(  i\right)  $. Note that in dimension $l=2$, we
expect $\epsilon_{n}\left(  i\right)  =O_{p}\left(  1/n^{1/2}\right)  $; this
time, with positive probability (uniformly as $n\rightarrow\infty$) we must
have that a point in the sample $\{Z_{1},...,Z_{m+n}\}\backslash X_{i}$ is
within $O_{p}\left(  1/n^{1/2}\right)  $ distance of $a_{\ast}\left(
i\right)  $ (because the probability that $X_{i}$ lies inside a ball of size
$1/n^{1/2}$ around a point $a$ is of order $O\left(  1/n^{1/2}\right)  $).
Therefore, in the case $l=2$ we formally have $\lambda_{\ast}\left(  n\right)
=-\bar{X}_{n}+O_{p}\left(  n^{-1/2}\right)  $, but we know from the CLT that
$\bar{X}_{n}=O_{p}\left(  n^{-1/2}\right)  $ so this time contribution of
$\epsilon_{i}\left(  n\right)  $ is non-negligible.

Similarly, when $l\geq3$ this simple analysis allows us to conclude that the
contribution of $\epsilon_{i}\left(  n\right)  =O\left(  n^{-1/l}\right)  $
will actually dominate the behavior of $\lambda_{\ast}\left(  n\right)  $ and
this explains why the rate of convergence depends on the dimension of the
vector $X_{i}$, namely, $l$. The specific rate depends on a delicate analysis
of the error being $\epsilon_{i}\left(  n\right)  $ which is performed in the
next sub-section. A key technical device introduced in our proof technique is
a Poisson point process which approximates the number of points in
$\{Z_{1},...,Z_{m+n}\}\backslash X_{i}$ which are within a distance of size
$O\left(  n^{-1/l}\right)  $ from the free optimizer $a_{\ast}\left(
i\right)  $ arising in (\ref{Eq_DET}).

The introduction of this point process, which in turn is required to analyze
$\epsilon_{i}\left(  n\right)  $, makes the proof of our result substantially
different from the standard approach used in the theory of Empirical
Likelihood (see
\cite{owen_empirical_1988,owen_empirical_1990,qin_empirical_1994}), which
builds on \cite{wilks_large-sample_1938}.

\subsection{Proof of Theorem \ref{SOSTheoremMean}}

The proof of Theorem \ref{SOSTheoremMean} is divided in several steps which we
will carefully record so that we can build from these steps in order to prove
the remaining results in the paper.

\subsubsection{\textbf{Step 1 (Dual Formulation and Lower Bound):}}

Using the same transformations introduced in (\ref{Primal_1}) we can obtain
the dual formulation of the SOS\ profile function (\ref{EWP_1}), which is a
natural adaptation of (\ref{WPDualFormulation}), namely
\begin{align*}
R_{n}^{W}(\theta_{\ast})  &  =\max_{\lambda,\gamma_{i}\geq0}\left\{
-\lambda\bar{X}_{n}-\frac{1}{n}\sum_{i=1}^{n}\gamma_{i}\right\} \\
\text{s.t.}\;-\lambda^{T}X_{i}-\gamma_{i}  &  \leq\inf_{j}\left\{
-\lambda^{T}Z_{j}+\left\Vert X_{i}-Z_{j}\right\Vert _{2}^{2}\right\}  ,\text{
for all }i\text{.}%
\end{align*}

Observe that the following lower bound applies by optimizing over $a\in R^{l}
$ instead of $a=Z_{j}\in\mathcal{Z}_{n}$, therefore obtaining the lower bound
\begin{align*}
\inf_{j}\left\{  -\lambda^{T}Z_{j}+\left\Vert X_{i}-Z_{j}\right\Vert _{2}%
^{2}\right\}   &  \geq\inf_{a}\left\{  -\lambda^{T}a+\left\Vert X_{i}%
-a\right\Vert _{2}^{2}\right\} \\
&  =-\lambda^{T}X_{i}-\left\Vert \lambda\right\Vert _{2}^{2}/4,
\end{align*}
with the optimizer $a_{\ast}\left(  X_{i},\lambda\right)  =X_{i}+\lambda/2$.

\subsubsection{\textbf{Step 2 (Auxiliary Poisson Point Processes):}}

Then, for each $i$ let us define a point process,
\[
N_{n}^{(i)}(t,\lambda)=\#\left\{  Z_{j}:\left\Vert Z_{j}-a_{\ast}\left(
X_{i},\lambda\right)  \right\Vert _{2}^{2}\leq t^{2/l}/n^{2/l},Z_{j}\neq
X_{i}\right\}  ,\label{PoP}%
\]
(recall that $Z_{j}\in R^{l}$). Observe that, actually, we have%
\[
N_{n}^{(i)}(t,\lambda)=N_{n}^{(i)}(t,\lambda,1)+N_{n}^{(i)}(t,\lambda,2),
\]
where
\begin{align*}
N_{n}^{(i)}(t,\lambda,1)  &  =\#\left\{  X_{j}:\left\Vert X_{j}-a_{\ast
}\left(  X_{i},\lambda\right)  \right\Vert _{2}^{2}\leq t^{2/l}/n^{2/l}%
,X_{j}\neq X_{i}\right\}  ,\\
N_{n}^{(i)}(t,\lambda,2)  &  =\#\left\{  Y_{j}:\left\Vert Y_{j}-a_{\ast
}\left(  X_{i},\lambda\right)  \right\Vert _{2}^{2}\leq t^{2/l}/n^{2/l}%
\right\}  .
\end{align*}

For any $X_{j}$ with $j\neq i$, conditional on $X_{i}$, due to the assumption
of density and the formula for the volume of $l$-dimensional ball
(\cite{rudin1964principles}), we have,
\begin{align*}
&  \mathbb{P}\left[  \left\Vert X_{j}-a_{\ast}\left(  X_{i},\lambda\right)
\right\Vert _{2}^{2}\leq t^{2/l}/n^{2/l}\middle|X_{i}\right] \\
&  =f_{X}\left(  a_{\ast}\left(  X_{i},\lambda\right)  \right)  \frac
{\pi^{l/2}}{\Gamma(l/2+1)}\frac{t}{n}+o_{p}(t/n)=f_{X}\left(  X_{i}%
+\lambda/2\right)  \frac{\pi^{l/2}}{\Gamma(l/2+1)}\frac{t}{n}+o_{p}(t/n).
\end{align*}
Similarly,
\[
\mathbb{P}\left[  \left\Vert Y_{j}-a_{\ast}\left(  X_{i},\lambda\right)
\right\Vert _{2}^{2}\leq t^{2/l}/n^{2/l}\middle|X_{i}\right]  =f_{Y}\left(
X_{i}+\lambda/2\right)  \frac{\pi^{l/2}}{\Gamma(l/2+1)}\frac{t}{n}%
+o_{p}(t/n).
\]
Since we have i.i.d. structure for the data points, thus we know, $N_{n}%
^{(i)}(t,\lambda,1)$ and $N_{n}^{(i)}(t,\lambda,2)$ conditional on $X_{i}$
follow binomial distributions,
\begin{align*}
N_{n}^{(i)}(t,\lambda,1)|X_{i}  &  \sim\text{Bin}\left(  f_{X}\left(
X_{i}+\lambda/2\right)  \frac{\pi^{l/2}}{\Gamma(d/2+1)}\frac{t}{n}%
+o_{p}(t/n),n-1\right)  ,\\
N_{n}^{(i)}(t,\lambda,2)|X_{i}  &  \sim\text{Bin}\left(  f_{Y}\left(
X_{i}+\lambda/2\right)  \frac{\pi^{l/2}}{\Gamma(l/2+1)}\frac{t}{n}%
+o_{p}(t/n),\left[  \kappa n\right]  \right)  ,\\
N_{n}^{\left(  i\right)  }\left(  t,\lambda\right)   &  =N_{n}^{(i)}%
(t,\lambda,1)+N_{n}^{(i)}(t,\lambda,2).
\end{align*}

Moreover, we have as $n\rightarrow\infty$,
\[
f_{X}\left(  X_{i}+\lambda/2\right)  \frac{\pi^{l/2}}{\Gamma(l/2+1)}\frac
{t}{n}\times(n-1)\rightarrow f_{X}\left(  X_{i}+\lambda/2\right)  \frac
{\pi^{l/2}}{\Gamma(l/2+1)}t.
\]
Thus, by Poisson approximation to the binomial distribution, we have the weak
convergence result
\[
N_{n}^{(i)}(\cdot,\lambda,1)|X_{i}\Rightarrow\text{Poisson}\left(
f_{X}\left(  X_{i}+\lambda/2\right)  \frac{\pi^{l/2}}{\Gamma(l/2+1)}%
\cdot\right)  ,
\]
in $D[0,\infty)$.

So we have that $N_{n}^{(i)}(\cdot,\lambda,1)$, conditional on $X_{i}$, is
asymptotically a time homogeneous Poisson process with rate $f_{X}\left(
X_{i}+\lambda/2\right)  \pi^{d/2}/\Gamma(d/2+1)$. Similar considerations apply
to $N_{n}^{(i)}(\cdot,\lambda,2)|X_{i}$ which yield that
\[
N_{n}^{(i)}(\cdot,\lambda)|X_{i}\Rightarrow\text{Poi}\left(  \Lambda\left(
X_{i},\lambda\right)  \cdot\right)  ,
\]
where
\[
\Lambda\left(  X_{i},\lambda\right)  =\left[  f_{X}\left(  X_{i}%
+\lambda/2\right)  +\kappa f_{Y}\left(  X_{i}+\lambda/2\right)  \right]
\frac{\pi^{l/2}}{\Gamma(l/2+1)}.
\]

Let us write $T_{i}\left(  n,\lambda\right)  $ to denote the first arrival time of
$N_{n}^{(i)}(\cdot,\lambda)$, that is,
\[
T_{i}\left(  n,\lambda\right)  =\inf\left\{  t\geq0:N_{n}^{(i)}(t,\lambda
)\geq1\right\}
\]
Then, we can specify the survival function for $T_{i}\left(  n\right)  $ to
be:
\begin{equation}
\mathbb{P}\left[  T_{i}\left(  n,\lambda\right)  >t\text{ }\middle|\text{ }%
X_{i}\right]  =\mathbb{P}\left[  N_{n}^{(i)}(t,\lambda)=0\text{ }%
\middle|\text{ }X_{i}\right]  =\exp\left(  -\Lambda\left(  X_{i}%
,\lambda\right)  t\right)  \left(  1+O\left(  1/n^{1/l}\right)  \right)  ,
\label{Rate}%
\end{equation}
uniformly on $t$ over compact sets. The error rate $O\left(  1/n^{1/l}\right)
$ is obtained by a simple Taylor expansion of the exponential function applied
to the middle term in the previous string of equalities. Motivated by the form
in the right hand side of (\ref{Rate}) we define $\tau_{i}\left(
X_{i},\lambda\right)  $ to be a random variable such that
\[
\mathbb{P}\left[  \tau_{i}\left(  X_{i},\lambda\right)  >t|X_{i}\right]  =\exp\left(
-\Lambda\left(  X_{i},\lambda\right)  t\right)  ,
\]
and we drop the dependence on $X_{i}$ and the subindex $i$ when we refer to
the unconditional version of $\tau_{i}\left(  X_{i},\lambda\right)  $, namely%
\[
\mathbb{P}\left[  \tau \left(\lambda\right)>t\right]  =\mathbb{E}\left[  \exp\left(  -\Lambda
\left(  X_{1},\lambda\right)  t\right)  \right]  .
\]

We finish Step 2 with the statement of two technical lemmas. The first
provides a rate of convergence for the Glivenko-Cantelli theorem associated to
the sequence $\{T_{i}\left(  n,\lambda\right)  \}_{i=1}^{n}$.

\begin{lemma}
\label{Lem_Tech}For any $T\in\left(  0,\infty\right)  $ (deterministic) and
$\alpha\in(0,2]$, we have that\newline%
\[
\overline{\lim}_{n\rightarrow\infty}\mathbb{E}\left(  \sup_{t\in\lbrack
0,T]}\left\vert \frac{1}{n^{1/2}}\sum_{i=1}^{n}\left(  I\left(  T_{i}\left(
n,\lambda\right)  \leq t\right)  -\mathbb{P[}T_{i}\left(  n,\lambda\right)  \leq t]\right)
\right\vert \right)  <\infty,
\]
and
\[
\overline{\lim}_{n\rightarrow\infty}\mathbb{E}\left(  \sup_{t\in\lbrack
0,T]}\left\vert \frac{1}{n^{1/2}}\sum_{i=1}^{n}\left(  \max\left(  t^{2}%
-T_{i}(n,\lambda)^{\alpha},0\right)  -\mathbb{E}\left[  \max\left(  t^{2}%
-T_{i}(n,\lambda)^{\alpha},0\right)  \right]  \right)  \right\vert \right)  <\infty.
\]

\end{lemma}

The second technical lemma deals with local properties of the distribution of
$T_{i}\left(  n,\lambda\right)  $. The proofs of both of these technical results are
given at the end of the proof of Theorem \ref{SOSTheoremMean}, in Section
\ref{Subsection_Technical_Lemmas}.

\begin{lemma}
\label{Lem_Tech2} For $X_{i}\in\mathbb{R}^{l}$ and any finite $t$, we have the
Poisson approximation to binomial as:
\[
\mathbb{P}\left[  T_{i}\left(  n,\lambda\right)  \leq t\right]  -\mathbb{P}\left[
\tau(\lambda)\leq t\right]  =O(t^{1+1/l}/{n^{1/l}}),
\]
and
\[
\mathbb{P}\left[  T_{i}\left(  n,\lambda\right)  \leq t\right]  -\mathbb{P}\left[
\tau(\lambda)\leq t\right]  = \mathbb{P}\left[  \tau> t\right]  O\left(  1/n^{l}%
\right)  .
\]

\end{lemma}

\subsubsection{\textbf{Step 3 (Closest Point and SOS Function
Simplification):}}

Note that the $i$-th constraint, namely,
\[
-\gamma_{i}\leq\lambda^{T}X_{i}+\inf_{j}\left\{  -\lambda^{T}Z_{j}+\left\Vert
X_{i}-Z_{j}\right\Vert _{2}^{2}\right\}  ,
\]
can be written as%
\begin{align*}
-\gamma_{i}  &  \leq\inf_{j}\left\{  -\lambda^{T}(Z_{j}-X_{i})+\left\Vert
X_{i}-Z_{j}\right\Vert _{2}^{2}\right\} \\
&  =-\left\Vert \lambda\right\Vert _{2}^{2}/4+\inf_{j}\left\{  \left\Vert
Z_{j}-(\lambda/2+X_{i})\right\Vert _{2}^{2}\right\} \\
&  =-\left\Vert \lambda\right\Vert _{2}^{2}/4+T_{i}^{2/l}\left(  n,\lambda\right)
/n^{2/l}.
\end{align*}
However, since $\gamma_{i}\geq0$ we must have that
\[
-\gamma_{i}\leq-\left\Vert \lambda\right\Vert _{2}^{2}/4+\min\left(
T_{i}^{2/l}\left(  n,\lambda\right)  /n^{2/l},\left\Vert \lambda\right\Vert _{2}%
^{2}/4\right)  .
\]
Therefore, the SOS profile function takes the form%
\[
R_{n}^{W}(\theta_{\ast})=\max_{\lambda}\left\{  -\lambda^{T}\bar{X}%
_{n}-\left\Vert \lambda\right\Vert _{2}^{2}/4+\frac{1}{n}\sum_{i=1}^{n}%
\min\left(  \frac{T_{i}^{2/l}\left(  n,\lambda\right)  }{n^{2/l}},\left\Vert
\lambda\right\Vert _{2}^{2}/4\right)  \right\}  .
\]

To simplify the notation, let us redefine $\lambda\longleftarrow2\lambda$ then
we have that the simplified SOS profile function becomes:
\begin{equation}
\label{SOS_EASY}R_{n}^{W}(\theta_{\ast})=\max_{\lambda}\left\{  -2\lambda
^{T}\bar{X}_{n}-\frac{1}{n}\sum_{i=1}^{n}\max\left(  \left\Vert \lambda
\right\Vert _{2}^{2}- \frac{T_{i}^{2/l}\left(  n,\lambda\right)  }{n^{2/l}},0\right)
\right\}  .
\end{equation}

\subsubsection{\textbf{Step 4 (Case }$l=1$\textbf{):}}

When $l=1$, let us denote $\sqrt{n}\bar{X}_{n}=Z_{n}$ and $\sqrt{n}%
\lambda=\zeta$, where by CLT we can show $Z_{n}\Rightarrow\tilde{Z} \sim
N(0,\sigma^{2})$, where when $l=1$ we have $\sigma^{2} = \Sigma$. Then, as
$n\to\infty$, we have:
\begin{align*}
nR_{n}^{W}(\theta_{\ast})  &  =\max_{\zeta}\left\{  -2\zeta Z_{n}-\frac{1}%
{n}\sum_{i=1}^{n}\max\left(  \zeta^{2}-T_{i}^{2}\left(  n,\zeta/\sqrt{n}\right)
n^{-1},0\right)  \right\} \\
&  =\max_{\zeta}\left\{  -2\zeta Z_{n}-\mathbb{E}\left[  \max\left(  \zeta
^{2}-T_{i}^{2}\left(  n,\zeta/\sqrt{n}\right)  n^{-1},0\right)  \right]  \right\}  +o_{p}(1)\\
& =  \max_{\zeta}\left\{  -2\zeta Z_{n}-\mathbb{E}\left[  \max\left(  \zeta
^{2}-T_{i}^{2}\left(  n,0\right)  n^{-1},0\right)  \right]  \right\}  +o_{p}(1)
\end{align*}
The second equation follows the estimate in (Lemma \ref{Lem_Tech}). 
Using the bonded derivative for the density condition and first order Taylor expansion, we can 
prove that 
$\mathbb{E}\left[T_{i}^{2}\left(  n,0\right)\right] - \mathbb{E}\left[T_{i}^{2}\left(  n,\zeta/\sqrt{n}\right)\right] \to 0$ as $n\to\infty$ for any fixed $\zeta$. Since max function is Lipschitz 
continuous function with constant $1$, and using the Dominating Convergence Theorem, the third equation above 
could be derived as
\begin{align}
&\mathbb{E}\left[  \max\left(  \zeta
^{2}-T_{i}^{2}\left(  n,0\right)  n^{-1},0\right)  \right] - \mathbb{E}\left[  \max\left(  \zeta
^{2}-T_{i}^{2}\left(  n,\zeta/\sqrt{n}\right)  n^{-1},0\right)  \right] \nonumber\\
\leq & 
\mathbb{E}\left[ \left\vert T_{i}^{2}\left(  n,\zeta/\sqrt{n}\right)  n^{-1} -T_{i}^{2}\left(  n,0\right) n^{-1}  \right\vert \right] =o_p(1).\label{EQN_DCT_T_i_n}
\end{align}

We know
the objective function as a function of $\zeta$ is a strictly convex function.
Since as $\zeta=b\left\vert Z_{n}\right\vert $ with $b\rightarrow\pm\infty$
implies that the objective function will tend to $-\infty$, we conclude that
the sequence of global optimizers is compact and each optimizer (i.e. for each
$n$) could be characterized by the first order optimality condition almost
surely. To make the analysis more clear, let us denote the expectation in the
maximization problem to be $g\left(  \zeta,n\right)  $, as a function of
$\zeta$, i.e.
\[
G\left(  \zeta,n\right)  = \mathbb{E}\left[  \max\left(  \zeta^{2}-T_{i}%
^{2}\left(  n,0\right)  n^{-1},0\right)  \right]  ,
\]
which is a deterministic function of $\zeta$ and for any $n$ it is convex.
Moreover, the derivative of $G\left(  \zeta,n\right)  $ is,
\[
g\left(  \zeta,n\right)  =\nabla_{\zeta} G\left(  \zeta,n\right)  =
2\zeta\mathbb{P}\left(  T_{i}^2\left(  n,0\right)  \leq n \zeta^{2}\right)  .
\]
We need to notice that while taking the derivative we require exchanging the
derivative and expectation, this can be done true hereby the Dominated
Convergence Theorem since
\[
\delta^{-1}\left\vert \max\left(  \left(  \zeta+\delta\right)  ^{2}-T_{i}
^{2}\left(  n,0\right)  n^{-1},0\right)  -\max\left(  \zeta^{2}-T_{i}^{2}\left(
n,0\right)  n^{-1},0\right)  \right\vert \leq2|\zeta|,
\]
for all $\delta>0$. We can take the derivative with respect to $\zeta$ in
$-2\zeta Z_{n} - G\left(  \zeta,n\right)  $ and set it to zero, as $n\to
\infty$ we obtain
\[
Z_{n}=-\zeta P\left(  T_{i}^2(n,0)\leq n\zeta^{2}\right)  =-\zeta P\left(
\tau^2(0)\leq n\zeta^{2}\right)  +o_{p}(1)=-\zeta+o_{p}(1).
\]
This estimate follows the second result of Lemma \ref{Lem_Tech2}. Therefore,
the optimizer $\zeta_{n}^{\ast}$, satisfies $\zeta_{n}^{\ast}=-Z_{n}+o_{p}(1)
$, as $n\to\infty$. Then, we plug it into the objective function to obtain
that the scaled SOS profile function satisfies
\[
nR_{n}^{W}(\theta_{\ast}) = 2 Z_{n}^{2} - G\left(  Z_{n},n\right)
+o_{p}\left(  1\right)  \text{ as }n\to\infty.
\]
We should notice $G\left(  Z_{n},n\right)  $ is a function defined via
expectation and evaluated at $Z_{n}$, thus it is a random variable that depends on
$Z_{n}$. By definition and $E\left[  \left\vert X\right\vert \right]
=\int_{0}^{\infty}\mathbb{P}\left[  \left\vert X\right\vert \geq t\right]
dt$, we know as $n\to\infty$,
\begin{align*}
G\left(  \zeta,n\right)   &  = \int_{0}^{\zeta^{2}} \mathbb{P}\left[
T_{i}^{2}\left(  n,0\right)  \leq n\left(  \zeta^{2} -t\right)  \right]  dt\\
&  = \int_{0}^{\zeta^{2}} \mathbb{P}\left[  \tau^{2}\left(  0\right)  \leq
n\left(  \zeta^{2} -t\right)  \right]  dt +o(1)\\
&  = \int_{0}^{\zeta^{2}} 1dt +o(1) = \zeta^{2}+o(1),
\end{align*}
where the second equality is derived from the second argument of Lemma
\ref{Lem_Tech2}. Then for the SOS profile function, it becomes,
\[
nR_{n}^{W}(\theta_{\ast}) =2Z_{n}^{2} - Z_{n}^{2}+o_{p}(1)= Z_{n}^{2}%
+o_{p}(1)\text{ as }n\to\infty.
\]
Applying the continuous mapping theorem and the Central Limit Theorem for
$Z_{n}$, we have
\[
nR_{n}^{W}(\theta_{\ast})\Rightarrow\sigma^{2}\chi_{1}^{2}.
\]

\subsubsection{\textbf{Step 5 (Case }$l=2$\textbf{):}}

Once again we introduce the substitution $\zeta=\sqrt{n}\lambda$ and $\sqrt
{n}\bar{X}_{n}=Z_{n}$ into (\ref{SOS_EASY}). Then, scaling the profile
function by $n$, as $n\to\infty$ we have
\begin{align}
nR_{n}^{W}(\theta_{\ast})=  &  \max_{\zeta}\left\{  -2\zeta^{T}Z_{n}-\frac
{1}{n}\sum_{i=1}^{n}\max\left(  \left\Vert \zeta\right\Vert _{2}^{2}%
-T_{i}\left(  n,\zeta/\sqrt{n}\right)  ,0\right)  \right\} \nonumber\\
=  &  \max_{\zeta}\left\{  -2\zeta^{T}Z_{n}-\mathbb{E}\left[  \max\left(
\left\Vert \zeta\right\Vert _{2}^{2}-T_{i}\left(  n,\zeta/\sqrt{n}\right)  ,0\right)
\right]  \right\}  +o_{p}(1)\nonumber\\
=  &  \max_{\zeta}\left\{  -2\zeta^{T}Z_{n}-\mathbb{E}\left[  \max\left(
\left\Vert \zeta\right\Vert _{2}^{2}-T_{i}\left(  n,0\right)  ,0\right)
\right]  \right\}  +o_{p}(1), \label{OP_l2}%
\end{align}
where the second equality is by applying Lemma \ref{Lem_Tech} (the
error is obtained by localizing $\zeta$ on a compact set, which is valid
because the sequence of global optimizers is easily seen to be tight), and the 
third equality is applying similar derivation as in \eqref{EQN_DCT_T_i_n},The
objective function is strictly convex as a function of $\zeta$ and we know
when $\left\Vert \zeta\right\Vert _{2}\rightarrow\infty$ the objective
function tends to $-\infty$, thus each global maximizer (for each $n$) can be
characterized by the first order optimality condition almost surely. Similar
as Case $l=1$, let us denote
\[
G\left(  \zeta,n\right)  = \mathbb{E}\left[  \max\left(  \left\Vert
\zeta\right\Vert _{2}^{2}-T_{i}\left(  n,0\right)  ,0\right)  \right]  .
\]
It is a continuous differentiable and convex function in $\zeta$ and with
derivative equals
\[
g\left(  \zeta,n\right)  = \nabla_{\zeta} G\left(  \zeta,n\right)  =
2\zeta\mathbb{P}\left[  \left\Vert \zeta\right\Vert _{2}^{2}\geq T_{i}\left(
n,0\right)  \right]  = 2\zeta\mathbb{P}\left[  \left\Vert \zeta\right\Vert _{2}
^{2}\geq\tau(0)\right]  +o(1) \text{ as }n\to\infty,
\]
where the first equality requires applying the Dominated Convergence Theorem,
as in the case $l=1$ and the second estimate follows the first argument in
Lemma \ref{Lem_Tech2}.
Combining the above estimation, we have the first order optimality condition
becomes
\begin{equation}
Z_{n}=-\zeta\mathbb{P}\left[  \left\Vert \zeta\right\Vert _{2}^{2}\geq
\tau(0)\right]  +o_{p}(1) = -\zeta\tilde{g}\left(  \zeta\right)  +o_{p}(1) \text{
as }n\to\infty, \label{2DimOpt}%
\end{equation}
where $\tilde{g}\left(  \zeta(0)\right)  = \mathbb{P}\left[  \left\Vert
\zeta\right\Vert _{2}^{2} \geq\tau\right]  $ is a deterministic function of
$\zeta$. Using equation (\ref{2DimOpt}), we conclude that the optimizer
$\zeta_{n}^{\ast}$, satisfies $\zeta_{n}^{\ast}=$ $-\rho Z_{n}+o_{p}\left(
1\right)  $, for some $\rho$. In turn, plugging in this representation into
equation (\ref{2DimOpt}), as $n\to\infty$ we have
\[
\left\Vert \zeta_{n}^{\ast}\right\Vert _{2}\tilde{g}\left(  \zeta_{n}^{\ast
}\right)  +o_{p}(1)=\left\Vert Z_{n}\right\Vert _{2}.
\]
Sending $n\rightarrow\infty$, we conclude that $\rho$ is the unique solution
to
\begin{equation}
\frac{1}{\rho}= \tilde{g}\left(  \rho\tilde{Z}\right)  . \label{Aux_EU}%
\end{equation}
%
Since the objective function is strictly convex and the above equation is
derived from first order optimality condition, we know the solution exists and
is unique (alternatively we can use the continuity and monotonicity of left
and right hand side of (\ref{Aux_EU}), to argue the existence and uniqueness).
Let us plug in the optimizer back to the objective function and we can see the
scaled SOS profile function becomes
\[
nR_{n}^{W}\left(  \theta_{\ast}\right)  =2\rho\left(  \left\Vert \tilde
{Z}\right\Vert _{2}^{2}\right)  \left\Vert Z_{n}\right\Vert _{2}^{2}-G\left(
\zeta_{n}^{\ast},n\right)  +o_{p}(1).
\]

For a positive random variable $Y$, we have: $\mathbb{E}\left[  Y\right]
=\int_{0}^{\infty}\mathbb{P}\left[  Y\geq t\right]  dt$. Therefore, for
$\zeta$ in a compact set, as $n\to\infty$ we have the following estimate
\begin{align*}
G\left(  \zeta,n\right)   &  = \int_{0}^{\left\Vert \zeta\right\Vert _{2}^{2}}
\mathbb{P}\left[  \left\Vert \zeta\right\Vert _{2}^{2}-T_{i}\left(  n,0\right)
\geq t\right]  dt\\
&  = \int_{0}^{\left\Vert \zeta\right\Vert _{2}^{2}}\mathbb{P}\left[
\left\Vert \zeta\right\Vert _{2}^{2}-\tau(0)\geq t\right]  dt+o(1)\\
&  = \left\Vert \zeta\right\Vert _{2}^{2} \int_{0}^{1} \mathbb{P}\left[  1-
{\tau(0)}/{\left\Vert \zeta\right\Vert _{2}^{2}}\geq s\right]  ds +o(1)\\
&  = \left\Vert \zeta\right\Vert _{2}^{2} \mathbb{E}\left[  \max\left(  1-
{\tau(0)}/{\left\Vert \zeta\right\Vert _{2}^{2}},0\right)  \right]  +o(1)\\
&  = \left\Vert \zeta\right\Vert _{2}^{2}\tilde{\eta}\left(  \zeta\right)
+o(1),
\end{align*}
where we define $\tilde{\eta}\left(  \zeta\right)  = \mathbb{E}\left[
\max\left(  1- {\tau(0)}/{\left\Vert \zeta\right\Vert _{2}^{2}},0\right)
\right]  $ is a deterministic continuous function of $\zeta$.
The second equation follows the first result of Lemma \ref{Lem_Tech2}. Finally
combine $G\left(  \zeta,n\right)  $ and the first term, using the CLT and
continuous mapping theorem, where we denote $Z_{n}\Rightarrow\tilde{Z}\sim
N(0,Var(X))$, as $n\to\infty$ we have:
\begin{align*}
nR_{n}^{W}(\theta_{\ast})  &  = 2\rho\left(  \tilde{Z}\right)  \left\Vert
{Z_{n}}\right\Vert _{2}^{2} - \rho\left(  \tilde{Z}\right)  ^{2} \tilde{\eta
}\left(  Z_{n}\right)  \left\Vert {Z_{n}}\right\Vert _{2}^{2} +o_{p}(1)\\
&  \Rightarrow2\rho\left(  \tilde{Z}\right)  \left\Vert {\tilde{Z}}\right\Vert
_{2}^{2} - \rho\left(  \tilde{Z}\right)  ^{2}\tilde{\eta}\left(  \tilde
{Z}\right)  \left\Vert \tilde{Z}\right\Vert _{2}^{2}.
\end{align*}

\subsubsection{\textbf{Step 6 (Case }$l\geq3$\textbf{):}}

For simplicity, let us write $\sqrt{n}\bar{X}_{n}=Z_{n}$ and $n^{\frac
{3}{2l+2}}\lambda=\zeta$, then as $n\to\infty$ we have%

\begin{align*}
&  n^{1/2+\frac{3}{2l+2}}R_{n}^{W}(\theta_{\ast})\\
&  =\max_{\zeta}\left\{  -2\zeta^{T}Z_{n}-n^{\left(  1/2+\frac{3}{2l+2}%
-\frac{2}{l}\right)  }\frac{1}{n}\sum_{i=1}^{n}\max\left(  \left\Vert
\frac{\zeta}{n^{\left(  \frac{3}{2l+2}-\frac{1}{l}\right)  }}\right\Vert _{2}
^{2}-T_{i}^{2/l}\left(  n,\zeta/n^{\frac{3}{2l+2}}\right)  ,0\right)  \right\} \\
&  =\max_{\zeta}\left\{  -2\zeta^{T}Z_{n}-n^{\left(  1/2+\frac{3}{2l+2}%
-\frac{2}{l}\right)  }\mathbb{E}\left[  \max\left(  \left\Vert \frac{\zeta
}{n^{\left(  \frac{3}{2l+2}-\frac{1}{l}\right)  }}\right\Vert _{2} ^{2}%
-T_{1}^{2/l}\left(  n,\zeta/n^{\frac{3}{2l+2}}\right)  ,0\right)  \right]  \right\}  +o_{p}(1)\\
&  =\max_{\zeta}\left\{  -2\zeta^{T}Z_{n}-n^{\left(  1/2+\frac{3}{2l+2}%
-\frac{2}{l}\right)  }\mathbb{E}\left[  \max\left(  \left\Vert \frac{\zeta
}{n^{\left(  \frac{3}{2l+2}-\frac{1}{l}\right)  }}\right\Vert _{2} ^{2}%
-T_{1}^{2/l}\left(  n,0\right)  ,0\right)  \right]  \right\}  +o_{p}(1).
\end{align*}
The estimate in second equation the previous display is due to an application of Lemma
\ref{Lem_Tech}, and the third equation follows the similar derivation as in \eqref{EQN_DCT_T_i_n}.
Similar as for the lower dimensional case, let us denote
\[
G\left(  \zeta,n\right)  = n^{\left(  1/2+\frac{3}{2l+2} -\frac{2}{l}\right)
}\mathbb{E}\left[  \max\left(  \left\Vert \frac{\zeta}{n^{\left(  \frac
{3}{2l+2}-\frac{1}{l}\right)  }}\right\Vert _{2} ^{2}-T_{1}^{2/l}\left(
n,0\right)  ,0\right)  \right]  ,
\]
being a deterministic function continuous and differentiable as a function of
$\zeta$. As we discussed for the case $l=2$ case, the objective function is
strictly convex in $\zeta$, the global optimizers are not only tight, but each
optimizer is also characterized by first order optimality conditions almost
surely. We can apply the Dominated Convergence Theorem, as we discussed for
$l=1$ and the gradient of $G\left(  \zeta,n\right)  $ has the following
estimate as $n\to\infty$,
\begin{align*}
g\left(  \zeta,n\right)   &  = \nabla_{\zeta}G\left(  \zeta,n\right)  =
2n^{\left(  1/2+\frac{3}{2l+2} -\frac{2}{l}\right)  }\zeta\mathbb{P}\left[
T_{i}\left(  n,0\right)  \leq\left\Vert \zeta n^{-\left(  \frac{3}{2l+2}%
-\frac{1}{l}\right)  }\right\Vert _{2} ^{l}\right] \\
&  = 2n^{\left(  1/2+\frac{3}{2l+2} -\frac{2}{l}\right)  }\zeta\mathbb{P}%
\left[  \tau(0)  \leq\left\Vert \zeta n^{-\left(  \frac{3}%
{2l+2}-\frac{1}{l}\right)  }\right\Vert _{2} ^{l}\right]  + o(1).
\end{align*}
The second equality estimate is considering $\zeta$ within a compact set and
the derivation follows the first argument in Lemma \ref{Lem_Tech2}. Then the
first order optimality condition for the SOS profile function becomes,
\begin{align*}
Z_{n}  &  = -n^{\left(  1/2+\frac{3}{2l+2} -\frac{2}{l}\right)  }%
\zeta\mathbb{P}\left[  \tau\left(  0\right)  \leq\left\Vert \zeta n^{-\left(
\frac{3}{2l+2}-\frac{1}{l}\right)  }\right\Vert _{2} ^{l}\right]  + o(1)
\text{ as }n\to\infty.
\end{align*}
For notation simplicity, let us define
\[
\kappa_{n} = \zeta n^{-\left(  \frac{3}{2l+2}-\frac{1}{l}\right)  }.
\]
We can observe for $\zeta$ in a compact set, $\left\Vert \zeta n^{-\left(
\frac{3}{2l+2}-\frac{1}{l}\right)  }\right\Vert _{2} ^{l} = \left\Vert
\kappa_{n}\right\Vert _{2}^{l}\rightarrow0$, as $n\to\infty$, then we can
write
\begin{align*}
&  \mathbb{P}\left[  \tau\left(  0\right)\leq\ \left\Vert \kappa_{n}\right\Vert _{2}%
^{l}\right]  =1-\mathbb{P}\left[  \tau\left(  0\right)> \left\Vert \kappa_{n}\right\Vert
_{2}^{l}\right]  =1-\mathbb{E}\left[  \mathbb{P}\left[  \tau\left(  0\right)> \left\Vert
\kappa_{n}\right\Vert _{2}^{l}\text{ }\middle|\text{ }X_{1}\right]  \right] \\
&  =\mathbb{E}\left[  1-\exp\left(  -\frac{\pi^{l/2}\left(  f_{X}\left(
X_{1}\right)  +\kappa f_{Y}\left(  X_{1}\right)  \right)
}{\Gamma(l/2+1)} \left\Vert \kappa_{n}\right\Vert _{2}^{l}\right)  \right] \\
&  =\mathbb{E}\left[  \frac{\pi^{l/2}}{\Gamma(l/2+1)}\left[  f_{X}\left(
X_{1}\right)  +\kappa f_{Y}\left(  X_{1}\right)  \right]  \left\Vert \kappa
_{n}\right\Vert _{2}^{l}\right]  +o_{p}\left(  n^{-(\frac{3l}{2l+2}-1)}\right)
\\
&  =C \left\Vert \kappa_{n}\right\Vert _{2}^{l}+o_{p}\left(  n^{-(\frac
{3l}{2l+2}-1)}\right)  ,
\end{align*}
where we denote
\[
C=\frac{\pi^{l/2}}{\Gamma(l/2+1)}\mathbb{E}\left[  f_{X}\left(  X_{1}\right)
+\kappa f_{Y}\left(  X_{1}\right)  \right]  .
\]
Plug it back into the optimizer, and as $n\to\infty$ we have:
\[
Z_{n}=-Cn^{(1/2-\frac{3}{2l+2})}n^{\left(  -\frac{3l}{2l+2}+1\right)  }%
\zeta\left\Vert \zeta\right\Vert _{2}^{l}+o_{p}(1)=-C\zeta\left\Vert
\zeta\right\Vert _{2}^{l}+o_{p}(1).
\]
We know that within the objective function, the second term is only based on
the $L_{2}$ norm of $\zeta$, thus to maximize the objective function we will
asymptotically select $\zeta_{n}^{\ast}=-c_{\ast}Z_{n}\left(  1+o\left(
1\right)  \right)  $, where $c_{\ast}>0$ is suitably chosen, thus, we conclude
that the optimizer takes the form,%
\[
\zeta_{n}^{\ast}=-Z_{n}{\left\Vert Z_{n}\right\Vert _{2}^{\left(  \frac
{1}{l+1}-1\right)  }}/{C^{\frac{1}{l+1}}}+o_{p}(1).
\]
Plugging-in the optimizer back into the objective function, as $n\to\infty$we
have:
\[
n^{1/2+\frac{3}{2l+2}}R_{n}^{W}(\theta_{\ast})=-2\zeta_{n}^{\ast\text{ }%
T}Z_{n}- G\left(  \zeta_{n}^{\ast},n\right)  +o_{p}(1).
\]

Let us focus on the analysis of $G\left(  \zeta,n\right)  $ in a compact set.
By definition, we can notice that inside the previous expectation there is a
positive random variable bounded by $\left\Vert \frac{\zeta}{n^{\left(
\frac{3}{2l+2}-\frac{1}{l}\right)  }}\right\Vert _{2}^{2}=\left\Vert
\kappa_{n}\right\Vert _{2}^{2}$, thus as $n\rightarrow\infty$ we have the
following estimate for the expectation as.
\begin{align*}
&  \mathbb{E}\left[  \max\left(  \left\Vert \kappa_{n}\right\Vert _{2}%
^{2}-T_{1}^{2/l}\left(  n,0\right)  ,0\right)  \right]  =\mathbb{E}\left[
\mathbb{E}\left[  \max\left(  \left\Vert \kappa_{n}\right\Vert _{2}^{2}%
-T_{1}^{2/l}\left(  n,0\right)  ,0\right)  \text{ }\middle|\text{ }X_{1}\right]
\right] \\
&  =\mathbb{E}\left[  \int_{0}^{\kappa_{n}}\mathbb{P}\left[  T_{1}\left(
n,0\right)  \leq\left(  \kappa_{n}-t\right)  ^{l/2}\text{ }\middle|\text{ }%
X_{1}\right]  dt\right] \\
&  =\mathbb{E}\left[  \int_{0}^{\left\Vert \kappa_{n}\right\Vert _{2}^{2}%
}\mathbb{P}\left[  \tau\left(  0\right)\leq\left(  \left\Vert \kappa_{n}\right\Vert _{2}%
^{2}-t\right)  ^{l/2}\text{ }\middle|\text{ }X_{1}\right]  +O\left(
1/n^{-1/2+1/l}\right)  dt\right] \\
&  =\mathbb{E}\left[  \int_{0}^{\left\Vert \kappa_{n}\right\Vert _{2}^{2}%
}\left(  1-\exp\left(  -\frac{\pi^{l/2}\left(  f_{X}\left(  X_{1}\right)  +\kappa f_{Y}\left(  X_{1}\right)  \right)  }{\Gamma(l/2+1)}\left(  \left\Vert \kappa
_{n}\right\Vert _{2}^{2}-t\right)  ^{l/2}\right)  \right)  dt\right] \\
&  \quad\quad\quad+O\left(  1/n^{-1/2+3/l-\frac{6}{2l+2}}\right) \\
&  =C\frac{2}{l+2}\left\Vert \frac{\zeta}{n^{\left(  \frac{3}{2l+2}-\frac
{1}{l}\right)  }}\right\Vert ^{l+2}+O\left(  1/n^{-1/2+3/l-\frac{6}{2l+2}%
}\right)
\end{align*}
The estimate in third equation follows by applying the first argument in Lemma
\ref{Lem_Tech2}. The final equality estimate is due to $\left\Vert \kappa
_{n}\right\Vert _{2}^{2}=\left\Vert \zeta n^{-\left(  \frac{3}{2l+2}-\frac
{1}{l}\right)  }\right\Vert _{2}^{2}\rightarrow0$ as $n\rightarrow\infty$.
Then, owing to the previous results, as $n\rightarrow\infty$ we have estimate
for $G\left(  \zeta,n\right)  $ as
\begin{align*}
G\left(  \zeta,n\right)   &  =-\frac{2C}{l+2}n^{\left(  1/2+\frac{3}%
{2l+2}-\frac{2}{l}\right)  }n^{\left(  -\frac{3l+6}{2l+2}+\frac{l+2}%
{l}\right)  }\left\Vert \zeta\right\Vert _{2}^{l+2}+o(1)\\
&  =-\frac{2C}{l+2}\left\Vert \zeta\right\Vert _{2}^{l+2}+o(1).
\end{align*}
Finally, we can know that, as $n\rightarrow\infty$, by the CLT we have
$Z_{n}\Rightarrow\tilde{Z}$, then using continuous mapping theorem, we have
that the scaled SOS profile function has the asymptotic distribution given by
\begin{align*}
&  \quad n^{1/2+\frac{5}{4l+2}}R_{n}^{W}(\theta_{\ast})\\
&  =2\left\Vert Z_{n}\right\Vert _{2}^{2}\frac{\left\Vert Z_{n}\right\Vert
_{2}^{\left(  \frac{1}{l+1}-1\right)  }}{C^{\frac{1}{l+1}}}-\frac{2}{l+2}%
\frac{\left\Vert Z_{n}\right\Vert _{2}^{1+\frac{1}{l+1}}}{C^{\frac{1}{l+1}}%
}+o_{p}(1)\\
&  =\frac{2l+2}{l+2}\frac{\left\Vert Z_{n}\right\Vert _{2}^{1+\frac{1}{l+1}}%
}{C^{\frac{1}{l+1}}}+o_{p}(1)\Rightarrow\frac{2l+2}{l+2}\frac{\left\Vert
\tilde{Z}\right\Vert _{2}^{1+\frac{1}{l+1}}}{C^{\frac{1}{l+1}}}.
\end{align*}

\subsubsection{\textbf{Proofs of Technical Lemmas in Step 2
\label{Subsection_Technical_Lemmas}}}

\begin{proof}
[\textbf{Proof of Lemma \ref{Lem_Tech}}]We shall introduce some notation which
will be convenient throughout our development. Define for $t\geq0$,
\begin{align*}
F_{n}\left(  t\right)   &  =P\left(  T_{i}\left(  n,\lambda\right)  \leq t\right)  ,\\
D_{i}\left(  t\right)   &  =I\left(  T_{i}\left(  n,\lambda\right)  \leq t\right)
,\text{ \ }\bar{D}_{i}\left(  t\right)  =I\left(  T_{i}\left(  n,\lambda\right)  \leq
t\right)  -F_{n}\left(  t\right)  ,\\
\bar{F}_{n}\left(  t\right)   &  =1+n^{-1/2}\sum_{i=1}^{n}\bar{D}_{i}\left(
t\right)  .
\end{align*}

Therefore, we are interested in studying
\[
\bar{F}_{n}\left(  t\right)  -1 =\frac{1}{n^{1/2}}\sum_{i=1}^{n}\left(
I\left(  T_{i}\left(  n,\lambda\right)  \leq t\right)  -F_{n}\left(  t\right)
\right)  .
\]
We will start by studying
\[
\mathbb{E}[\sup\{\bar{F}_{n}\left(  t\right)  :t\in\lbrack0,T]\}].
\]
First, we define
\[
h_{n}\left(  t\right)  =\frac{\bar{F}_{n}\left(  t_{-}\right)  }{\left(
\bar{F}_{n}^{\ast}\left(  t_{-}\right)  ^{2}+[\bar{F}_{n}]\left(
t_{-}\right)  \right)  ^{1/2}},
\]
where, for a given function $\{g\left(  t\right)  :t\in\lbrack0,T]\}$, we
define
\begin{align*}
g^{\ast}\left(  t\right)   &  =\sup\{g\left(  s\right)  :s\in\lbrack0,t]\},\\
\lbrack g]\left(  t\right)   &  =\int_{0}^{t}\left(  dg\left(  s\right)
\right)  ^{2}.
\end{align*}
In addition, $\lbrack g]\left(  t\right)  $ is defined as the quadratic
variational process, i.e.,
\[
\lbrack g]\left(  t\right)  = \lim_{n\to\infty} \sum_{i=1}^{n}\left[  g\left(
\frac{i\times t}{n}\right)  - g\left(  \frac{(i-1)\times t}{n}\right)
\right]  ^{2}.
\]
In particular,
\[
\lbrack\bar{F}_{n}]\left(  t\right)  =\frac{1}{n}\sum_{i=1}^{n}I\left(
T_{i}\left(  n,\lambda\right)  \leq t\right)  .
\]

We observe that $\bar{F}_{n}^{\ast}\left(  t\right)  \geq1$ , therefore
$h_{n}\left(  t\right)  $ is well defined; moreover, note that
\[
h_{n}\left(  t\right)  ^{2}\leq1.
\]
We invoke Theorem 1.2 of \cite{beiglbock_pathwise_2015} and conclude that
\[
\sup_{0\leq t\leq T}\bar{F}_{n}\left(  t\right)  \leq6\sqrt{[\bar{F}%
_{n}]\left(  T\right)  }+2\int_{0}^{T}h_{n}\left(  t\right)  d\bar{F}%
_{n}\left(  t\right)  .
\]
Now we analyze the integral in the right hand side of the previous display.
Observe that
\begin{align}
\mathbb{E}\left(  \int_{0}^{T}h_{n}\left(  t\right)  d\bar{F}_{n}\left(
t\right)  \right)   &  =\frac{1}{n^{1/2}}\sum_{i=1}^{n}\mathbb{E}\left(
\int_{0}^{T}h_{n}\left(  t\right)  d\bar{D}_{i}\left(  t\right)  \right)
\nonumber\\
&  =n^{1/2}\mathbb{E}\left(  \int_{0}^{T}h_{n}\left(  t\right)  d\bar{D}%
_{1}\left(  t\right)  \right)  . \label{Eq_UE_1}%
\end{align}
Let us write
\[
_{1}\bar{F}_{n}\left(  t\right)  =\bar{F}_{n}\left(  t\right)  -\bar{D}%
_{1}\left(  t\right)  /n^{1/2},
\]
that is, we simply remove the last term in the sum defining $\bar{F}%
_{n}\left(  t\right)  $. We have that
\[
h_{n}\left(  t\right)  =\frac{_{1}\bar{F}_{n}\left(  t_{-}\right)  +\bar
{D}_{1}\left(  t_{-}\right)  /n^{1/2}}{\left(  \bar{F}_{n}^{\ast}\left(
t_{-}\right)  ^{2}+[_{1}\bar{F}_{n}]\left(  t_{-}\right)  +\left[
D_{1}\right]  \left(  t_{-}\right)  /n\right)  ^{1/2}},
\]
moreover,
\[
\left\vert _{1}\bar{F}_{n}^{\ast}\left(  t\right)  -\bar{F}_{n}^{\ast}\left(
t\right)  \right\vert \leq1/n^{1/2}.
\]
We then can write
\begin{align}
h_{n}\left(  t\right)   &  =\frac{_{1}\bar{F}_{n}\left(  t_{-}\right)
+\bar{D}_{1}\left(  t_{-}\right)  /n^{1/2}}{\left(  \bar{F}_{n}^{\ast}\left(
t_{-}\right)  ^{2}+[_{1}\bar{F}_{n}]\left(  t_{-}\right)  +\left[
D_{1}\right]  \left(  t_{-}\right)  /n\right)  ^{1/2}}\label{Exp_h_aux}\\
&  =\frac{_{1}\bar{F}_{n}\left(  t_{-}\right)  }{\left(  _{1}\bar{F}_{n}%
^{\ast}\left(  t_{-}\right)  ^{2}+[_{1}\bar{F}_{n}]\left(  t_{-}\right)
\right)  ^{1/2}}\left(  1+\frac{L_{n}\left(  t_{-}\right)  }{n^{1/2}}\right)
,\nonumber
\end{align}
where we can select a deterministic constant $c\in\left(  0,\infty\right)  $
such that $\left\vert L_{n}\left(  t\right)  \right\vert \leq c$ for
$j=0\text{ and }1$ assuming $n\geq4$ (this constrain in $n$ is imposed so that
a Taylor expansion for the function $1/(1-x)$ can be developed for
$x\in\left(  0,1\right)  $). We now insert (\ref{Exp_h_aux}) into
(\ref{Eq_UE_1}) and conclude that if we define
\[
\bar{h}_{n}\left(  t\right)  =\frac{_{1}\bar{F}_{n}\left(  t_{-}\right)
}{\left(  _{1}\bar{F}_{n}^{\ast}\left(  t_{-}\right)  ^{2}+[_{1}\bar{F}%
_{n}]\left(  t_{-}\right)  \right)  ^{1/2}},
\]
it suffices to verify that
\[
n^{1/2}\mathbb{E}\left(  \int_{0}^{T}\bar{h}_{n}\left(  t\right)  d\bar{D}%
_{1}\left(  t\right)  \right)  <\infty.
\]
Define $\widetilde{h}_{n}\left(  t\right)  $ to be a copy of $\bar{h}%
_{n}\left(  t\right)  $, independent of $X_{1}$ and $T_{1}\left(  n\right)  $.
In particular, $\widetilde{h}_{n}\left(  t\right)  $ is constructed by using
all of the $X_{j}$'s except for $X_{1}$, which might be replaced by an
independent copy, $X_{1}^{\prime}$, of $X_{1}$. Observe that the number of
processes $\{\bar{D}_{i}\left(  t\right)  :t\leq T\}$ that depend on
$T_{1}\left(  n\right)  $ and $X_{1}$ is smaller than $N_{n}\left(
T,\lambda,1\right)  $. Therefore, similarly as we obtained from the analysis
leading to the definition of $\bar{h}_{n}\left(  \cdot\right)  $, we have that
a random variable $\bar{L}_{N_{n}\left(  T,\lambda,1\right)  }$ can be defined
so that $\left\vert \bar{L}_{N_{n}\left(  T,\lambda,1\right)  }\right\vert
\leq c(1+N_{n}\left(  T,\lambda,1\right)  )$ for some (deterministic) $c>0$
and $n\geq4$ and satisfying
\begin{align*}
&  \mathbb{E}\left(  \int_{0}^{T}\bar{h}_{n}\left(  t\right)  d\bar{D}%
_{1}\left(  t\right)  \right) \\
&  =\mathbb{E}\left(  \bar{h}_{n}\left(  T_{1}\left(  n\right)  \right)
I\left(  T_{1}\left(  n\right)  \leq T\right)  \right)  -\mathbb{E}\left(
\widetilde{h}_{n}\left(  T_{1}\left(  n\right)  \right)  I\left(  T_{1}\left(
n\right)  \leq T\right)  \right) \\
&  =\mathbb{E}\left(  \widetilde{h}_{n}\left(  T_{1}\left(  n\right)  \right)
I\left(  T_{1}\left(  n\right)  \leq T\right)  \right)  -\mathbb{E}\left(
\widetilde{h}_{n}\left(  \tau_{i}\left(  X_{i}\right)  \right)  I\left(
\tau_{i}\left(  X_{i}\right)  \leq T\right)  \right) \\
&  +\mathbb{E}\left(  \bar{L}_{N_{n}\left(  T,\lambda,1\right)  }%
/n^{1/2}\right) \\
&  =\mathbb{E}\left(  \bar{L}_{N_{n}\left(  T,\lambda,1\right)  }%
/n^{1/2}\right)  .
\end{align*}

We have that
\[
\left\vert \mathbb{E}\left(  \bar{L}_{N_{n}\left(  T,\lambda,1\right)
}/n^{1/2}\right)  \right\vert \leq\left\vert \mathbb{E}\left(  c(1+N_{n}%
\left(  T,\lambda,1\right)  )\right)  \right\vert /n^{1/2}=O\left(
1/n^{1/2}\right)  .
\]
Consequently, we conclude that
\[
n^{1/2}\mathbb{E}\left(  \int_{0}^{T}h_{n}\left(  t\right)  d\bar{D}%
_{1}\left(  t\right)  \right)  =O\left(  1\right)  ,
\]
as $n\rightarrow\infty$, as required. Thus we proved that the first part of
the lemma holds. For the second part, we observe that
\begin{align*}
&  \quad\overline{\lim}_{n\rightarrow\infty}\mathbb{E}\left(  \sup
_{t\in\lbrack0,T]}\left\vert \frac{1}{n^{1/2}}\sum_{i=1}^{n}\left(
\max\left(  t^{2}-T_{i}(n,\lambda)^{\alpha},0\right)  -\mathbb{E}\left[  \max\left(
t^{2}-T_{i}(n,\lambda)^{\alpha},0\right)  \right]  \right)  \right\vert \right) \\
&  =\overline{\lim}_{n\rightarrow\infty}\mathbb{E}\left(  \sup_{t\in
\lbrack0,T]}\left\vert \int_{0}^{t}\frac{1}{n^{1/2}}\sum_{i=1}^{n}\left(
2sI\left(  T_{i}^{\alpha}\left(  n,\lambda\right)  \leq s^{2}\right)  -2s\mathbb{P[}%
T_{i}^{\alpha}\left(  n,\lambda\right)  \leq s^{2}]\right)  ds\right\vert \right) \\
&  \leq\overline{\lim}_{n\rightarrow\infty}\int_{0}^{T}\mathbb{E}\left(
\sup_{t\in\lbrack0,T]}\left\vert \frac{1}{n^{1/2}}\sum_{i=1}^{n}\left(
2tI\left(  T_{i}^{\alpha}\left(  n,\lambda\right)  \leq t^{2}\right)  -2t\mathbb{P[}%
T_{i}^{\alpha}\left(  n,\lambda\right)  \leq t^{2}]\right)  \right\vert \right)  dt\\
&  \leq2T^{2}\overline{\lim}_{n\rightarrow\infty}\mathbb{E}\left(  \sup
_{t\in\lbrack0,T]}\left\vert \frac{1}{n^{1/2}}\sum_{i=1}^{n}\left(  I\left(
T_{i}\left(  n,\lambda\right)  \leq t\right)  -\mathbb{P[}T_{i}\left(  n,\lambda\right)  \leq
t]\right)  \right\vert \right)  <\infty.
\end{align*}
Hence, applying the result for the first part of the lemma, we conclude the
second part as well.
\end{proof}

\begin{proof}
[\textbf{Proof of Lemma \ref{Lem_Tech2}}]%
\begin{align*}
\mathbb{P}\left[  T_{i}\left(  n,\lambda\right)  \leq t\right]   &  =\mathbb{P}\left(
Bin\left(  \mathbb{P}\left(  \left\Vert X_{i}-a\left(  X_{i},\lambda\right)
\right\Vert _{2}\leq t^{1/l}/n^{1/l}\right)  , n-1 \right)  \geq1\right) \\
&  =1-\left(  1-\mathbb{P}\left(  \left\Vert X_{i}-a\left(  X_{i}%
,\lambda\right)  \right\Vert _{2}\leq t^{1/l}/n^{1/l}\right)  \right)  ^{n}.
\end{align*}
Then, as $n\to\infty$ and $t\to0^{+}$
\[
\mathbb{P}\left(  \left\Vert X_{j}-a\left(  X_{i},\lambda\right)  \right\Vert
_{2}\leq t^{1/l}/n^{1/l}\right)  =c_{0}t/n+c_{1}t/n\cdot t^{1/l}%
/n^{1/l}+o\left(  t^{1+1/l}/n^{1+1/l}\right)  .
\]
Therefore by the Poisson approximation to the Binomial distribution we know:
\begin{align*}
\mathbb{P}\left[  T_{i}\left(  n,\lambda\right)  \leq t\right]   &  =1-\exp\left(
-c_{0}t\right)  +O\left(  t^{1+1/l}/n^{1/l}\right)  ,\\
\mathbb{P}\left[  \tau(\lambda)\leq t\right]   &  =1-\exp\left(  -c_{0}t\right)  .
\end{align*}
Thus we proved the first claim:
\[
\mathbb{P}\left[  T_{i}\left(  n,\lambda\right)  \leq t\right]  -\mathbb{P}\left[
\tau(\lambda)\leq t\right]  =O\left(  t^{1+1/l}/n^{1/l}\right)  .
\]
The second claim follows the definition of $\tau$ and equation \eqref{Rate},
where as $n\to\infty$ we have
\begin{align*}
&  \mathbb{P}\left[  T_{i}\left(  n,\lambda\right)  \leq t\right]  -\mathbb{P}\left[
\tau(\lambda)\leq t\right]  = \mathbb{P}\left[  T_{i}\left(  n,\lambda\right)  > t\right]
-\mathbb{P}\left[  \tau (\lambda)> t\right] \\
&  = \mathbb{E}\left[  \exp\left(  -\Lambda\left(  \lambda,X_{1}\right)
\right)  \right]  \left(  1 + O\left(  1/n^{l}\right)  \right)  -
\mathbb{E}\left[  \exp\left(  -\Lambda\left(  \lambda,X_{1}\right)  \right)
\right] \\
&  =\mathbb{P}\left[  \tau(\lambda)> t\right]  O\left(  1/n^{l}\right)  .
\end{align*}

\end{proof}

\subsection{Proofs of Additional Theorems}\label{Sec_additional_proof}

In this subsection, we are going to provide the proofs of the remaining
theorems and corollaries (Theorem \ref{SOSGeneralImplicit}, Theorem
\ref{SOSGeneralExplicit}, Corollary \ref{SOSGeneralImplicitPlugIn} and
Corollary \ref{SOSGeneralExplicitPlugIn}). We are going to follow closely the
proof of Theorem \ref{SOSTheoremMean} and discuss the differences inside each
of its steps.

\subsubsection{\textbf{Proofs of SOS Theorems for General Estimation}}

We will first prove the corresponding theorems for general estimating
equations. As we discussed before, Theorem \ref{SOSGeneralImplicit} is the
direct generalization of Theorem \ref{SOSTheoremMean} and we are going to only
discuss the proof of Theorem \ref{SOSGeneralExplicit} in this part.

\begin{proof}
[\textbf{Proof of Theorem \ref{SOSGeneralExplicit}}]Let us first denote
$\bar{h}_{n}\left(  \theta\right)  =\frac{1}{n}\sum_{i=1}^{n}h\left(
\theta,X_{i}\right)  $. The analogue of \textbf{Step 1}, namely, the dual
formulation takes the form
\[
R_{n}^{W}(\theta_{\ast})=\max_{\lambda}\left\{  -\lambda^{T}\bar{h}_{n}\left(
\theta_{\ast}\right)  -\frac{1}{n}\sum_{i=1}^{n}\max_{j}\left\{  \lambda
^{T}h\left(  \theta_{\ast},Z_{j}\right)  -\lambda^{T}h\left(  \theta_{\ast
},X_{i}\right)  -\left\Vert X_{i}-Z_{j}\right\Vert _{2}^{2}\right\}
^{+}\right\}  .
\]
\textbf{Step 2 and Step 3} are given as follows, for $l=1$ and $l=2$, let us denote
$\sqrt{n}\bar{h}_{n}\left(  \theta_{\ast}\right)  =Z_{n}$ and $\sqrt{n}%
\lambda=2\zeta$, we can scale the SOS\ profile function by $n$, arriving to%
\[
nR_{n}^{W}(\theta_{\ast})=\max_{\zeta}\left\{  -2\zeta^{T}Z_{n}-\frac{1}%
{n}\sum_{i=1}^{n}n\max_{j}\left\{  2\frac{\zeta^{T}}{\sqrt{n}}h\left(
\theta_{\ast},Z_{j}\right)  -2\frac{\zeta^{T}}{\sqrt{n}}h\left(  \theta_{\ast
},X_{i}\right)  -\left\Vert X_{i}-Z_{j}\right\Vert _{2}^{2}\right\}
^{+}\right\}  .
\]
For each $i$, let us consider the maximization problem
\begin{equation}
\max_{j}\left\{  2\frac{\zeta^{T}}{\sqrt{n}}h\left(  \theta_{\ast}%
,Z_{j}\right)  -2\frac{\zeta^{T}}{\sqrt{n}}h\left(  \theta_{\ast}%
,X_{i}\right)  -\left\Vert X_{i}-Z_{j}\right\Vert _{2}^{2}\right\}  .
\label{OP_IN}%
\end{equation}
Similar as Step 1 of the proof for Theorem \ref{SOSTheoremMean}, we would like
to solve the maximization problem (\ref{OP_IN}) by first minimizing over $z$
(as a free variable), instead of over $j$ and then quantify the gap. Observe
that the uniform bound $\left\Vert D_{x}^{2}h\left(  \theta_{\ast}%
,\cdot\right)  \right\Vert <\tilde{K}$ stated in BE1) implies that for all $n$
large enough (in particular, $n^{1/2}>2\tilde{K}\left\Vert \zeta\right\Vert $)
implies that%
\begin{equation}
\max_{z}\left\{  2\frac{\zeta^{T}}{\sqrt{n}}h\left(  \theta_{\ast},z\right)
-2\frac{\zeta^{T}}{\sqrt{n}}h\left(  \theta_{\ast},X_{i}\right)  -\left\Vert
X_{i}-z\right\Vert _{2}^{2}\right\}  , \label{OP_IN_2}%
\end{equation}
has an optimizer in the interior. Therefore, by the differentiability
assumption stated in BE1) we know that any global minimizer, $\bar{a}_{\ast
}\left(  X_{i},\zeta\right)  $, of the problem (\ref{OP_IN_2}) satisfies
\begin{align}
\bar{a}_{\ast}\left(  X_{i},\zeta\right)   &  =X_{i}+D_{x}h\left(
\theta_{\ast},\bar{a}_{\ast}\left(  X_{i},\zeta\right)  \right)  ^{T}%
\cdot\frac{\zeta}{n^{1/2}}\nonumber\\
&  =X_{i}+D_{x}h\left(  \theta_{\ast},X_{i}\right)  ^{T}\cdot\frac{\zeta
}{n^{1/2}}+O\left(  \frac{\left\Vert \zeta\right\Vert _{2}^{2}}{n}\left\Vert
D_{x}h\left(  \theta_{\ast},\bar{a}_{\ast}\left(  X_{i},\zeta\right)  \right)
\right\Vert _{2}\right)  . \label{Eq_Ex_opt}%
\end{align}
Moreover, owing to BE1), we obtain that
\begin{equation}
\left\Vert D_{x}h\left(  \theta_{\ast},\bar{a}_{\ast}\left(  X_{i}%
,\zeta\right)  \right)  -D_{x}h\left(  \theta_{\ast},X_{i}\right)  \right\Vert
_{2}\leq\tilde{K}\frac{\left\Vert \zeta\right\Vert _{2}}{n^{1/2}}.
\label{Eq_Grad_Aux}%
\end{equation}
Consequently, if we define%
\[
a_{\ast}\left(  X_{i},\zeta\right)  =X_{i}+D_{x}h\left(  \theta_{\ast}%
,X_{i}\right)  ^{T}\cdot\frac{\zeta}{n^{1/2}},
\]
we obtain due to (\ref{Eq_Ex_opt}) and (\ref{Eq_Grad_Aux}) that%
\[
\left\Vert a_{\ast}\left(  X_{i},\zeta\right)  -\bar{a}_{\ast}\left(
X_{i},\zeta\right)  \right\Vert _{2}=O\left(  \frac{\left\Vert \zeta
\right\Vert _{2}^{2}}{n}\left(  \left\Vert D_{x}h\left(  \theta_{\ast}%
,X_{i}\right)  \right\Vert _{2}+\frac{\left\Vert \zeta\right\Vert _{2}%
}{n^{1/2}}\right)  \right)  .
\]
Then, after performing a Taylor expansion and applying inequality
(\ref{Eq_Grad_Aux}) we obtain that
\begin{align*}
&  \quad2\frac{\zeta^{T}}{\sqrt{n}}h\left(  \theta_{\ast},X_{i}\right)
-2\frac{\zeta^{T}}{\sqrt{n}}h\left(  \theta_{\ast},\bar{a}_{\ast}\left(
X_{i},\zeta\right)  \right)  +\left\Vert X_{i}-\bar{a}_{\ast}\left(
X_{i},\zeta\right)  \right\Vert _{2}^{2}\\
&  =2\frac{\zeta^{T}}{\sqrt{n}}h\left(  \theta_{\ast},X_{i}\right)
-2\frac{\zeta^{T}}{\sqrt{n}}h\left(  \theta_{\ast},a_{\ast}\left(  X_{i}%
,\zeta\right)  \right)  +\left\Vert X_{i}-a_{\ast}\left(  X_{i},\zeta\right)
\right\Vert _{2}^{2}\\
&  \quad+O\left(  \frac{\left\Vert \zeta\right\Vert ^{3}}{n^{3/2}}\right)
+O\left(  \frac{\left\Vert D_{x}h\left(  \theta_{\ast},X_{i}\right)
\right\Vert _{2}^{2}\left\Vert \zeta\right\Vert _{2}^{3}}{n^{3/2}}\right)  .
\end{align*}
In turn, a direct calculation gives that, as $n\rightarrow\infty$
\begin{align*}
-\frac{\zeta^{T}V_{i}\zeta}{n}  &  =2\frac{\zeta^{T}}{\sqrt{n}}h\left(
\theta_{\ast},X_{i}\right)  -2\frac{\zeta^{T}}{\sqrt{n}}h\left(  \theta_{\ast
},a_{\ast}\left(  X_{i},\zeta\right)  \right)  +\left\Vert X_{i}-a_{\ast
}\left(  X_{i},\zeta\right)  \right\Vert _{2}^{2}\\
&  +O\left(  \frac{\left\Vert D_{x}h\left(  \theta_{\ast},X_{i}\right)
\right\Vert ^{2}\left\Vert \zeta\right\Vert ^{3}}{n^{3/2}}\right)  .
\end{align*}

Similarly as in Step 2 of the proof of Theorem \ref{SOSTheoremMean} we can
define the point process $N^{(i)}\left(  t,\zeta\right)  $ and $T_{i}\left(
n,\lambda\right)  $. We know the gap between freeing the variable $z$ and restricting
the maximization over the $Z_{j}$'s (i.e. the difference between
(\ref{OP_IN_2}) and (\ref{OP_IN}))\ is
\begin{align*}
&  \max_{j}\left\{  \frac{1}{n}\zeta^{T}V_{i}\zeta-\left(  2\frac{\zeta^{T}%
}{\sqrt{n}}h\left(  \theta_{\ast},Z_{j}\right)  -2\frac{\zeta^{T}}{\sqrt{n}%
}h\left(  \theta_{\ast},X_{i}\right)  -\left\Vert X_{i}-Z_{j}\right\Vert
_{2}^{2}\right)  \right\} \\
&  +O\left(  \frac{\left\Vert D_{x}h\left(  \theta_{\ast},X_{i}\right)
\right\Vert ^{2}\left\Vert \zeta\right\Vert ^{3}}{n^{3/2}}\right)  .
\end{align*}
By the definition of $T_{i}\left(  n,\lambda\right)  $, we can write the profile
function for $l=1$ as
\begin{align*}
&  nR_{n}^{W}(\theta_{\ast})\\
&  =\max_{\zeta}\left\{  -2\zeta^{T}Z_{n}-\frac{1}{n}\sum_{i=1}^{n}\max\left(
\zeta^{T}V_{i}\zeta-\frac{T_{i}^{2}\left(  n,\lambda\right)  }{n}+O\left(
\frac{\left\Vert D_{x}h\left(  \theta_{\ast},X_{i}\right)  \right\Vert
^{2}\left\Vert \zeta\right\Vert ^{3}}{n^{1/2}}\right)  ,0\right)  \right\}  .
\end{align*}
Note that the sequence of global optimizers is tight as $n\rightarrow\infty$
because $\mathbb{E}\left(  V_{i}\right)  $ is assumed to be strictly positive
definite with probability one. In turn, from the previous expression we
obtain, following a similar derivation as in the proof of Theorem
\ref{SOSTheoremMean} (invoking Lemma \ref{Lem_Tech}) and using the fact that
$\zeta$ can be restricted to compact sets, that as $n\rightarrow\infty$
\[
nR_{n}^{W}(\theta_{\ast})=\max_{\zeta}\left\{  -2\zeta^{T}Z_{n}-\mathbb{E}%
\left[  \max\left(  \zeta^{T}V_{1}\zeta-\frac{T_{1}^{2}\left(  n,\lambda\right)  }%
{n}\right)  \right]  \right\}  +o_{p}\left(  1\right)  .
\]
Then, for $l=2$, as $n\rightarrow\infty$ we have estimate for the profile
function as
\[
nR_{n}^{W}(\theta_{\ast})=\max_{\zeta}\left\{  -2\zeta^{T}Z_{n}-\mathbb{E}%
\left[  \max\left(  \zeta^{T}V_{1}\zeta-T_{1}^{2}\left(  n,\lambda\right)  \right)
\right]  \right\}  +o_{p}\left(  1\right)  .
\]

When $l\geq3$, let us denote $\sqrt{n}\bar{h}_{n}\left(  \theta_{\ast}\right)
=Z_{n}$ and $n^{\frac{3}{2l+2}}\lambda=2\zeta$, we can scale profile function
by $n^{\frac{1}{2}+\frac{3}{2l+2}}$ and write it as
\begin{align*}
&  n^{\frac{1}{2}+\frac{3}{2l+2}}R_{n}^{W}(\theta_{\ast})\\
&  =\max_{\zeta}\left\{  -2\zeta^{T}Z_{n}-\frac{1}{n}\sum_{i=1}^{n}n^{\frac
{1}{2}+\frac{3}{2l+2}}\max_{j}\left\{  2\frac{\zeta^{T}}{n^{\frac{3}{2l+2}}%
}h\left(  \theta_{\ast},Z_{j}\right)  -2\frac{\zeta^{T}}{n^{\frac{3}{2l+2}}%
}h\left(  \theta_{\ast},X_{i}\right)  -\left\Vert X_{i}-Z_{j}\right\Vert
_{2}^{2}\right\}  ^{+}\right\}  .
\end{align*}
By applying same derivation as for $l=1\text{ and }2$ above, we can define a
point process $N^{(i)}\left(  t,\zeta\right)  $ and $T_{i}\left(  n\right)  $
as in the proof of Theorem \ref{SOSTheoremMean}. As $n\to\infty$, we have the
estimate for profile function becomes
\begin{align*}
&  n^{\frac{1}{2}+\frac{3}{2l+2}}R_{n}^{W}(\theta_{\ast})\\
&  =\max_{\zeta}\left\{  -2\zeta^{T}Z_{n}-n^{\frac{1}{2}+\frac{3}{2l+2}%
-\frac{2}{l}}\frac{1}{n}\sum_{i=1}^{n}\max\left(  n^{-\left(  \frac{6}%
{2l+2}-\frac{2}{l}\right)  }\zeta^{T}V_{i}\zeta-T_{i}^{2/l}\left(  n,\lambda\right)
,0\right)  \right\}  +o_{p}\left(  1\right) \\
&  =\max_{\zeta}\left\{  -2\zeta^{T}Z_{n}-n^{\frac{1}{2}+\frac{3}{2l+2}%
-\frac{2}{l}}\mathbb{E}\left[  \max\left(  n^{-\left(  \frac{6}{2l+2}-\frac
{2}{l}\right)  }\zeta^{T}V_{1}\zeta-T_{1}^{2/l}\left(  n,\lambda\right)  ,0\right)
\right\}  \right]  +o_{p}\left(  1\right)  .
\end{align*}
The final estimation follows as in the proof for Theorem \ref{SOSTheoremMean}
(i.e. applying Lemma \ref{Lem_Tech}).

In \textbf{Step 4} for $l=1$, as $n\to\infty$ the objective function is
\[
nR_{n}^{W}(\theta_{\ast})=\max_{\zeta}\left\{  -2\zeta^{T}Z_{n}\left(
\theta_{\ast}\right)  -\mathbb{E}\left[  \max\left(  \zeta^{T}V_{1}\zeta
-\frac{T_{1}^{2}\left(  n,\lambda\right)  }{n},0\right)  \right]  \right\}
+o_{p}(1).
\]
Let us denote $G:\mathbb{R}^{l}\to\mathbb{R}$ to be a deterministic continuous
function, defined as
\[
G\left(  \zeta,n\right)  = \mathbb{E}\left[  \max\left(  \zeta^{T}V_{1}%
\zeta-\frac{T_{1}^{2}\left(  n,\lambda\right)  }{n},0\right)  \right]  .
\]
We know $\Upsilon=\mathbb{E}\left[  V_{1}\right]  $ is symmetric strictly
positive definite matrix, then the objective function is strictly convex and
differentiable in $\zeta$. Thus the (unique) global maximizer is characterized
by the first order optimality condition almost surely. We take derivative
w.r.t. $\zeta$ and set it to be $0$, applying the same estimation in the
original proof the first order optimality condition becomes
\begin{equation}
Z_{n}=-\Upsilon\zeta+o_{p}(1) \text{ as }n\to\infty. \label{l1FirstOrder}%
\end{equation}
Since $\Upsilon$ is invertible, for any $n$ we can solve optimal $\zeta
_{n}^{\ast}=-\Upsilon^{-1}Z_{n}+o_{p}(1)$. Plugging $\zeta_{n}^{\ast}$ in the
objective function, as $n\to\infty$ we have
\[
nR_{n}^{W}(\theta_{\ast})=2{Z}_{n}^{T}\Upsilon^{-1}{Z}_{n}- G\left(
-\Upsilon^{-1}Z_{n},n\right)  +o_{p}(1).
\]
As $n\rightarrow\infty$, we can apply the same estimation in the proof of
Theorem \ref{SOSTheoremMean}, it becomes
\[
nR_{n}^{W}(\theta_{\ast})\Rightarrow\tilde{Z}^{T}\Upsilon^{-1}\tilde{Z}.
\]
Thus we proof the claim for $l=1$.

In \textbf{Step 5} for $l=2$, as $n\to\infty$ the objective function has
estimate
\[
nR_{n}^{W}(\theta_{\ast})=\max_{\zeta}\left\{  -2\zeta^{T}Z_{n}\left(
\theta_{\ast}\right)  -\mathbb{E}\left[  \max\left(  \zeta^{T}V_{1}%
\zeta-{T_{1}\left(  n,\lambda\right)  },0\right)  \right]  \right\}  +o_{p}(1).
\]
Still, we denote $G\left(  \zeta,n\right)  $ to be a deterministic function
given as,
\[
G\left(  \zeta,n\right)  = \mathbb{E}\left[  \max\left(  \zeta^{T}V_{1}
\zeta-{T_{1}\left(  n,\lambda\right)  },0\right)  \right]  .
\]
Same as discussed in for $l=1$, the objective function is strictly convex and
differentiable in $\zeta$, thus the (unique) global maximizer could be
characterized via first order optimality condition almost surely. We take
derivative w.r.t. $\zeta$ and set it to be $0$, applying same estimation in
the proof of Theorem \ref{SOSTheoremMean} the first order optimality condition
becomes
\begin{equation}
Z_{n}=-\mathbb{E}\left[  V_{1}1_{\left(  \tau(0)\leq\zeta^{T}V_{1}\zeta\right)
}\right]  \zeta+o_{p}(1) \text{ as }n\to\infty. \label{l2FirstOrder}%
\end{equation}
We know the objective function is strictly convex differentiable, then for
fixed $Z_{n}$ there is a unique $\zeta_{n}^{\ast}$ that satisfies the first
order optimality condition (\ref{l2FirstOrder}). We plug in the optimizer and
the objective function becomes
\[
nR_{n}^{W}\left(  \theta_{\ast}\right)  =-2Z_{n}^{T}\zeta_{n}^{\ast}- G\left(
\zeta_{n}^{\ast},n\right)  +o_{p}(1)\text{ as }n\to\infty.
\]
As $n\rightarrow\infty$, we can apply the same estimation in the proof of
Theorem \ref{SOSTheoremMean}, we have
\[
nR_{n}^{W}\left(  \theta_{\ast}\right)  \Rightarrow-2\tilde{Z}^{T}\tilde
{\zeta}-\tilde{\zeta}^{T}\tilde{G}\left(  \tilde{\zeta}\right)  \tilde{\zeta
},
\]
where $\tilde{G}:\mathbb{R}^{q}\to\mathbb{R}^{q}\times\mathbb{R}^{q}$ is a
deterministic continuous mapping defined as,
\[
\tilde{G}\left(  \zeta\right)  = \mathbb{E}\left[  V_{1}\max\left(  1-{\tau(0)
}/{\left(  \zeta^{T}V_{1}\zeta\right)  },0\right)  \right]  ,
\]
and $\tilde{\zeta}:=\tilde{\zeta}\left(  \tilde{Z}\right)  $ is the unique
solution to
\[
\tilde{Z}=-\zeta\mathbb{E}\left[  V_{1}1_{\left(  \tau(0)\leq\zeta^{T}V_{1}%
\zeta\right)  }\right]  .
\]
Then we proved the claim for $l=2$.

Finally, in \textbf{Step 6} for $l\geq3$, as $n\to\infty$ the objective
function is
\begin{align*}
&  n^{1/2+\frac{3}{2l+2}}R_{n}^{W}(\theta_{\ast})\\
&  =\max_{\zeta}\left\{  -2\zeta^{T}Z_{n}-n^{\left(  1/2+\frac{3}{2l+2}%
-\frac{2}{l}\right)  }\mathbb{E}\left[  \max\left(  n^{-\left(  \frac{6}%
{2l+2}-\frac{2}{l}\right)  }\zeta^{T}V_{1}\zeta-T_{1}^{2/l}\left(  n,\lambda\right)
,0\right)  \right]  \right\}  +o_{p}(1).
\end{align*}
We denote $G\left(  \zeta,n\right)  $ to be a deterministic function defined
as,
\[
G\left(  \zeta,n\right)  = n^{\left(  1/2+\frac{3}{2l+2} -\frac{2}{l}\right)
}\mathbb{E}\left[  \max\left(  n^{-\left(  \frac{6} {2l+2}-\frac{2}{l}\right)
}\zeta^{T}V_{1}\zeta-T_{1}^{2/l}\left(  n,\lambda\right)  ,0\right)  \right]  .
\]
Follows the same discussion above for $l=1\text{ and }2$, we know the
objective function is strictly convex differentiable in $\zeta$ and the global
maximizer is characterized by first order optimality condition almost surely.
We take derivative of the objective function w.r.t. $\zeta$ and set it to be
$0$. We apply the same technique as in the proof of Theorem
\ref{SOSTheoremMean}, the first order optimality condition becomes
\begin{equation}
Z_{n}=-\mathbb{E}\left[  V_{1}\frac{\pi^{l/2}\left(  f_{X}\left(
X_{1}\right)  +\kappa f_{Y}\left(  X_{1}\right)  \right)  }{\Gamma
(l/2+1)}V_{1}\left(  \zeta^{T}V_{1}\zeta\right)  ^{l}\right]  \zeta
+o_{p}(1).\text{ as }n\to\infty\label{l3FirstOrder}%
\end{equation}
The objective condition is strictly convex differentiable and for fixed
$Z_{n}$ there is a unique $\zeta_{n}^{\ast}$ satisfying the first optimality
condition (\ref{l3FirstOrder}). We plug $\zeta_{n}^{\ast}$ into the objective
function and it becomes
\[
n^{1/2+\frac{3}{2l+2}}R_{n}^{W}(\theta_{\ast})=-2Z_{n}^{T}\zeta_{n}^{\ast
}-G\left(  \zeta_{n}^{\ast},n\right)  +o_{p}(1) \text{ as }n\to\infty.
\]
As $n\rightarrow\infty$, we can apply same estimate in the proof of Theorem
\ref{SOSTheoremMean}, we have
\[
n^{1/2+\frac{3}{2l+2}}R_{n}^{W}(\theta_{\ast})\Rightarrow-2\tilde{Z}^{T}
\tilde{\zeta}-\frac{2}{l+2}\tilde{G}\left(  \tilde{\zeta}\right)  ,
\]
where $\tilde{G}:\mathbb{R}^{q}\to\mathbb{R}$ is a deterministic continuous
function given as,
\[
\tilde{G}\left(  \zeta\right)  = \mathbb{E}\left[  \frac{\pi^{l/2}\left(
f_{X}(X_{1})+\kappa f_{Y}\left(  X_{1}\right)  \right)  }{\Gamma\left(
l/2+1\right)  }\left(  \zeta^{T}V_{1}\zeta\right)  ^{l/2+1}\right]  ,
\]
and $\tilde{\zeta}:=\tilde{\zeta}\left(  \tilde{Z}\right)  $ is the unique
solution to
\[
\tilde{Z}=-\mathbb{E}\left[  V_{1}\frac{\pi^{l/2}\left(  f_{X}\left(
X_{1}\right)  +\kappa f_{Y}\left(  X_{1}\right)  \right)  }{\Gamma
(l/2+1)}V_{1}\left(  \zeta^{T}V_{1}\zeta\right)  ^{l}\right]  \zeta.
\]
We proved the claim for $l\geq3$ and finish the proof for Theorem
\ref{SOSGeneralExplicit}.
\end{proof}

\subsubsection{\textbf{Proofs of SOS Theorems for General Estimation with
Plug-In}}

The proofs of the plug-in version of SOS theorems for general estimation
equation also mainly follows the proof of Theorem \ref{SOSTheoremMean}, we are
going to discuss the different steps here.

\begin{proof}
[\textbf{Proof of Corollary \ref{SOSGeneralImplicitPlugIn}}]For implicit
formulation, as we discussed for Theorem \ref{SOSGeneralImplicit}, we can
redefine $X_{i}\leftarrow h\left(  \gamma_{\ast},\nu_{n},X_{i}\right)  $,
$Z_{k}\leftarrow h\left(  \gamma_{\ast},\nu_{n},Z_{k}\right)  $, $X_{i}%
(\ast)\leftarrow h\left(  \gamma_{\ast},\nu_{\ast},X_{i}\right)  $ and
$Z_{k}(\ast)\leftarrow h\left(  \gamma_{\ast},\nu_{\ast},X_{i}\right)  $. Then
the proof for the implicit formulation with plug-in goes as follows.

In \textbf{Step 1}, the dual formulation is similar given as
\begin{align*}
R_{n}^{W}(\gamma_{\ast})  &  =\max_{\lambda,\gamma_{i}\geq0}\left\{
-\lambda\bar{X}_{n}-\frac{1}{n}\sum_{i=1}^{n}\gamma_{i}\right\} \\
\text{s.t.}\;-\gamma_{i}  &  \leq\min_{j}\left\{  \lambda^{T}X_{i}-\lambda
^{T}Z_{j}+\left\Vert X_{i}-Z_{j}\right\Vert _{2}^{2}\right\}  ,\text{ for all
}i\text{.}%
\end{align*}
We can apply first order Taylor expansion to $h\left(  \gamma_{\ast},\nu
_{n},X_{i}\right)  $ w.r.t. $\nu$, then we have
\[
h\left(  \gamma_{\ast},\nu_{n},X_{i}\right)  =h\left(  \gamma_{\ast},\nu
_{\ast},X_{i}\right)  +O_{p}\left(  \frac{\left\Vert D_{\nu}h\left(
\gamma_{\ast},\bar{\nu}_{n},X_{i}\right)  \right\Vert }{n^{1/2}}\right)  ,
\]
where $\bar{\nu}_{n}$ is a point between $\nu_{n}$ and $\nu_{\ast}$. By our
change of notation for $X_{i}$, $X_{i}(\ast)$, $Z_{k}$ and $Z_{k}(\ast)$ and
the above Taylor expansion, we can observe
\[
Z_{k}=Z_{k}(\ast)+\epsilon_{n}\left(  Z_{k}\right)  ,
\]
where $\epsilon_{n}\left(  Z_{k}\right)  =O_{p}\left(  \left\Vert D_{\nu
}h\left(  \gamma_{\ast},\bar{\nu}_{n},Z_{k}\right)  \right\Vert /n^{1/2}%
\right)  $.

In \textbf{Step 2} we can define a point process $N_{n}^{(i)}\left(
t,\lambda\right)  $ and $T_{i}\left(  n\right)  $ as in the proof of Theorem
\ref{SOSTheoremMean}, but the rate becomes
\[
\Lambda\left(  X_{i},\lambda\right)  =\left[  f_{X}\left(  X_{i}%
+\lambda/2+\epsilon_{n}\left(  X_{i}\right)  \right)  +\kappa f_{Y}\left(
X_{i}+\lambda/2+\epsilon_{n}\left(  X_{i}\right)  \right)  \right]  \frac
{\pi^{l/2}}{\Gamma\left(  l/2+1\right)  }.
\]
As $n\rightarrow\infty$, same as in the proof of Theorem \ref{SOSTheoremMean}
and Theorem \ref{SOSGeneralExplicit} we can argue $\lambda\rightarrow0$. Then
we can define $\tau(0)$ same as in the proof of Theorem \ref{SOSTheoremMean} and
has the with same distribution
\[
\mathbb{P}\left[  \tau(0)\geq t\right]  =\mathbb{E}\left[  \exp\left(  -\left(
f_{X}\left(  X_{1}\right)  +\kappa f_{Y}\left(  X_{1}\right)  \right)
\frac{\pi^{l/2}}{\Gamma\left(  l/2+1\right)  }\right)  \right]  .
\]
Then the rest of the proof in \textbf{Step 3, 4, 5 and 6} stay the same as
that of Theorem \ref{SOSTheoremMean}, but replacing the CLT for $Z_{n}$ by the
asymptotic distribution given in C2).
\end{proof}

\begin{proof}
[\textbf{Proof of Corollary \ref{SOSGeneralExplicitPlugIn}}]For explicit
formulation, the proof follows more closely the proof of Theorem
\ref{SOSGeneralExplicit} and we are discussing the differences as follows.

In \textbf{Step 1}, the dual formulation takes the form
\begin{align*}
&  R_{n}^{W}(\theta_{\ast})\\
&  =\max_{\lambda}\left\{  -\lambda^{T}\bar{h}_{n}\left(  \gamma_{\ast}%
,\nu_{n}\right)  -\frac{1}{n}\sum_{i=1}^{n}\max_{j}\left\{  \lambda
^{T}h\left(  \gamma_{\ast},\nu_{n},Z_{j}\right)  -\lambda^{T}h\left(
\gamma_{\ast},\nu_{n},X_{i}\right)  -\left\Vert X_{i}-Z_{j}\right\Vert
_{2}^{2}\right\}  ^{+}\right\}  .
\end{align*}
\textbf{Step 2 and Step 3} Follows the same as for the proof of Theorem
\ref{SOSGeneralExplicit} however we need to notice that difference is the
definition of $\bar{a}_{\ast}\left(  X_{i},\zeta\right)  $, for $l=1\text{ and
}2$ we have
\begin{align}
\bar{a}_{\ast}\left(  X_{i},\zeta\right)   &  =X_{i}+D_{x}h\left(
\gamma_{\ast},\nu_{n},\bar{a}_{\ast}\left(  X_{i},\zeta\right)  \right)
\cdot\frac{\zeta}{n^{1/2}}\\
&  =X_{i}+D_{x}h\left(  \gamma_{\ast},\nu_{n},X_{i}\right)  \cdot\frac{\zeta
}{n^{1/2}}+O\left(  \frac{\left\Vert \zeta\right\Vert _{2}^{2}}{n}\left\Vert
D_{x}h\left(  \gamma_{\ast},\nu_{n},\bar{a}_{\ast}\left(  X_{i},\zeta\right)
\right)  \right\Vert _{2}\right) \nonumber\\
&  =X_{i}+D_{x}h\left(  \gamma_{\ast},\nu_{\ast},X_{i}\right)  \cdot
\frac{\zeta}{n^{1/2}}+O\left(  \frac{\left\Vert \zeta\right\Vert _{2}^{2}}%
{n}\left\Vert D_{x}h\left(  \gamma_{\ast},\nu_{n},\bar{a}_{\ast}\left(
X_{i},\zeta\right)  \right)  \right\Vert _{2}\right) \nonumber\\
&  +O\left(  \frac{\left\Vert \zeta\right\Vert _{2}}{n^{1/2}}\left\Vert
\nu_{n}-\nu_{\ast}\right\Vert _{2}\left\Vert D_{x}h\left(  \gamma_{\ast}%
,\nu_{n},\bar{a}_{\ast}\left(  X_{i},\zeta\right)  \right)  \right\Vert
_{2}\left\Vert D_{\nu}D_{x}h\left(  \gamma_{\ast},\bar{\nu}_{n},\bar{a}_{\ast
}\left(  X_{i},\zeta\right)  \right)  \right\Vert _{2}\right)  ,
\end{align}
where $\bar{\nu}_{n}$ is a point between $\nu_{n}$ and $\nu_{\ast}$. By
assumption C5)-C7) we can notice the rest of step 2 and 3 stay the same as in
the proof of Theorem \ref{SOSGeneralExplicit}. In \textbf{Step 4, 5 and 6} we
use $Z_{n}=\frac{1}{n^{1/2}}\sum_{i=1}^{n}h\left(  \gamma_{\ast},\nu_{n}%
,X_{i}\right)  \Rightarrow\tilde{Z}^{\prime}$ given in C2).
\end{proof}

\subparagraph{ACKNOWLEDGMENT}

{The authors are grateful to the editorial team and the referees for their careful review of the paper and their constructive suggestions, which greatly helped to improve it. J. Blanchet gratefully acknowledges support from the following NSF grants  
1915967, 1820942, 1838576. The authors are also grateful to Karthyek Murthy and Henry Lam
for helpful discussions.}

\bibliographystyle{plain}
\bibliography{SoS_2}

\end{document}